\documentclass[nohyperref]{article}

\usepackage{microtype}
\usepackage{graphicx}
\usepackage{subfigure}
\usepackage{booktabs} % for professional tables

\usepackage{hyperref}

\usepackage[accepted]{icml2023}

\usepackage{amsmath}
\usepackage{amssymb}
\usepackage{mathtools}
\usepackage{amsthm}

\usepackage[capitalize,noabbrev]{cleveref}

\theoremstyle{plain}
\newtheorem{theorem}{Theorem}

\newtheorem{lemma}[theorem]{Lemma}
\newtheorem{corollary}[theorem]{Corollary}
\theoremstyle{definition}
\newtheorem{definition}[theorem]{Definition}

\theoremstyle{remark}
\newtheorem{remark}[theorem]{Remark}

\usepackage[textsize=tiny]{todonotes}
\usepackage{graphicx}

\newcommand{\cut}[1]{{}}

\makeatletter
\newcommand*{\rom}[1]{\expandafter\@slowromancap\romannumeral #1@}
\makeatother

\newcommand{\vc}{{\mathbf{c}}}

\newcommand{\ve}{{\mathbf{e}}}

\newcommand{\vu}{{\mathbf{u}}}

\newcommand{\vw}{{\mathbf{w}}}
\newcommand{\vx}{{\mathbf{x}}}
\newcommand{\vy}{{\mathbf{y}}}
\newcommand{\vz}{{\mathbf{z}}}

\newcommand{\cH}{{\mathcal{H}}}

\newcommand{\cO}{{\mathcal{O}}}

\newcommand{\cR}{{\mathcal{R}}}

 \usepackage[bb=boondox]{mathalfa}

\usepackage{xparse}
\DeclareFontFamily{U}{ntxmia}{}
\DeclareFontShape{U}{ntxmia}{m}{it}{<-> ntxmia }{}
\DeclareFontShape{U}{ntxmia}{b}{it}{<-> ntxbmia }{}
\DeclareSymbolFont{lettersA}{U}{ntxmia}{m}{it}
\SetSymbolFont{lettersA}{bold}{U}{ntxmia}{b}{it}

\ExplSyntaxOn
\NewDocumentCommand{\varmathbb}{m}
 {
  \tl_map_inline:nn { #1 }
  {
    \use:c { varbb##1 }
  }
 }
\tl_map_inline:nn { ABCDEFGHIJKLMNOPQRSTUVWXYZ }
 {
  \exp_args:Nc \DeclareMathSymbol{varbb#1}{\mathord}{lettersA}{\int_eval:n { `#1+67 }}
 }
\exp_args:Nc \DeclareMathSymbol{varbbk}{\mathord}{lettersA}{169}
\ExplSyntaxOff
\newcommand{\RR}{\mathbb{R}}

\renewcommand{\SS}{{\mathbb{S}}}

\newcommand{\supp}{{\mathrm{supp}}} % support
\newcommand{\tr}{{\mathrm{tr}}} % trace
\newcommand{\Proj}{\Pi}

\DeclareMathOperator*{\argmin}{argmin}

\newcommand*{\fix}{\mathrm{Fix}\,}

\newcommand{\prox}{\mathrm{Prox}}

\newcommand{\dom}{\mathrm{dom}\,} % domain
\newcommand{\reals}{\mathbb{R}}

\newcommand{\bigO}{\mathcal{O}}

\newcommand{\ones}{\mathbf{1}}

\definecolor{lightgrey}{gray}{0.8}
\definecolor{medgrey}{gray}{0.6}
\definecolor{darkgrey}{gray}{0.4}
\newcommand{\opF}{{\varmathbb{F}}}

\newcommand{\opH}{{\varmathbb{H}}}
\newcommand{\opI}{{\varmathbb{I}}}

\newcommand{\opS}{{\varmathbb{S}}}
\newcommand{\opT}{{\varmathbb{T}}}

\newcommand{\lagrange}{\mathbf{L}}

\newcommand{\Span}{\mathrm{span}}

\newcommand{\al}{\opF}
\newcommand{\conv}{\mathrm{conv}\,}
\newcommand{\clconv}{\overline{\mathrm{conv}}\,}

\icmltitlerunning{Accelerated Infeasibility Detection of Constrained Optimization and Fixed-Point Iterations}

\begin{document}

\twocolumn[
\icmltitle{Accelerated Infeasibility Detection of \\
           Constrained Optimization and Fixed-Point Iterations}



\icmlsetsymbol{equal}{*}

\begin{icmlauthorlist}
\icmlauthor{Jisun Park}{snumath}
\icmlauthor{Ernest K. Ryu}{snumath}
\end{icmlauthorlist}

\icmlaffiliation{snumath}{Department of Mathematical sciences, Seoul National University, Seoul, Republic of Korea}

\icmlcorrespondingauthor{Ernest K. Ryu}{ernestryu@snu.ac.kr}

\icmlkeywords{monotone operator,fixed-point iteration,splitting methods,accelerated proximal point method,Halpern iteration,acceleration,complexity lower bound,anchor acceleration,convex optimization,infeasibility detection,PG-EXTRA,decentralized optimization}

\vskip 0.3in
]



\printAffiliationsAndNotice{}  % leave blank if no need to mention equal contribution

\begin{abstract}
As first-order optimization methods become the method of choice for solving large-scale optimization problems, optimization solvers based on first-order algorithms are being built. Such general-purpose solvers must robustly detect infeasible or misspecified problem instances, but the computational complexity of first-order methods for doing so has yet to be formally studied. In this work, we characterize the optimal accelerated rate of infeasibility detection. We show that the standard fixed-point iteration achieves a $\mathcal{O}(1/k^2)$ and $\mathcal{O}(1/k)$ rates, respectively, on the normalized iterates and the fixed-point residual converging to the infimal displacement vector, while the accelerated fixed-point iteration achieves $\mathcal{O}(1/k^2)$ and $\tilde{\mathcal{O}}(1/k^2)$ rates. We then provide a matching complexity lower bound to establish that $\Theta(1/k^2)$ is indeed the optimal accelerated rate.
\end{abstract}

\section{Introduction}
\label{sec:intro}

First-order optimization methods have become the method of choice for solving the large-scale optimization problems of the modern era. As first-order methods scale more favorably than classical interior-point methods \cite{ODonoghue2016SCS,Stellato2020OSQP,garstka2021cosmo}, new optimization solvers based on first-order algorithms are being built with the goal of replacing classical solvers based on interior-point methods or simplex methods in large-scale applications.

However, these new first-order solvers are far less equipped to robustly detect infeasible or misspecified problem instances. A general-purpose solver must robustly detect infeasible problem instances arising from user misspecification or from applications such as embedded application, mixed-integer optimization with branch-and-bound technique, or combinatorial optimization \cite{naik2017embedded,de2012computation}. 
Classical solvers based on interior point methods or simplex methods, in their first phase, determines whether the problem is feasible or infeasible by finding a feasible point.
The behavior of such classical solvers under pathologies is well understood through extensive theoretical research and through the decades-long deployment of open-source and commercial solvers. The analysis of first-order algorithms such as Douglas-Rachford splitting (DRS) and ADMM applied to pathological problem instances has started to gain attention. However, the computational complexity of determining the infeasibility of a given problem instance has yet to be formally studied.

In this work, we characterize the optimal accelerated rate of infeasibility detection by analyzing the convergence rates of fixed-point iterations towards the infimal displacement vector, which serves as a certificate of infeasibility. We show that the standard fixed-point iteration achieves a $\mathcal{O}(1/k^2)$ and $\mathcal{O}(1/k)$ rates, respectively, on the normalized iterates and the fixed-point residual, while the accelerated fixed-point iteration achieves $\mathcal{O}(1/k^2)$ and $\tilde{\mathcal{O}}(1/k^2)$ rates. We then provide a matching complexity lower bound to establish that $\Theta(1/k^2)$ is indeed the optimal accelerated rate.

\subsection{Preliminaries and notations}
\label{sec:preliminary}

We use the standard notations in \citet{ryu2022large}.

\paragraph{Sets and operators.}
Let $\cH$ be a real Hilbert space.
For a set $C\subseteq\cH$, we denote by $\conv C$ a convex hull of $C$, $\overline{C}$ a closure of $C$, and $\clconv C$ a closure of a convex hull of $C$.
If $C$ is a nonempty closed convex set, for any $x\in\cH$, there exists a unique vector $z\in C$ such that $z \in \argmin_{y\in C} \|x-y\|^2$, which is denoted as $\Proj_C(x)$ and called a projection of $x$ onto $C$.
We also let $\delta_C$ refer to the indicator function of set $C$.
We denote by $\SS^n$ the set of all $n\times n$ symmetric matrices and by $\SS^n_+$ the set of all $n\times n$ symmetric positive semidefinite matrices.
We say $X \succeq Y$ if $X-Y\in\SS^n_+$ for $X,Y\in\SS^n$.

Let $\opT\colon \cH\rightrightarrows\cH$ be a set-valued operator.
$\dom\opT = \{x \mid \opT x \neq \emptyset\}$ is called the {domain} of $\opT$,
and $\cR(\opT)=\{y \mid \exists x \in \cH \textit{ s.t. } y \in \opT x\}$ is called the {range} of $\opT$.
For a single-valued operator $\opT\colon\cH\to\cH$, it is called {nonexpansive} if $\|\opT x - \opT y\| \le \|x-y\|$ for all $x,y\in\dom\opT$ and 
$\gamma$-{contractive} with $0<\gamma<1$ if $\|\opT x - \opT y\| \le \gamma \|x-y\|$ holds for all $x,y\in\dom\opT$, and {$\theta$-averaged} if there exists nonexpansive operator $\opS$ and identity operator $\opI$ such that $\opT = \theta \opS + (1-\theta) \opI$.
$\opT$ is called {maximal} nonexpansive (contractive) if $\dom\opT = \cH$.
If $x_\star\in\cH$ is a point such that $x_\star=\opT x_\star$, we call $x_\star$ a {fixed point} of $\opT$.
$\fix\opT\subseteq\cH$ denotes a set of fixed points of $\opT$.

\paragraph{Fixed-point iteration.}
Classical Banach fixed-point theorem illustrates that if $\opT\colon\cH\to\cH$ is a contraction, then $\fix\opT$ is nonempty and the \emph{Picard iteration} \eqref{eq:fpi}
\begin{equation}
    x^{k+1} = \opT x^k, \qquad k=0,1,\ldots
    \tag{Picard}
    \label{eq:fpi}
\end{equation}
starting from $x^0\in\cH$ converges to some $x_\star\in\fix\opT$.
When $\opT$ is nonexpansive but not necessarily contractive, \eqref{eq:fpi} may not converge to the fixed point of $\opT$.
In such cases, to guarantee the convergence, one may use \emph{Krasnosel'ski{\u\i}-Mann iteration} \cite{krasnosel1955two,mann1953mean} or \emph{Halpern iteration} \cite{halpern1967fixed}, whose forms are described in \cref{sec:convergencerate}.

\paragraph{Constrained optimization and fixed-point iterations.}
Consider a constrained optimization problem
\[
    \begin{array}{ll}
        \underset{x \in \cH}{\mbox{minimize}} & f(x)\\
        {\mbox{subject to}}&x\in C,
    \end{array}
\]
where $f\colon\reals^n\to\reals$ is convex and $C \subseteq \reals^n$ is a nonempty closed convex set.
Problems of this type can be solved with various first-order methods including projected gradient method, proximal gradient method, alternating direction method of multipliers (ADMM), and primal-dual hybrid gradient (PDHG).
These methods can be understood and analyzed as nonexpansive fixed-point iterations \cite{ryu2022large}.
Therefore, the analysis of fixed-point iteration broadly applies to this broad class of first-order methods.
% For instance, ADMM is a convergent algorithm, and it can be proven from the fact that ADMM is a convergent fixed-point iteration.

\paragraph{Inconsistent operators.}
We say $\opT$ is {consistent} if $\fix\opT\neq\emptyset$ and {inconsistent} if $\fix\opT=\emptyset$.
$\opT$ is inconsistent if and only if $0 \notin \cR(\opI-\opT)$.
From the well-known fact that $\overline{\cR(\opI-\opT)}$ is closed and convex \citep[Lemma 4]{Pazy1971},
$v = \Proj_{\overline{\cR(\opI-\opT)}}(0)$ is well-defined, and $v$ is called {infimal displacement vector}.
If $\opT$ is consistent, then $v=0$.
If $v\neq 0$, then $\opT$ is inconsistent.

Any nonexpansive operator $\opT$ is of exactly one of these three cases: (i) $\fix\opT\neq\emptyset$, (ii) $v\neq0$, and (iii) $\fix\opT=\emptyset$ with $v=0$.
In convex optimization, (i) corresponds to the case where primal and dual solution exist and primal-dual gap being $0$, and (ii) corresponds to the case where either the primal problem or the dual problem is infeasible.
(iii) corresponds to pathological weakly feasible and weakly infeasible cases \cite{BanjacGoulart2019,LiuRyuYin2019,RyuLinYin2019douglas}.
Our main focus will be on cases (i) and (ii).

\subsection{Prior work}
\label{sec:priorworks}

\paragraph{Inconsistent fixed-point iteration.}
\citet{Browder1966} first proved that the iterates of Picard iteration is bounded if and only if $\opT$ is consistent, and later followed by work of \citet{Pazy1971} showing the convergence of $-x^k/k$ to the infimal displacement vector.
This result has been extended to Banach space setup by \citet{reich1973asymptotic}.
If $\opT$ is more than just a nonexpansive operator, then difference of iterates $\opT^k x^0 - \opT^{k+1} x^0$ also converges to the infimal displacement vector; see \citet{bailion1978asymptotic}, \citet{ReichSharfrir1987}, and \citet{bruck1977weak} for averaged, firmly-nonexpansive, and strongly nonexpansive operators in Banach spaces.
For more on general settings, see \citet{Reich1981asymptotic,Reich1982asymptotic,PlantReich1983asymptotics,ariza2014firmly,Nicolae2013asymptotic}.
Despite its numerous appearance, it was not until in late 1990s where the term `minimal displacement vector' was coined \cite{BauschkeBorweinLewis1997}.
It was later called `infimal displacement vector' \cite{BauschkeHareMoursi2014}, and its properties have been analyzed with depth as well \cite{BauschkeDouglasMoursi2016,Ryu2018,BauschkeMoursi2018magnitude,BauschkeMoursi2020minimal}.

\paragraph{First-order numerical solvers.}
The interior point method \cite{nesterov1994interior} has been successful in solving convex optimization problems, and a number of numerical solvers based on this exists \cite{nesterov1994interior,sturm1999using,gurobi,aps2019mosek,mattingley2012cvxgen}.
Recently, first-order method solving conic optimization programs has gained huge interest, due to its scalability to very large and high-dimensional problems.
ADMM-based solvers such as SCS \cite{ODonoghue2016SCS,sopasakis2019superscs}, OSQP \cite{Stellato2020OSQP}, and COSMO \cite{garstka2021cosmo}, and also include PDHG-based solver PDLP \cite{chambolle2011first,Applegate2021PDLP} are first-order numerical solvers.

\paragraph{Constrained optimization and infeasibility.}
For convex feasibility problem, primary choice of methods are cyclic projection \cite{von1951functional}, Dykstra's algorithm \cite{dykstra1983algorithm}, AAR method \cite{bauschke2004finding}, and so on.
These methods have been analyzed extensively \cite{boyle1986method,bauschke1994dykstra,BauschkeBorweinLewis1997,ArtachoBorweinTam2014,BorweinTam2015,ArtachoBorweinTam2016}.
For general constrained convex optimization problem, Douglas-Rachford splitting (DRS) \cite{lions1979splitting} and alternating direction method of multipliers (ADMM) \cite{Glowinski1975admm,gabay1976dual} are popular choices of algorithm, and their behavior on infeasible primal or dual problems has been recently analyzed \cite{eckstein1992douglas,BauschkeHareMoursi2014,raghunathan2014infeasibility,BanjacGoulart2019,LiuRyuYin2019,BauschkeMoursi2020,Banjac2021,BanjacLygeros2021asymptotic,BauschkeMoursi2021,ODonoghue2021operator,moursi2022douglas}.
Recently, PDHG \cite{chambolle2011first} has been used as a first-order algorithm solving possibly inconsistent LP and QP \cite{Applegate2021}.

\paragraph{Accelerated fixed-point iterations.}
Picard iteration converges when the operator $\opT$ is contractive, but does not converge with nonexpansivity alone.
If $\opT$ is averaged, fixed-point residual of Picard iteration converges in $\bigO(1/k)$ rate \cite{davis2015convergence}.
But rather than adding conditions on operators, interpolation or extrapolation schemes \cite{krasnosel1955two,mann1953mean,anderson1965,ishikawa1976fixed,xu2004viscosity,mainge2008convergence,dong2018modified,shehu2018convergence,themelis2019supermann,reich2021inertial,walker2011anderson,zhang2020globally,shehu2020inertial,shehu2020iterative,scieur2020regularized,barre2020convergence} may result in faster convergence rate, which is the case for Halpern iteration \cite{halpern1967fixed}, which exhibits $\cO(1/k^2)$ rate \cite{sabach2017first,Lieder2021halpern}.

For the inconsistent fixed-point iteration, the rate of convergence to infimal displacement vector is measured.
Unlike the convergence itself \cite{bailion1978asymptotic}, the $\bigO(1/k)$ rate of convergence was not known until late 2010s \cite{LiuRyuYin2019}.
Another sequence converging to infimal displacement vector is normalized iterates, and it is proven to converge in $\bigO(1/k^2)$ rate \cite{Applegate2021}.

\paragraph{Complexity lower bound.}
Using the information-based complexity framework \cite{nemirovski1992information}, lower bounds to the iteration complexity has been thoroughly studied for first-order convex optimization methods \cite{Nesterov2004convex,drori2017exact,carmon2020stationary1,drori2020stoc-complexity,carmon2021stationary2,dragomir2022optimal,drori2022oracle,yoon2021accelerated,ParkRyu2022}.
For the fixed-point iterations, \citet{Diako2020halpern} first proved $\Omega(1/k^2)$-lower bound, and \citet{ParkRyu2022} later closed the constant gap by showing that Halpern iteration of \citet{Lieder2021halpern} has exactly matching $\Theta(1/k^2)$-complexity to the lower bound of \citet{ParkRyu2022}.
However, these works are restricted to the consistent fixed-point iterations.

\paragraph{Performance estimation problem (PEP).}
From the seminal work of \citet{drori2014performance}, performance estimation problem (PEP) has been widely used to obtain the worst-case complexity of algorithms, including first-order methods \cite{kim2017convergence,taylor2017smooth,de2017worst,kim2018another,taylor2018exactpgm,barre2020complexity,de2020worst,kim2021optimizing,abbaszadehpeivasti2022conditions,abbaszadehpeivasti2022exact,barre2022principled,kamri2022worst,rotaru2022tight,gupta2023nonlinear}, operator splitting methods \cite{ryu2020operator}, minimax algorithms \cite{abbaszadehpeivasti2021rate,gorbunov2022last,zamani2022convergence}, proximal point methods \cite{gu2020tight,kim2021accelerated,gu2022tight,gu2023tight}, decentralized methods \cite{colla2021automated,colla2022automated,colla2022automatic}, coordinate descent methods \cite{abbaszadehpeivasti2022convergence}, and even the continuous-time models \cite{moucer2022systematic}.
PEP also finds the optimal method with optimal worst-case complexity \cite{drori2016optimal-kelley,kim2016optimized,drori2020efficient,taylor2021optimal,kim2021accelerated,ParkRyu2022}, and is even used to construct the Lyapunov function for the proof of convergence \cite{taylor2018lyapunov} and complexity lower bound \cite{dragomir2022optimal}.
All these works assume the existence of the solution or optimal value.
% , and to the best of our knowledge, our work is the first to provide a PEP approach for infeasible problems.

\subsection{Contribution}
We summarize the contribution of this work as follows. First, we prove upper bounds on the rates of convergence of certain sequences to the infimal displacement vector, which can serve as a certificate of infeasibility. In particular, we establish a $\bigO(1/k^2)$-rate for the normalized iterates and $\tilde{\bigO}(1/k^2)$-rate for the fixed-point residual of the Halpern iteration. Second, we extend the performance estimation problem (PEP) methodology to inconsistent fixed-point iterations based on a new interpolability result and demonstrate how we used this methodology to discover the upper bounds. Third, we prove a matching $\Omega(1/k^2)$-complexity lower bound and thereby establish that the  $\bigO(1/k^2)$ upper bound is the optimal accelerated rate. Finally, we complement our theoretical results with a numerical experiment on a decentralized semidefinite program (SDP).

\section{Measure of optimality}
\label{sec:measure}
Consider a nonexpansive operator $\opT\colon\cH\to\cH$.
Then $\fix\opT=\emptyset$ if and only if $0\notin\cR(\opI-\opT)$.
In such case, since $\overline{\cR(\opI-\opT)}$ is a closed convex set, it has a unique minimum element $v = \argmin_{y\in\overline{\cR(\opI-\opT)}} \|y\|^2$.
Roughly, $v$ represents the distance from ${\cR(\opI-\opT)}$ to containing $0$, or $\opT$ being consistent.
As long as $v$ remains nonzero, $\opT$ will never have a fixed point.
For an operator $\opT$ which we do not have full access to, if we are able to obtain $v$ approximately from only a sufficient number of first-order oracle calls, then this will save resources including time and computational power.

Given a nonexpansive operator $\opT\colon\cH\to\cH$, we measure the rate of convergence to $v$ for following sequences.

\begin{definition}
    We call $\frac{x^k-x^0}{\alpha_k}$ with proper scaling factor $\alpha_k>0$ a \emph{normalized iterate}, and call $x^k-\opT x^k$ a \emph{fixed-point residual}.
\end{definition}

Normalized iterate of Picard iteration converges to $-v$ \cite{Applegate2021}, and fixed-point residual of Picard iteration with averaged operator converges to $v$ \cite{RyuLinYin2019douglas}.
Following lemma states that when $v^k$ is either normalized iterate or fixed-point residual at iteration $k$, strong (norm) convergence of $v^k$ to $v$ is equivalent to the convergence of $\|v^k\|$ to $\|v\|$.
Therefore, we measure the rate of convergence for both $\|v^k-v\|^2\to 0$ and $\|v^k\|-\|v\| \to 0$.

\begin{lemma} \label{lem:normconvequiv}
    Let $\opT\colon\cH\to\cH$ be a nonexpansive operator and $v$ be its infimal displacement vector.
    If $v^k$ for $k\in\mathbb{N}$ is either $- \frac{x^k-x^0}{\alpha_k}$ or $x^k-\opT x^k$ with assumption that $\alpha_k>0$ satisfies $- \frac{x^k-x^0}{\alpha_k} \in {\overline{\cR(\opI-\opT)}}$ for all $k\in\mathbb{N}$, then
    \[
        \langle v^k,\, v \rangle \ge \|v\|^2, \quad k=1,2,\dots
    \]
    and
    \[
        \lim_{k\to\infty} v^k = v
        \quad\Leftrightarrow\quad
        \lim_{k\to\infty} \|v^k\| = \|v\|.
    \]
\end{lemma}

Proof of \cref{lem:normconvequiv} is deferred to \cref{sec:proofmeasure}.

\subsection{Comparison of two optimality measures}

In \cref{sec:convergencerate}, we show upper bounds on the two optimality measures $\|v^k-v\|^2$ and $(\|v^k\|-\|v\|)^2$. Since
\[
\|v^k-v\|\ge \|v^k\|-\|v\|
\]
by the triangle inequality, the former is the more rigorous optimality measure in the sense that it is no easier to reduce. This makes intuitive sense as $\|v^k-v\|^2$ corresponds to characterizing the rate of $v^k\rightarrow v$, which is the convergence of both the magnitude and direction of the vectors, while $(\|v^k\|-\|v\|)^2$ corresponds to characterizing the rate of $\|v^k\|\rightarrow \|v\|$, which is the convergence of only the magnitude of the vectors.

% In \cref{sec:lowerbound}, we show lower bounds on $\|v^k-v\|^2$ and $(\|v^k\|-\|v\|)^2$. The lower bounds match the upper bound and we conclude that $\mathcal{O}(1/k^2)$ is the optimal accelerated rate for both optimality measures.

The relative difference $\|v^k-v\|\ge \|v^k\|-\|v\|$ do manifest in terms of different constants. For both optimality measures $\|v^k-v\|^2$ and $(\|v^k\|-\|v\|)^2$, the best known upper bound, presented in \cref{cor:optimalkm}, is $\frac{4}{k^2}\|x^0-x_\star\|^2$. On the other hand, the best lower bound for $\|v^k-v\|^2$ is $\frac{4}{k^2}\|x^0-x_\star\|^2$, while for $(\|v^k\|-\|v\|)^2$ it is $\frac{1}{2k^2}\|x^0-x_\star\|^2$.
So a conclusion of this work is that the two optimality measures are equivalent (up to a constant factor of at most $8$) in their optimal worst-case computational complexity.

% \subsection{Comparing $\|x^k-\opT x^k-v\|$ and $\|x^k-\opT x^k\|-\|v\|$}

% {\color{gray}
% Also note that
% \[
%     \left( \|v^k\| - \|v\| \right)^2 \le \|v^k - v\|^2 \le \|v^k\|^2 - \|v\|^2,
% \]
% which implies that $\|v^k\|\to\|v\|$ is a weaker notion of convenge in terms of measuring the rate of convergence, although the convergence to $v$ itself is guaranteed for all three quantities above.
% }
% \begin{lemma}
%     Let $\opT\colon\cH\to\cH$ be a nonexpansive operator with infimal displacement vector $v$.
%     For any $x\in\dom\opT$,
%     \begin{align*}
%         &\left( \|x-\opT x\| - \|v\| \right)^2 \\
%         &\le
%         \|x-\opT x-v\|^2
%         \le
%         \|x-\opT x\|^2 - \|v\|^2 \\
%         &=
%         \left( \|x-\opT x\| - \|v\| \right) \left( \|x-\opT x\| + \|v\| \right).
%     \end{align*}
% \end{lemma}

% Above lemma shows that the performance measure $\|x-\opT x\|-\|v\|$ is weaker in the sense that the convergence $\|x-\opT x-v\|\to 0$ always implies $\|x-\opT x\|-\|v\|\to 0$, and the convergence rate of the latter directly follows from the former.

\section{Rate of convergence to $v$}
\label{sec:convergencerate}

We study the rate of convergence to $v$ for normalized iterate and fixed-point residual of \eqref{eq:kmiteration} and \eqref{eq:halpern}.
% , and the normalized iterate of general Mann iteration.
In the last part, we deal with the normalized iterate of general Mann iteration.

\subsection{Convergence of KM iteration}
\label{subsec:km}

Consider the \emph{Krasnosel'ski{\u\i}-Mann iteration} \eqref{eq:kmiteration}
\begin{equation} \label{eq:kmiteration}
    x^{k+1} = \lambda_{k+1} x^k + (1-\lambda_{k+1}) \opT x^k, \quad k=0,1,\ldots,
    \tag{KM}
\end{equation}
where $x^0\in \cH$ is a starting point and $\lambda_{k+1}\in[0,1)$.

\begin{theorem}[Convergence rate of normalized iterate] \label{thm:kmaverconv}
    Let $\{x^k\}_{k\in\mathbb{N}}$ be the iterates of \eqref{eq:kmiteration} starting from $x^0\in\cH$.
    For any $\varepsilon>0$ and $x_\varepsilon\in\cH$ such that $\|x_\varepsilon-\opT x_\varepsilon-v\| \le \min \left\{ \frac{\varepsilon^2}{2\|v\|+1},\, 1,\, \varepsilon \right\}$,
    \[
        \left\| \frac{x^k-x^0}{ {\scriptstyle \sum_{i=1}^{k}} (1-\lambda_i)} + v \right\|^2
        \le
        \left( \frac{2}{ {\scriptstyle \sum_{i=1}^{k}} (1-\lambda_i)} \|x^0-x_\varepsilon\| + \varepsilon \right)^2
    \]
    for all $k = 1,2,\dots$.
    If we further assume that $v\in\cR(\opI-\opT)$, then there exists $x_\star\in\cH$ such that $x_\star-\opT x_\star=v$ and 
    \[
        \left\| \frac{x^k-x^0}{ {\scriptstyle \sum_{i=1}^k} (1-\lambda_i)} + v \right\|^2
        \le
        \frac{4}{{\left( {\scriptstyle \sum_{i=1}^k} (1-\lambda_i) \right)}^2} \|x^0-x_\star\|^2
    \]
    for all $k = 1,2,\dots$.
\end{theorem}

\begin{theorem}[Convergence rate of fixed-point residual] \label{thm:kmfprconv}
    Let $\{x^k\}_{k\in\mathbb{N}}$ be the iterates of \eqref{eq:kmiteration} starting from $x^0\in\cH$ and
    $k_0 = \min\{ i\in\mathbb{N} \mid \lambda_i > 0 \}$.
    For any $\varepsilon>0$ and $x_\varepsilon\in\cH$ such that $\|x_\varepsilon-\opT x_\varepsilon-v\| \le \min \left\{ \frac{\varepsilon^2}{2\|v\|+1},\, 1,\, \varepsilon \right\}$,
    \begin{align*}
        & \left( \sum_{i=0}^{k} \frac{\lambda_{i+1} (1-\lambda_{i+1})}{\sum_{i=0}^{k} \lambda_{i+1} (1-\lambda_{i+1})} \|x^{i} - \opT x^{i} - v\| \right)^2 \\
        &\qquad\qquad \le
        \left( { \frac{1}{\sqrt{ {\scriptstyle \sum_{i=0}^{k}} \lambda_{i+1} (1-\lambda_{i+1})}}} \|x^0-x_\varepsilon\| + \varepsilon \right)^2
    \end{align*}
    and
    \begin{align*}
        &\left( \| x^k-\opT x^k \| - \|v\| \right)^2 \\
        &\qquad\qquad \le
        \left( {\frac{1}{ \sqrt{{\scriptstyle \sum_{i=0}^{k}}\lambda_{i+1} (1-\lambda_{i+1})}}} \|x^0-x_\varepsilon\| + \varepsilon \right)^2
    \end{align*}
    for $k\ge k_0$.
    If we further assume that $v\in\cR(\opI-\opT)$, then there exists $x_\star\in\cH$ such that $x_\star-\opT x_\star=v$, 
    \begin{align*}
        & \left( \sum_{i=0}^{k} \frac{\lambda_{i+1} (1-\lambda_{i+1})}{\sum_{i=0}^{k} \lambda_{i+1} (1-\lambda_{i+1})} \|x^{i} - \opT x^{i} - v\| \right)^2 \\
        &\qquad\qquad \le
        \frac{1}{{\scriptstyle \sum_{i=0}^{k}} \lambda_{i+1} (1-\lambda_{i+1})} \|x^0-x_\star\|^2,
    \end{align*}
    and
    \begin{align*}
        \left( \| x^k-\opT x^k\| - \|v\| \right)^2
        &\le
        \frac{1}{{ \scriptstyle \sum_{i=0}^{k}} \lambda_{i+1} (1-\lambda_{i+1})} \|x^0-x_\star\|^2
    \end{align*}
    for $k\ge k_0$.
\end{theorem}

We defer the proofs to \cref{sec:proofkm}.
Note that Theorems~\ref{thm:kmaverconv} and \ref{thm:kmfprconv} imply the convergence of normalized iterate and fixed-point residual to $v$ respectively when $\sum_{k=1}^\infty (1-\lambda_k) = \infty$ and $\sum_{k=1}^\infty \lambda_k (1-\lambda_k) = \infty$.
We also point out that the bound on the Ces\`aro mean in \cref{thm:kmfprconv} is practically useful when we use the randomized iterate selection technique of \citet{ghadimi2013stochastic,ghadimi2016accelerated}:
choosing $\bar{k} \in \{1,2,\dots,k\}$ with probability proportional to $\lambda_{\bar{k}+1} (1-\lambda_{\bar{k}+1})$, fixed-point residual $x^{\bar{k}} - \opT x^{\bar{k}}$ of $\bar{k}$-th iterate will yield the same rate of convergence as \cref{thm:kmfprconv}.

\begin{corollary}
\label{cor:optimalkm}
    Let $\{x^k\}_{k\in\mathbb{N}}$ be the iterates of \eqref{eq:kmiteration} starting from $x^0\in\cH$.
 Assume that $v\in\cR(\opI-\opT)$.
 The bound of \cref{thm:kmaverconv} is optimized at $\lambda_k=0$ for all $k\in\mathbb{N}$ with
    \[
        \left\| \frac{x^k-x^0}{k} + v \right\|^2
        \le
        \frac{4}{k^2} \|x^0-x_\star\|^2.
    \]
The bound of \cref{thm:kmfprconv} is optimized at $\lambda_k=\frac{1}{2}$ for all $k\in\mathbb{N}$ with
    \[
        \frac{1}{k+1} \sum_{i=0}^k \| x^i-\opT x^i - v\|^2
        \le
        \frac{4}{k+1} \|x^0-x_\star\|^2
    \]
    and
    \[
        (\|x^k-\opT x^k\| - \|v\|)^2
        \le
        \frac{4}{k+1} \|x^0-x_\star\|^2.
    \]
\end{corollary}

    \cref{cor:optimalkm} recovers the rates of \citep[Theorem~3]{LiuRyuYin2019} and \citep[Theorem~3]{Applegate2021}.
    % \cref{thm:kmfprconv} generalizes the result for consistent fixed-point iteration proven by \citet{DavisYin2016} to the inconsistent setup.
To clarify, we view the results of Sections~\ref{subsec:halpern} and \ref{sec:lowerbound} to be the major contributions of this work.
Our contribution of \cref{subsec:km}, presented in Theorems~\ref{thm:kmaverconv} and \ref{thm:kmfprconv}, is to generalize the results of \citep{LiuRyuYin2019,Applegate2021}
to the KM iteration with $\{\lambda_k\}_{k\in \mathbb{N}}$ that varies with $k$.

\paragraph{Counterexample.}
Theorems~\ref{thm:kmaverconv} and \ref{thm:kmfprconv} show that convergence of normalized iterates requires ${\sum_{k=1}^\infty(1-\lambda_k) = \infty}$, while convergence of fixed-point residual requires the stronger condition ${\sum_{k=1}^\infty\lambda_k(1-\lambda_k) = \infty}$. The following demonstrates that it is possible for the normalized iterates to converge while the fixed-point residual diverges.

Define $\opT \colon \RR^3 \to \RR^3$ as $\opT(x,y,z) = (-y, x, z-1)$.
Then $\cR(\opI-\opT) = \RR^2 \times \{1\}$ and $v = (0,0,1)$.
Let $\{ (x^k,y^k,z^k) \}_{k \in \mathbb{N} \cup \{0\}}$ be a sequence of iterates generated by \eqref{eq:kmiteration} with $\opT$ and $\lambda_k=0$ for all $k\in\mathbb{N}$ starting from $(x^0,y^0,z^0) = (1,0,0)$.
Then
\[
    (x^k, y^k, z^k)
    = \left(
        \cos \frac{k\pi}{2},\,
        \sin \frac{k\pi}{2},\,
        -k
    \right),
\]
and the normalized iterates converge to $-v$.
However, 
% \begin{align*}
%     &(x^k,y^k,z^k) - \opT (x^k,y^k,z^k) \\
%     &=
%     \left(
%         \cos\frac{k\pi}{2} + \sin\frac{k\pi}{2},\,
%         - \cos\frac{k\pi}{2} + \sin\frac{k\pi}{2},\,
%         1
%     \right).
% \end{align*}
$\left\| (x^k,y^k,z^k) - \opT (x^k,y^k,z^k) - v\right\| = \sqrt{2}$ for all $k\in\mathbb{N}$, so the fixed-point residual does not converge to $v$.

\subsection{Convergence of Halpern iteration}
\label{subsec:halpern}

Consider the \emph{Halpern iteration} \eqref{eq:halpern}
\begin{equation}
    x^{k+1} = \lambda_{k+1} x^0 + (1-\lambda_{k+1}) \opT x^k, \quad k=0,1,\ldots,
    \tag{Halpern}
    \label{eq:halpern}
\end{equation}
where $x^0\in \cH$ is a starting point and $\lambda_{k+1}\in[0,1)$.
Note that \eqref{eq:fpi} corresponds $\lambda_k\equiv0$ and OHM \cite{Lieder2021halpern} corresponds to $\lambda_k=\frac{1}{k+1}$.
Define $\theta_0=0$ and 
    \[
        \theta_k = \sum_{n=1}^k (1-\lambda_k)(1-\lambda_{k-1})\cdots(1-\lambda_{k-n+1})
    \]
    for $k=1,2,\dots$.
\begin{lemma} \label{lem:haltheta}
    For $k=0,1,\dots$, 
    \[
        \theta_{k+1} = (1-\lambda_{k+1})(1+\theta_k).
    \]
    If $\lambda_k\equiv0$, then $\theta_k=k$.
    If $\lambda_k=\frac{1}{k+1}$, then $\theta_k=\frac{k}{2}$.
\end{lemma}

\begin{theorem}[Convergence rate of normalized iterate]
    \label{thm:halaverconv}
        Let $\{x^k\}_{k\in\mathbb{N}}$ be the iterates of \eqref{eq:halpern} starting from $x^0\in\cH$.
    % For any $\varepsilon>0$, there exists $x_\varepsilon\in\cH$ such that $\|x_\varepsilon-\opT x_\varepsilon-v\|<\varepsilon$ and
    For any $\varepsilon>0$ and $x_\varepsilon\in\cH$ such that $\|x_\varepsilon-\opT x_\varepsilon\|^2 - \|v\|^2 \le \varepsilon^2$,
    \[
        \left\| \frac{x^k-x^0}{\theta_k} + v \right\|^2
        \le
        \left( \frac{2}{\theta_k} \|x^0-x_\varepsilon\| + \varepsilon \right)^2
    \]
    for $k=1,2,\dots$.
    If we further assume that $v\in\cR(\opI-\opT)$, then there exists $x_\star\in\cH$ such that $x_\star-\opT x_\star=v$ and 
    \[
        \left\| \frac{x^k-x^0}{\theta_k} + v \right\|^2
        \le
        \frac{4}{\theta_k^2} \|x^0-x_\star\|^2
    \]
    for $k=1,2,\dots$.
\end{theorem}

% \begin{corollary}[Convergence of normalized iterate]
%     If $\{x^k\}_{k\in\mathbb{N}}$ is a sequence of iterates generated by \eqref{eq:halpern} starting from $x^0\in\cH$ and $\lim_{k\to\infty} \theta_k = \infty$, then
%     \[
%         \lim_{k\to\infty} \frac{x^k-x^0}{\theta_k} = - v.
%     \]
% \end{corollary}

We defer the proofs to \cref{sec:proofhalpern}.
Note that the normalized iterates converge to $-v$ if $\theta_k \rightarrow  \infty$, which, in particular,  happens if $\lambda_k \rightarrow 0$.
    See \cref{lem:halpernlambda}.

% \begin{remark} \label{rem:halpernlambda}
%     The condition $\lim_{k\to\infty} \theta_k = \infty$ can be implied from a simpler condition $\lim_{k\to\infty} \lambda_k = 0$.
%     See \cref{lem:halpernlambda}.
% \end{remark}

% Choosing $\lambda_k$ as in \citet{Lieder2021} results in the convergence of fixed-point residual $x^k-\opT x^k$ to $v$ as well.
\begin{theorem}[Convergence rate of fixed-point residual] \label{thm:halfprconv}
    Let $\{x^k\}_{k\in\mathbb{N}}$ be the iterates of \eqref{eq:halpern} starting from $x^0\in\cH$  with $\lambda_k = \frac{1}{k+1}$.
    % For any $\varepsilon>0$, there exists $x_\varepsilon\in\cH$ such that $\|x_\varepsilon-\opT x_\varepsilon-v\|<\varepsilon$,
    For any $\varepsilon>0$ and $x_\varepsilon\in\cH$ such that $\|x_\varepsilon - \opT x_\varepsilon\|^2 - \|v\|^2 \le \cO(\varepsilon^2)$,
    we have
    \[
        \left( \| x^k - \opT x^k\| - \|v\| \right)^2 
        \le
        \left( \frac{4}{k} \|x^0-x_\varepsilon\| + \varepsilon \right)^2
    \]
    and
    \begin{align*}
        &\|x^{k}-\opT x^{k}-v\|^2 \\  
        % \|x^{k}-\opT x^{k}-v\|^2
        &\le
        \left(\frac{\sqrt{\sum_{n=1}^{k} \frac{1}{n} + 4}+1}{k+1}\right)^2 \|x^0-x_\varepsilon\|^2 + \varepsilon
        % \\ &\le
        % \frac{2\left(\sum_{n=1}^k \frac{1}{n} + 5 \right)}{(k+1)^2} \|x^0-x_\varepsilon\|^2 + \varepsilon.
    \end{align*}
    for $k=1,2,\dots$.
    If we further assume that $v\in\cR(\opI-\opT)$, then there exists $x_\star\in\cH$ such that $x_\star-\opT x_\star=v$, 
    \[
        \left( \| x^k - \opT x^k\| - \|v\| \right)^2
        \le
        \frac{16}{k^2} \|x^0-x_\star\|^2,
    \]
    and
    \begin{align*}
        \|x^{k}-\opT x^{k}-v\|^2
        &\le
        \left(\frac{\sqrt{\sum_{n=1}^{k} \frac{1}{n} + 4}+1}{k+1}\right)^2 \|x^0-x_\star\|^2
        % \\ &\le
        % \frac{2\left(\sum_{n=1}^k \frac{1}{n} + 5 \right)}{(k+1)^2} \|x^0-x_\star\|^2.
    \end{align*}
    for $k=1,2,\dots$.
    % This implies $\lim_{k\to\infty} (x^k-\opT x^k) = v$.
\end{theorem}

\begin{proof}[Proof outline]
    Consider a potential function $V^k$ defined as
    \begin{align*}
        &\,V^k \\
        &=
        (k+1)\left\{ k\|x^k-\opT x^k\|^2 + 2\langle x^k-\opT x^k, x^k-x^0 \rangle \right\} \\
        &
        + k(k+1) \left\langle - \frac{2}{k}(x^k-x^0) - \left(x_\varepsilon-\opT x_\varepsilon\right), x_\varepsilon-\opT x_\varepsilon \right\rangle \\
        &
        + \frac{2(k+1)}{k} \left\|x^k-x_\varepsilon + \frac{k}{2} \left(x_\varepsilon-\opT x_\varepsilon\right) \right\|^2 \\
        &\quad
        -  \left( {\sum_{n=1}^{k}} \frac{1}{n} \right) \|x^0-x_\varepsilon\|^2
    \end{align*}
    for all $k\in\mathbb{N}$.
    We can show $V^k\le V^{k-1}\le \cdots\le V^1$.
    From $V^k \le V^1 \le 3\|x^0-x_\varepsilon\|^2$, we obtain the desired convergence rate.
    When there exists $x_\star$ such that $v=x_\star - \opT x_\star$, use $x_\star$ instead of $x_\varepsilon$.
    The detailed proof is deferred to \cref{sec:proofhalpern}.
\end{proof}

The precise form of the $\cO(\varepsilon^2)$-term in \cref{thm:halfprconv} is stated in the proof, which is deferred to \cref{sec:proofhalpern}.
Note that in $V^k$, the first term, written as $(k+1)\left\{\cdots \right\}$, is the potential function that was used in prior work \cite{Diako2020halpern,ParkRyu2022} to analyze the convergence of consistent fixed-point iterations.
So the first term is known to be nonincreasing, and the three additional terms are required to adapt the proof to the inconsistent case.

\begin{corollary} \label{cor:hallinearconv}
    Let $\{x^k\}_{k\in\mathbb{N}}$ be the iterates of \eqref{eq:halpern} starting from $x^0\in\cH$  with $\lambda_k = \frac{1}{k+1}$.
    For any $\varepsilon>0$, there exists $x_\varepsilon\in\cH$ such that $\|x_\varepsilon-\opT x_\varepsilon-v\|<\varepsilon$,
    \[
        \left\| \frac{2(x^k-x^0)}{k} + v \right\|^2 \le \left( \frac{4}{k} \|x^0-x_\varepsilon\| + \varepsilon \right)^2,
    \]
    \[
        \left( \| x^k - \opT x^k\| - \|v\| \right)^2
        \le
        \left( \frac{4}{k} \|x^0-x_\varepsilon\| + \varepsilon \right)^2,
    \]
    and
    \begin{align*}
        &\|x^{k}-\opT x^{k}-v\|^2 \\
        &\le
        \left(\frac{\sqrt{\sum_{n=1}^{k} \frac{1}{n} + 4}+1}{k+1}\right)^2 \|x^0-x_\varepsilon\|^2 + \varepsilon
        % \\ &\le
        % \frac{2\left(\sum_{n=1}^k \frac{1}{n} + 5 \right)}{(k+1)^2} \|x^0-x_\varepsilon\|^2 + \varepsilon.
    \end{align*}
        for $k=1,2,\dots$.
    If we further assume that $v\in\cR(\opI-\opT)$, then there exists $x_\star\in\cH$ such that $x_\star-\opT x_\star=v$,
    \[
        \left\| \frac{2(x^k-x^0)}{k} + v \right\|^2 \le \frac{16}{k^2} \|x^0-x_\star\|^2,
    \]
    \[
        \left( \| x^k - \opT x^k\| - \|v\| \right)^2 \le \frac{16}{k^2} \|x^0-x_\star\|^2,
    \]
    and
    \begin{align*}
        \|x^{k}-\opT x^{k}-v\|^2
        &\le
        \left(\frac{\sqrt{\sum_{n=1}^{k} \frac{1}{n} + 4}+1}{k+1}\right)^2 \|x^0-x_\star\|^2
    \end{align*}
    for $k=1,2,\dots$.
\end{corollary}

An observation we point out is that when $\opT$ is an affine operator, the normalized iterate $- \frac{x^{k+1} - x^0}{k+1}$ of Picard iteration coincides with the fixed-point residual $x^k - \opT x^k$ of \eqref{eq:halpern} with $\lambda_k = \frac{1}{k+1}$. 
See \cref{lem:picardequiv}.

\subsection{Convergence of Mann iteration}
\label{subsec:mann}

The \emph{Mann iteration} \eqref{eq:manniteration}
\begin{equation*}
    x^k = {\small \sum_{i=0}^{k-1}} \nu^k_i \, \opT x^{i-1},
    \tag{Mann} \label{eq:manniteration}
\end{equation*}
where $\nu_i^k \ge 0$ for $i=0,\dots,k$ and $k=1,2,\dots$, $\sum_{i=0}^k \nu^k_i = 1$ for $k=1,2,\dots$, and $\opT x^{-1} := x^0$,
is a further general class of iterations including \eqref{eq:kmiteration} and \eqref{eq:halpern}. \cref{lem:mannnormalizingfactor} of \cref{subsec:proofmann} shows that there exists positive sequence $\{\alpha_k\}_{k\in\mathbb{N}}$ that depends on $\{\nu_i^k\}_{\substack{0\le i\le k,k\in \mathbb{N}}}$ such that 
\[
    - \frac{x^k-x^0}{\alpha_k} \in \overline{\cR(\opI-\opT)}, \quad k=1,2,\dots.
\]
Furthermore, 
\cref{thm:mannaverconv} of \cref{subsec:proofmann} shows that
\[
    \left\| \frac{x^k-x^0}{\alpha_k} + v \right\|^2 \le \left( \frac{2}{\alpha_k} \|x^0-x_\varepsilon\| + \varepsilon \right)^2
\]
and the normalized iterate converges to $-v$ if $\alpha_k \rightarrow \infty$.
This result generalizes the convergence results of Theorems~\ref{thm:kmaverconv} and \ref{thm:halaverconv} respectively for \eqref{eq:kmiteration} and \eqref{eq:halpern}.

\section{PEP with possibly infeasible operators}
\label{sec:pep}
Instrumental in the discovery of the results of \cref{sec:convergencerate} was the use of the performance estimation problem (PEP) \cite{drori2014performance,taylor2017smooth}. Loosely speaking, the PEP is a computer-assisted methodology for finding optimal methods by numerically solving semidefinite programs \cite{drori2014performance,drori2016optimal-kelley,kim2016optimized,taylor2018exactpgm,drori2020efficient,kim2021optimizing,kim2021accelerated,ParkRyu2022}. In prior work, PEP had been utilized in the analysis of \emph{consistent} monotone inclusion and fixed-point problems \citep{ryu2020operator,kim2021accelerated,ParkRyu2022}.
In this section, we describe how to apply the PEP methodology in the analysis of algorithms for \emph{inconsistent} problems.

\subsection{Interpolation result}
\label{subsec:interpolation}

The performance estimation problem framework relies on certain interpolation results. The following result strengthens the prior interpolation result of \citep[Fact 2]{ryu2020operator} by additionally restricting the range of the extension and thereby allows us to control the infimal displacement vector of the interpolation.

\begin{theorem}[Interpolability]
    \label{thm:interpolate}
    Let $\{ (x_i,\, y_i) \}_{i\in I} \subset \cH\times\cH$ be a set of vectors with index set $I$ such that
    \[
        \|y_i - y_j\| \le \|x_i - x_j\|,
        \quad \forall \,i, j \in I.
    \]
    Let $C=\clconv \{x_i - y_i\}_{i\in I} \subseteq \cH$, where $\clconv$ denotes the closure of the convex hull.
    \begin{itemize}
    \item[(i)] There exists a nonexpansive $\widetilde{\opT}\colon\cH\to\cH$ such that
    \[
        y_i = \widetilde{\opT} x_i, \quad \forall\, i\in I
    \]
    and $v = \Proj_{C}(0)$ is its infimal displacement vector.
    
    \item[(ii)] If we further assume that $v = x_\star - y_\star$, $\star\in I$ and
    \[
        \langle x_i - y_i,\, v \rangle \ge \|v\|^2, \quad \forall\, i\in I,
    \]
    then there exists a nonexpansive $\widetilde{\opT}\colon\cH\to\cH$ such that
    \[
        y_i = \widetilde{\opT} x_i, \quad \forall\, i\in I
    \]
    and $v$ is its infimal displacement vector.
    \end{itemize}
\end{theorem}

We defer the proof to \cref{sec:proofinterpolation}. The key insight is to use the range/domain-restricting extension of \citep{reich2005,bauschke2007fenchel}, construction of which, in turn, relies on the Fitzpatrick function \citep{Fitzpatrick1988_representing}.

\subsection{PEP formulation}
\label{subsec:pepformulation}

We now describe the PEP formulation with inconsistent operators through an example.
Consider \eqref{eq:halpern} with $\lambda_k = \frac{1}{k+1}$, which we refer to as the optimized Halpern method (OHM) of \cite{Lieder2021halpern}.
Let $k\in\mathbb{N}$ and define the index set $I = \{0,1,\ldots,k,\star\}$.
We consider nonexpansive operators $\opT$ that have an infimal displacement vector $v$ and a point $x_\star\in\cH$ such that $v = x_\star-\opT x_\star$.
The goal is to find the worst-case instance of $\opT$ such that $\|x^k - \opT x^k - v\|^2$ is maximized.

We start from the infinite-dimensional performance estimation problem 
\[
    \begin{array}{ll}
        \underset{\opT}{\mbox{maximize}} & \|x^k - \opT x^k - v\|^2 \\
        \mbox{subject to} & \text{$\opT\colon\cH\to\cH$ is nonexpansive} \\
        & v = \Pi_{\overline{\cR(\opI-\opT)}}(0) = x_\star - \opT x_\star \\
        & x^{n+1} = {\frac{n+1}{n+2}} \opT x^n + \frac{1}{n+2} x^0 \\
        & \|x^0-x_\star\|^2 \le R^2
    \end{array}
\]
% \begin{align*}
%     \underset{\opT}{\mbox{maximize}}
%     &\quad
%     \|x^k - \opT x^k - v\|^2 \\
%     \mbox{subject to} 
%     &\quad
%     \text{$\opT\colon\cH\to\cH$ is nonexpansive} \\
%     &\quad
%     v = \Pi_{\overline{\cR(\opI-\opT)}}(0) = x_\star - \opT x_\star \\
%     &\quad
%     x^{n+1} = {\frac{n+1}{n+2}} \opT x^n + \frac{1}{n+2} x^0\\
%     % , \quad n=0,1,\ldots,k-1 \\
%     &\quad
%     \|x^0-x_\star\|^2 \le R^2,
% \end{align*}
where $n=0,1,\dots,k-1$.
Using \cref{thm:interpolate} and scaling by $R$, we get the equivalent non-convex finite-dimensional problem
\[
    \begin{array}{ll}
        \underset{(x^i,\,y^i)_{i\in I}}{\mbox{maximize}} 
        & \|x^k-y^k-v\|^2 \\
        \mbox{subject to} 
        & \| y^i - y^j \|^2 \le \| x^i - x^j \|^2, \quad \forall i,j\in I,\, i\neq j \\
        & v = x_\star - y_\star \\
        & \left\langle x^i-y^i,\, v \right\rangle \ge \|v\|^2, \quad \forall i \in I \\
        & x^{n+1} = \frac{n+1}{n+2} y^n + \frac{1}{n+2} x^0 \\
        & \|x^0-x_\star\|^2 \le 1
    \end{array}
\]
where $n=0,1,\dots,k-1$.
Next, consider the following Gram matrix $Z = G^\intercal G \in \SS^{k+3}_+$, where
\begin{equation}
    G = 
    \begin{bmatrix}
        v^0 & \cdots & v^k & v & x^0-x_\star
    \end{bmatrix}
    \label{eq:g-def}
\end{equation}
with $v^i = x^i-y^i$ for $i=0,1,\ldots,k$.
% Gram matrix formulation drops the dependency on the dimension $d$ of the underlying space $\cH$ of operator $\opT$, only imposing such condition to hold implicitly.
% We close this section with the lemma guaranteeing the equivalency between the last formulation and the final semidefinite program formulation to be formalized at the end.
We finally obtain the following equivalent (convex) semidefinite program, 
\[
    \begin{array}{ll}
        \underset{Z\in\SS^{k+3}_+}{\mbox{maximize}} & \tr(C_k Z) \\
        \mbox{subject to} & \tr(A_{i,j} Z) \ge 0, \quad \forall\, i,j \in I\setminus\{\star\},\, i\neq j \\
        & \tr(A_{i,\star} Z) \ge 0, \quad \forall\, i\in I \setminus\{\star\} \\
        & \tr(B_{i} Z) \le 0, \quad \forall\, i\in I\setminus\{\star\} \\
        & \tr(D_0 Z) \le 1,
    \end{array}
\]
where $A_{i,j}$, $A_{i,\star}$, and $B_i$ for $i,j\in I\setminus\{\star\}$, $C_k$, and $D_0$ in $\SS^{k+3}$ are all defined in \cref{subsec:proofpepformulation}. The details and the subtleties of deriving the SDP representation are also further discussed in \cref{subsec:proofpepformulation}.

\begin{figure*}[ht]
    \centering
    \includegraphics[width=.32\textwidth]{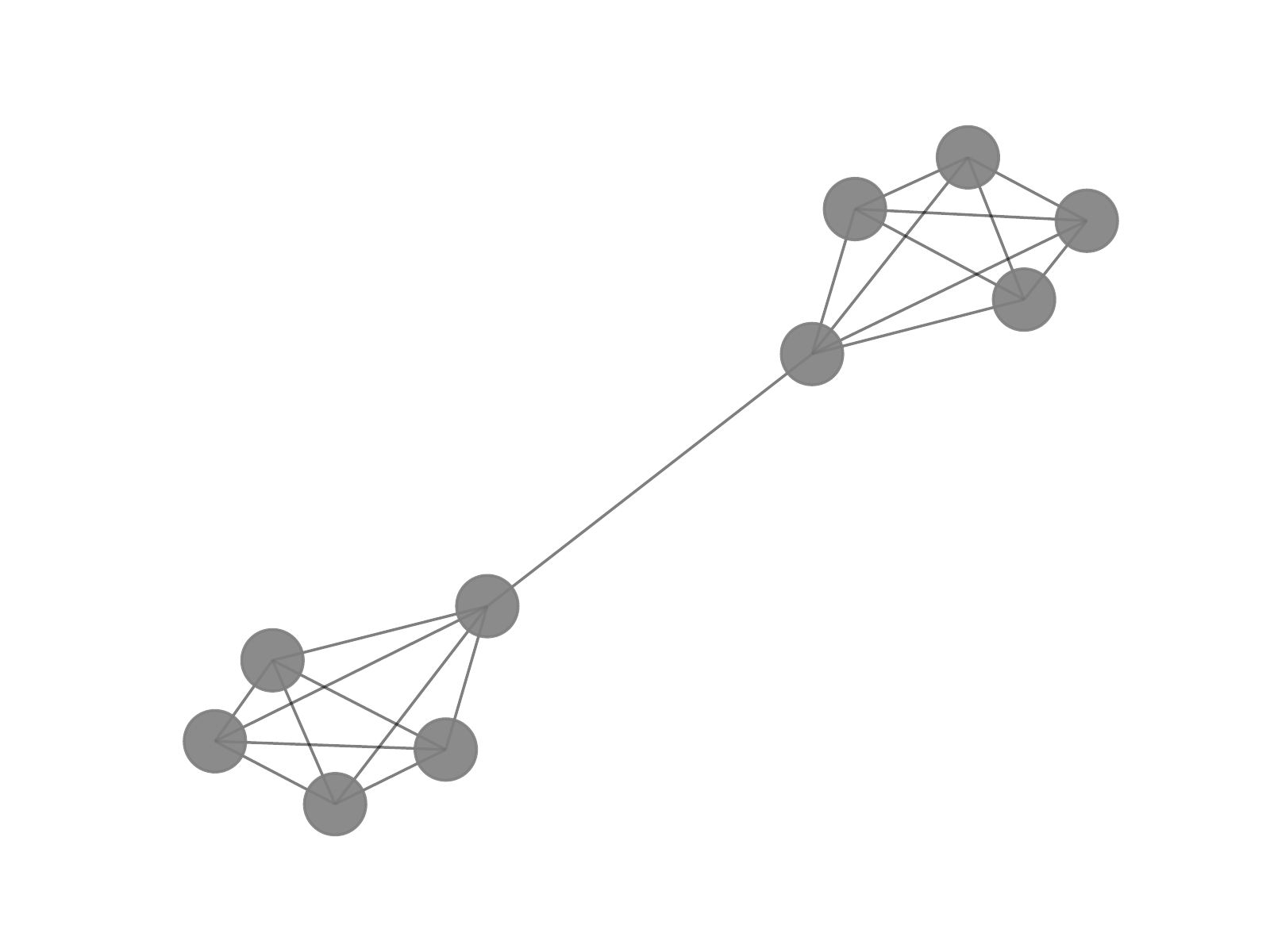}
    \hfill
    \includegraphics[width=.32\textwidth]{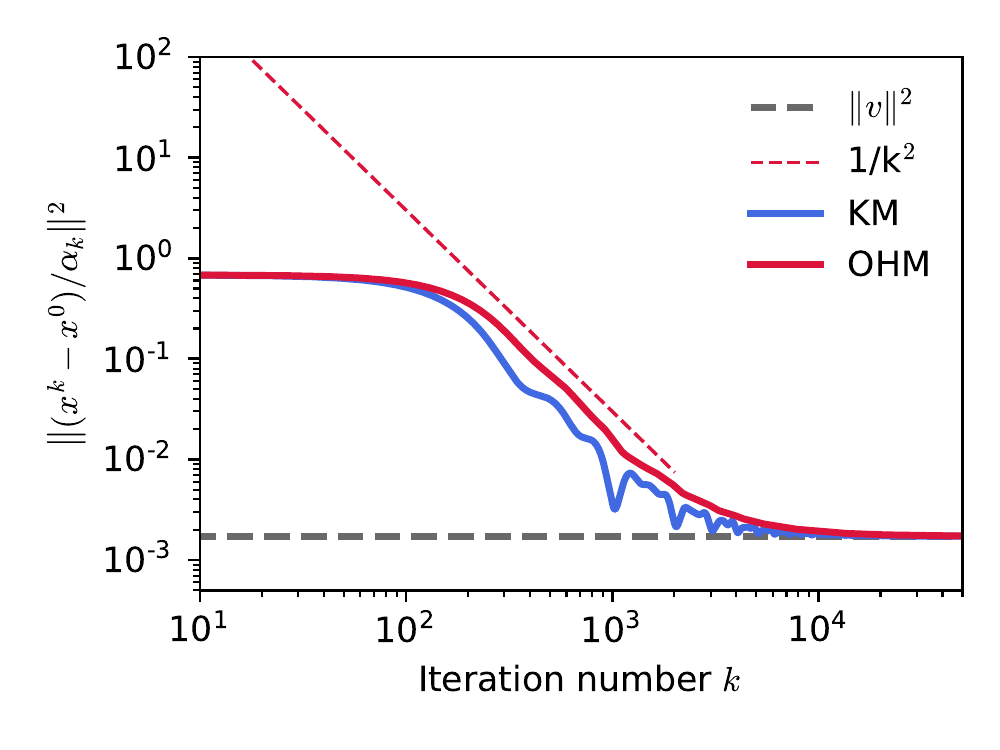}
    \hfill
    \includegraphics[width=.32\textwidth]{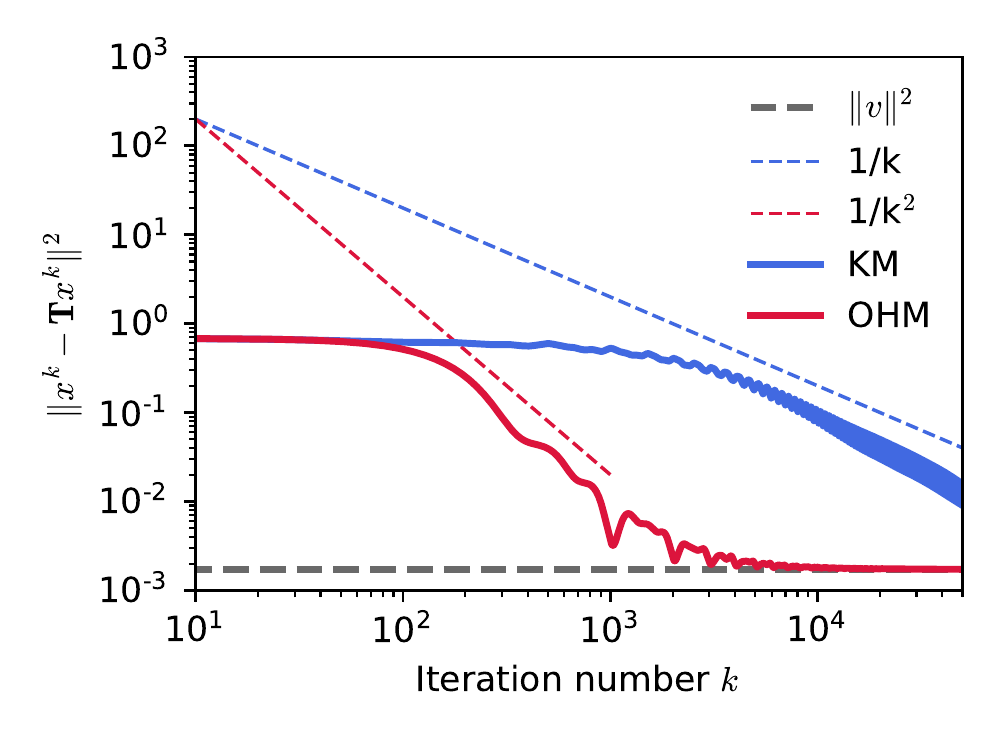}
    \caption{
        Solving SDP with $50,000$ iterations of PG-EXTRA \eqref{eq:fpi} and OHM with PG-EXTRA \eqref{eq:halpern}. We use an infeasible instance, whose setups are described in \cref{sec:proofexperiments}. Parameters are $n=10$, $m=11$, $p=10$ with $\alpha=\beta=0.01$.
        (Left) Network graph.
        (Middle) Squared norm of normalized iterate $\|(x^k-x^0)/\alpha_k\|^2$.
        (Right) Squared norm of fixed-point residual $\|x^k-\opT x^k\|^2$.
        }
    \label{fig:sdp_experiment}
\end{figure*}

\section{Complexity lower bound}
\label{sec:lowerbound}
In this section, we establish a lower bound on the computational complexity of approximating the infimal displacement vector $v$.
Following the information-based complexity framework \cite{nemirovski1992information}, we begin by considering algorithms satisfying the linear span condition
\[
    x^{k+1} = x^0 + \Span\{ x^0-\opT x^0, x^1-\opT x^1, \ldots, x^k - \opT x^k \},
    \tag{span}
    \label{def:span}
\]
which covers a broad range of fixed-point iterations including \eqref{eq:kmiteration}, \eqref{eq:halpern}, \eqref{eq:manniteration}, Anderson acceleration and many more.
We then remove the linear span assumption and expand the class of algorithm to all ``deterministic fixed-point iterations.''

\begin{theorem} \label{thm:lowerboundspan}
    Let $k\in\mathbb{N}$, $x^0=0\in\cH$, and $v\in\cH$, where $\dim\cH \ge k+1$.
    Then, there exists a nonexpansive operator $\opT\colon\cH\to\cH$ and $x_\star\in\cH$ such that $v = x_\star - \opT x_\star$, $v$ becomes the infimal displacement vector of $\opT$, and
    \[
        \left(\left\| \sum_{i=0}^{k-1} \nu_i (x^i-\opT x^i) \right\| - \|v\|\!\right)^{\!2}
        \!\!\ge
        \frac{1}{2k^2} \|x^0-x_\star\|^2
    \]
    and
    \[
        \left\| \sum_{i=0}^{k-1} \nu_i (x^i-\opT x^i) - v \right\|^2
        \ge
        \frac{4}{k^2} \|x^0-x_\star\|^2
    \]
    hold for any iterates $\{x^n\}_{n=0}^{k-1}$ satisfying \eqref{def:span}
    and any choice of real numbers $\{\nu_i\}_{i=0}^{k-1}$ such that $\sum_i \nu_i = 1$.
\end{theorem}

\begin{proof}[Proof outline]
    We construct a nonexpansive operator $\opT\colon\RR^{k+1}\to\RR^{k+1}$ with its infimal displacement $\tilde{v}=(0,\ldots,0,\|v\|)$.
    % , using the worst-case operator of \citet{ParkRyu2022}.
    Then, we choose an orthogonal matrix $U\in\RR^{(k+1)\times (k+1)}$ such that $U^\intercal U = \opI$ and $U^\intercal \tilde{v} = v$ to construct a nonexpansive operator $\opT_U = U \opT U^\intercal$ whose infimal displacement vector is $v$.
    Our specific construction, inspired by \citet{ParkRyu2022}, is
    \[
        x - \opT x = 
        \underbrace{\begin{bmatrix}
            1 & 0 & 0 & \hdots & 0 & 1 & 0 \\
            -1 & 1 & 0 & \hdots & 0 & 0 & 0 \\
            0 & -1 & 1 & \hdots & 0 & 0 & 0 \\
            \vdots & \vdots & \vdots & \ddots & \vdots & \vdots & \vdots \\
            0 & 0 & 0 & \hdots & 1 & 0 & 0 \\
            0 & 0 & 0 & \hdots & -1 & 1 & 0 \\
            0 & 0 & 0 & \hdots & 0 & 0 & 0
        \end{bmatrix}}_{ \in \mathbb{R}^{(k+1)\times (k+1)}} x
        + \begin{bmatrix}
            \alpha \\ 0 \\ 0 \\ \vdots \\ 0 \\ 0 \\ \|v\|
        \end{bmatrix}
    \]
    for all $x \in \mathbb{R}^{k+1}$ with $\alpha\neq0$.
    % Note that if $\|v\|=0$, $\opT$ has a unique fixed point $x_\star = - \frac{\alpha}{2} (1, 1, \dots, 1)$ due to the invertibility of $M$.
    % Then $\opT$ is a nonexpansive operator with infimal displacement vector $\tilde{v}=(0,\dots,0,\|v\|)$.
    We provide the detailed proof in \cref{sec:prooflowerbound}.
\end{proof}

\paragraph{Matching upper and lower bounds.}
The $\frac{4}{k^2}\|x^0-x_\star\|^2$ upper bound on \eqref{eq:fpi} for $- \frac{x^k-x^0}{k} \to v$ of \cref{cor:optimalkm} exactly matches the $\frac{4}{k^2}\|x^0-x_\star\|^2$ lower bound of \cref{thm:lowerboundspan}.
The upper bounds on \eqref{eq:kmiteration} with $\lambda_k\equiv\lambda\in(0,1)$ and \eqref{eq:halpern} with $\lambda_k=\frac{1}{k+1}$ of \cref{cor:hallinearconv} match the lower bound up to a constant.

The $\mathcal{O}(\frac{\log k}{k^2})$ upper bound of \eqref{eq:halpern} for $x^k-\opT x^k \to v$ matches the lower bound up to logarithmic factors, and this is the fastest known rate for the convergence of the fixed-point residual $x^k-\opT x^k$ to $v$.
However, the $\mathcal{O}(1/k^2)$ upper bound of \eqref{eq:halpern} for $\|x^k-\opT x^k\| \to \|v\|$ does match the lower bound up to a constant.

Finally, the upper bound
\[
    (\|v^k\|-\|v\|)^2\le (\|v^k-v\|)^2\le
\frac{4}{k^2}\|x^0-x_\star\|^2
\]
of \cref{cor:optimalkm} and the $\frac{1}{2k^2}\|x^0-x_\star\|^2$ lower bound of $\|v^k\| \to \|v\|$ of \cref{thm:lowerboundspan} match only up to a constant factor $8$ (where $v^k$ are as defined in  \cref{thm:lowerboundspan} and \cref{cor:optimalkm}).
Reducing this gap may be an interesting direction of future work.

\paragraph{Lower bounds for deterministic iterations.}
Finally, we use the resisting oracle technique of \citet{nemirovski1983problem} to extend the complexity lower bound to general deterministic fixed-point iterations, an algorithm class we formally define in \cref{sec:prooflowerbound}. The following result no longer requires the linear span assumption \eqref{def:span}.

\begin{theorem} \label{thm:lowerboundgeneral}
    Let $k\in\mathbb{N}$, $x^0\in\cH$, and $v\in\cH$, where $\dim\cH \ge 2k-1$.
    Then, there exists a nonexpansive operator $\opT\colon\cH\to\cH$ and $x_\star\in\cH$ such that $v = x_\star - \opT x_\star$, the infimal displacement vector of $\opT$ is $v$, and
    \[
        \left( \left\| \sum_{i=0}^{k-1} \nu_i (x^i-\opT x^i) \right\| - \|v\| \right)^2
        \ge
        \frac{1}{{2} k^2} \|x^0-x_\star\|^2
    \]
    and
    \[
        \left\| \sum_{i=0}^{k-1} \nu_i (x^i-\opT x^i) - v \right\|^2
        \ge
        \frac{4}{k^2} \|x^0-x_\star\|^2
    \]
    hold for iterates $\{x^n\}_{n=0}^{k-1}$ generated by any deterministic fixed-point iteration and any choice of real numbers $\{\nu_i\}_{i=0}^{k-1}$ such that $\sum_i \nu_i = 1$.
\end{theorem}
The proof of \cref{thm:lowerboundgeneral} is deferred to \cref{sec:prooflowerbound}.

\section{Experiments}
\label{sec:experiments}

Consider an infeasible semidefinite problem (SDP)
\[
    \begin{array}{ll}
        \underset{x\in\RR^d}{\mbox{minimize}} &
        \sum_{i=1}^p c_i^\intercal x \\
        \mbox{subject to} &
        \mathcal{A}_i[x] = \sum_{j=1}^{d} A_i^j x_j \preceq B_i, 
        \quad 1\le i\le p,
    \end{array}
\]
where $A_i^j,\, B_i \in \SS^n$ and $\mathcal{A}_i\colon \mathbb{R}^d\rightarrow\mathbb{S}^n$ is a linear operator defined by $\mathcal{A}_i[x] = \sum_{j=1}^{d} A_i^j x_j  $.

Consider a setup where each objective function $c_i^\intercal x$ and $i$-th constraint $\mathcal{A}_i[x] \preceq B_i$ are private to the local agent $i\in\{1,\dots,p\}$.
Assume that they communicate only with their neighbors, which are represented in the graph as connected nodes.
% Each agent knows whether its local constraint is feasible, but does not know whether there exists a joint feasible region.
This SDP can be solved in decentralized manner with PG-EXTRA
% \begin{align*}
%     \!\!U_i^{k+1} &\!= \Pi_{-\SS^n_+} \left( U_i^k + \beta(B_i - \mathcal{A}_i[x^k]) \right) \\
%     % \!\!\vw^{k+1} &\!= \vw^k + \frac{1}{2}(I-W)\vx^k
%     \!\! w^{k+1}_i &\!= w^k_i + \frac{1}{2} \left( x^k_i - {\scriptstyle \sum_{j\in N_i}} W_{i,j} x^k_j \right)
%     \tag{PG-EXTRA}
%     \\
%     \!\!x_i^{k+1} &\!= x_i^k - \alpha\beta (2 w_i^{k+1} - w_i^k)+\alpha \left( 
%     \mathcal{A}_i^*
%     [ 2 U_i^{k+1} - U_i^k] - c_i \right)
% \end{align*}
of \citet{shi2015proximal}.
% Here, $\alpha,\,\beta>0$,  $W\in\RR^{p\times p}$ a mixing matrix with $(i,j)$-entry $W_{i,j}$, $N_i$ the neighbors of agent $i$, and $\mathcal{A}_i^*\colon \mathbb{S}^n\rightarrow\mathbb{R}^d$ the adjoint linear operator of $\mathcal{A}_i$.
% Also, each agent $i$ keeps its local copy of $u_i^k$, $w_i^k$, and $x_i^k$ at iteration $k$.
See \cref{sec:proofexperiments} for the details of infeasible SDP instance, derivation of PG-EXTRA for SDP, and the choices of parameters.

\cref{fig:sdp_experiment} compares the results of PG-EXTRA and PG-EXTRA combined with OHM. Both algorithms' normalized iterates and fixed-point residuals converged to $v$, but OHM is faster for fixed-point residual, as our theory suggests.

\section{Conclusions}
\label{sec:conclusions}

In this work, we analyzed the convergence rates of fixed-point iterations towards the infimal displacement vector. By providing matching upper and lower bounds, we established the optimal accelerated complexity to be $\mathcal{O}(1/k^2)$. The discovery of our upper bounds was assisted by the performance estimation problem (PEP) methodology, which we extended to accommodate inconsistent problem setups.

In our view, the analysis of optimization algorithms applied to inconsistent problems is a necessary step in designing robust general-purpose solvers. Carrying out similar analyses for different algorithms under different inconsistent problems is an interesting direction of future work, and we expect our newly extended PEP methodology to be broadly useful in such endeavors.

% Interestingly, the upper bounds on $\|v^k-v\|^2$ are usually consistently weaker than those of $(\|v^k\|-\|v\|)^2$.
% For example, in \cref{thm:kmfprconv}, the rate on $\|v^k-v\|^2$ is a best-iterate rate, while the rate on $(\|v^k\|-\|v\|)^2$ is the last-iterates rate, which is a stronger type of result.
% As another example, in \cref{thm:halfprconv}, the rate on $\|v^k-v\|^2$ is $\mathcal{O}(\log k/k^2)$, while the rate of $(\|v^k\|-\|v\|)^2$ is $\mathcal{O}(1/k^2)$.

% \section*{Acknowledgements}

% \textbf{Do not} include acknowledgements in the initial version of
% the paper submitted for blind review.

% If a paper is accepted, the final camera-ready version can (and
% probably should) include acknowledgements. In this case, please
% place such acknowledgements in an unnumbered section at the
% end of the paper. Typically, this will include thanks to reviewers
% who gave useful comments, to colleagues who contributed to the ideas,
% and to funding agencies and corporate sponsors that provided financial


\newpage
\bibliography{main}
\bibliographystyle{icml2023}

%%%%%%%%%%%%%%%%%%%%%%%%%%%%%%%%%%%%%%%%%%%%%%%%%%%%%%%%%%%%%%%%%%%%%%%%%%%%%%%
%%%%%%%%%%%%%%%%%%%%%%%%%%%%%%%%%%%%%%%%%%%%%%%%%%%%%%%%%%%%%%%%%%%%%%%%%%%%%%%
% APPENDIX
%%%%%%%%%%%%%%%%%%%%%%%%%%%%%%%%%%%%%%%%%%%%%%%%%%%%%%%%%%%%%%%%%%%%%%%%%%%%%%%
%%%%%%%%%%%%%%%%%%%%%%%%%%%%%%%%%%%%%%%%%%%%%%%%%%%%%%%%%%%%%%%%%%%%%%%%%%%%%%%
\newpage
\appendix
\onecolumn

\section{Omitted proof of \cref{sec:measure}}
\label{sec:proofmeasure}

\begin{proof}[Proof of \cref{lem:normconvequiv}]
    $x^k - \opT x^k \in \overline{\cR(\opI-\opT)}$, so $v^k \in \overline{\cR(\opI-\opT)}$.
    From the property of the projection, as $v^k\in\overline{\cR(\opI-\opT)}$,
    \[
        \langle v^k, v \rangle \ge \|v\|^2, 
        \quad 
        \forall k\in\mathbb{N}.
    \]
    Then we have
    \[
        \|v^k-v\|^2 = \|v^k\|^2 - 2\langle v^k, v \rangle + \|v\|^2 
        \le
        \|v^k\|^2 - \|v\|^2.
    \]
    If $\lim_{k\to\infty} v^k = v$, then obviously, $\lim_{k\to\infty} \|v^k\| = \|v\|$.
    If $\lim_{k\to\infty} \|v^k\| = \|v\|$, then $\lim_{k\to\infty} \|v^k-v\|^2 = 0$ from above inequality, so $\lim_{k\to\infty} v^k = v$.
\end{proof}

\section{Omitted proofs of \cref{sec:convergencerate}}
\label{sec:IDwithresiduals}

\subsection{Omitted proofs of \cref{subsec:km}}
\label{sec:proofkm}

Following lemmas will be used in the proof of \cref{thm:kmaverconv} and \cref{thm:kmfprconv}.
\begin{lemma} \label{lem:kmmonotone}
    If $\{x^k\}_{k\in\mathbb{N}}$ and $\{y^k\}_{k\in\mathbb{N}}$ are sequences of iterates generated by \eqref{eq:kmiteration} starting from $x^0\in\cH$ and $y^0\in\cH$ respectively, for any $k\in\mathbb{N}\cup\{0\}$,
    \[
        \|x^{k+1}-\opT x^{k+1}\| \le \|x^k-\opT x^k\|
    \]
    and
    \[
        \|x^{k+1}-y^{k+1}\| \le \|x^k-y^k\|.
    \]
\end{lemma}

\begin{proof}
    \begin{align*}
        \|x^{k+1}-\opT x^{k+1}\|
        &=
        \|x^{k+1}-\opT x^k+\opT x^k-\opT x^{k+1}\| \\
        &\le
        \|x^{k+1}-\opT x^k\| + \|x^k-x^{k+1}\| \\
        &=
        \lambda_{k+1} \|x^k-\opT x^k\| + (1-\lambda_{k+1}) \|x^k-\opT x^k\| \\
        &=
        \|x^k-\opT x^k\|
    \end{align*}
    and
    \begin{align*}
        \|x^{k+1}-y^{k+1}\|
        &=
        \|(1-\lambda_{k+1})(\opT x^k-\opT y^k) + \lambda_{k+1}(x^k-y^k)\| \\
        &\le
        (1-\lambda_{k+1})\|x^k-y^k\| + \lambda_{k+1}\|x^k-y^k\| \\
        &=
        \|x^k-y^k\|.
    \end{align*}
\end{proof}

\begin{lemma} \label{lem:kmmonotonestar}
    For any $\varepsilon>0$, there exists $x_\varepsilon\in\cH$ such that
    \[
        \|x_\varepsilon - \opT x_\varepsilon - v\| \le \varepsilon.
    \]
    And for any $k\in\mathbb{N}\cup\{0\}$,
    \[
        \|x_\varepsilon^k - \opT x_\varepsilon^k\| - \|v\| \le \varepsilon.
    \]
\end{lemma}
\begin{proof}
    Since $v \in \overline{\cR(\opI-\opT)}$, for any $\varepsilon>0$, we may choose $y_\varepsilon \in \cR(\opI-\opT)$ such that $\|y_\varepsilon - v\| \le \varepsilon$.
    As $y_\varepsilon \in \cR(\opI-\opT)$, there exists $x_\varepsilon\in\cH$ such that $y_\varepsilon = x_\varepsilon-\opT x_\varepsilon$, so
    \[
        \|x_\varepsilon - \opT x_\varepsilon - v\| \le \varepsilon.
    \]

    We know that from \cref{lem:kmmonotone} that for any $k\in\mathbb{N}$,
    \[
        \|x_\varepsilon^k - \opT x_\varepsilon^k\| \le \|x_\varepsilon^{k-1}-\opT x_\varepsilon^{k-1}\|.
    \]
    Therefore,
    \begin{align*}
        \|x_\varepsilon^k - \opT x_\varepsilon^k\| - \|v\|
        &\le 
        \|x_\varepsilon - \opT x_\varepsilon\| - \|v\| \\
        &\le
        \|x_\varepsilon - \opT x_\varepsilon - v\|
        \le \varepsilon.
    \end{align*}
\end{proof}

We now prove our main results of this section.
\begin{proof}[Proof of \cref{thm:kmaverconv}]
    For $\varepsilon>0$, define $\tilde{\varepsilon}$ as
    \[
        \tilde{\varepsilon} = 
        \min \left\{ \frac{\varepsilon^2}{2\|v\|+1}, 1, \varepsilon \right\}
    \]
    and let $x_\varepsilon\in\cH$ be a vector in $\cH$ such that
    \[
        \|x_\varepsilon-\opT x_\varepsilon - v\| \le \tilde{\varepsilon},
    \]
    whose existence is guaranteed from \cref{lem:kmmonotonestar}.
    
    Now let $\{x_\varepsilon^k\}_{k\in\mathbb{N}}$ be a sequence of iterates generated by \eqref{eq:kmiteration} starting from $x_\varepsilon$.
    Expanding the $x^k$ term, we get
    \begin{align*}
        &\frac{x^k-x^0}{\sum_{i=1}^{k} (1-\lambda_i)} + v \\
        &=
        \frac{1}{\sum_{i=1}^{k} (1-\lambda_i)} \left\{
            (x^k-x^k_\varepsilon)
            - (x^0-x_\varepsilon)
            - \left(x_\varepsilon - x^k_\varepsilon - \left(\sum_{i=1}^{k} (1-\lambda_i) \right) v \right)
        \right\} \\
        &=
        \frac{1}{\sum_{i=1}^{k} (1-\lambda_i)} \left\{
            (x^k-x^k_\varepsilon)
            - (x^0-x_\varepsilon)
            - \sum_{i=1}^{k} (1-\lambda_i) \left( x^{i-1}_\varepsilon - \opT x^{i-1}_\varepsilon - v \right)
        \right\}
    \end{align*}
    and taking its norm,
    \begin{align*}
        &\left\|
            \frac{x^k-x^0}{\sum_{i=1}^{k} (1-\lambda_i)} + v
        \right\| \\
        &\le
        \frac{1}{\sum_{i=1}^{k} (1-\lambda_i)} \left(
            \|x^k-x^k_\varepsilon\| + \|x^0-x_\varepsilon\|
        \right)
        + \sum_{i=1}^{k} \frac{(1-\lambda_i)}{\sum_{i=1}^{k} (1-\lambda_i)} \left\| x^{i-1}_\varepsilon - \opT x^{i-1}_\varepsilon - v \right\| \\
        &\le
        \frac{2}{\sum_{i=1}^{k} (1-\lambda_i)} \|x^0-x_\varepsilon\| +  \sum_{i=1}^{k} \frac{(1-\lambda_i)}{\sum_{i=1}^{k} (1-\lambda_i)} \|x^{i-1}_\varepsilon - \opT x^{i-1}_\varepsilon - v\|
        \tag{$\because$ \cref{lem:kmmonotone}}
    \end{align*}
    Since
    \[
        \|x - \opT x - v\|^2
        = \|x-\opT x\|^2 - 2\langle x-\opT x, v \rangle + \|v\|^2
        \le
        \|x-\opT x\|^2 - \|v\|^2,
        \quad \forall x \in \cH,
    \]
    we get
    \begin{align*}
        \|x^{i}_\varepsilon - \opT x^{i}_\varepsilon - v\|^2
        &\le
        \|x^{i}_\varepsilon-\opT x^{i}_\varepsilon\|^2 - \|v\|^2 \\
        &\le
        \|x_\varepsilon - \opT x_\varepsilon\|^2 - \|v\|^2 
        \tag{$\because$ \cref{lem:kmmonotone}} \\
        &=
        \left( \|x_\varepsilon - \opT x_\varepsilon\| - \|v\| \right) \left( \|x_\varepsilon - \opT x_\varepsilon\| + \|v\| \right) \\
        &\le
        \tilde{\varepsilon} (2\|v\| + \tilde{\varepsilon})
        \tag{$\because$ \cref{lem:kmmonotonestar}} \\
        &\le
        \tilde{\varepsilon} (2\|v\| + 1)
        \le \varepsilon^2
    \end{align*}
    for any $i\in\mathbb{N}\cup\{0\}$.
    Gathering all facts above, we get
    \begin{align*}
        \left\|
            \frac{x^k-x^0}{\sum_{i=1}^{k} (1-\lambda_i)} + v
        \right\|
        &\le
        \frac{2}{\sum_{i=1}^{k} (1-\lambda_i)} \|x^0-x_\varepsilon\|
        + \varepsilon
    \end{align*}
    for any $k\in\mathbb{N}$.

    If $v\in\cR(\opI-\opT)$, there exists $x_\star\in\cH$ such that $v=x_\star-\opT x_\star$.
    The proof above applies well with $\varepsilon=0$ and $x_\varepsilon = x_\star$, so we are done.
\end{proof}

According to \cref{thm:kmaverconv}, the normalized iterate of \eqref{eq:kmiteration} converges to $-v$ when $\sum_{i=1}^\infty (1-\lambda_i) = \infty$.

\begin{corollary}
    Let $\{x^k\}_{k\in\mathbb{N}}$ be the iterates of \eqref{eq:kmiteration} starting from $x^0\in\cH$.
    If $\sum_{i=1}^\infty (1-\lambda_i) = \infty$, then
    \[
        \lim_{k\to\infty} \frac{x^k-x^0}{\sum_{i=1}^k (1-\lambda_i)} = - v.
    \]
\end{corollary}

\begin{proof}
    According to the first claim, for any $\varepsilon>0$, there exists $x_\varepsilon\in\cH$ such that
    \[
        \left\|
            \frac{x^k-x^0}{\sum_{i=1}^{k} (1-\lambda_i)} + v
        \right\|
        \le
        \frac{2}{\sum_{i=1}^{k} (1-\lambda_i)} \|x^0-x_\varepsilon\|
        + \varepsilon.
    \]
    Therefore, given $\sum_{i=1}^\infty (1-\lambda_i) = \infty$,
    \[
        0
        \le
        \limsup_{k\to\infty} \left\|
            \frac{x^k-x^0}{\sum_{i=1}^{k} (1-\lambda_i)} + v
        \right\|
        \le
        \varepsilon
    \]
    for any $\varepsilon>0$.
    We may conclude that $\frac{x^k-x^0}{\sum_{i=1}^{k} (1-\lambda_i)}$ converges to $-v$ in norm.
\end{proof}

Convergence of the fixed-point residual $x^k-\opT x^k$ to $v$ requires a stronger assumption, which is $\sum_{k=0}^\infty \lambda_k(1-\lambda_k) = \infty$.
This is a stronger condition than that of \cref{thm:kmaverconv} in a sense that
\[
    \sum_{k=0}^\infty \lambda_i(1-\lambda_i) = \infty
    \quad\implies\quad
    \sum_{k=0}^\infty \lambda_i = \infty.
\]
In case of $\fix\opT \neq \emptyset$, The iterates $\{x^k\}$ generated by \eqref{eq:kmiteration} exhibits Fejer-monotonicity with respect to $\fix\opT$ \citep[Chapter~5]{bauschke2011convex}, which is a useful concept in proving the convergence of \eqref{eq:kmiteration} in terms of $x^k-\opT x^k\to0$ and $x^k \to x_\star$.
However, when $\fix\opT=\emptyset$, such analysis is impossible.
% , and we provide a concept beyond Fejer-monotonicity.

Consider a sequence $\{\lambda_k\}_{k\in\mathbb{N}\cup\{0\}}$ of stepsizes to \eqref{eq:kmiteration}.
Define $\opT_k\colon\cH\to\cH$, for each $k\in\mathbb{N}$ as
\[
    \opT_k := (1-\lambda_{k}) \opT + \lambda_{k} \opI.
\]
Then if $\{x^k\}_{k\in\mathbb{N}}$ is a sequence of iterates generated by \eqref{eq:kmiteration} with $\{\lambda_k\}_{k\in\mathbb{N}\cup\{0\}}$ starting from $x^0\in\cH$,
\[
    x^{k+1} = \opT_{k+1} x^k.
\]

\begin{lemma} \label{lem:kmfejer}
    If $\{x^k\}_{k\in\mathbb{N}}$ and $\{y^k\}_{k\in\mathbb{N}}$ are sequences of iterates generated by \eqref{eq:kmiteration} starting from $x^0\in\cH$ and $y^0\in\cH$ respectively, for any $k\in\mathbb{N}\cup\{0\}$,
    \[
        \|x^k-y^k\|^2 - \|x^{k+1}-y^{k+1}\|^2
        \ge
        \lambda_{k+1}(1-\lambda_{k+1}) \|(x^k-\opT x^k) - (y^k-\opT y^k)\|^2
    \]
\end{lemma}
\begin{proof}
    First of all, if $\lambda_{k+1} = 0$ or $1$, the theorem trivially holds from the fact that $\opT_{k+1}$ is a nonexpansive operator.

    Suppose $\lambda_{k+1}\in(0,1)$.
    \begin{align*}
        &\|(x^k-x^{k+1}) - (y^k-y^{k+1})\|^2 \\
        &=
        \|x^k-y^k\|^2 + \|x^{k+1}-y^{k+1}\|^2 - 2 \langle x^{k+1}-y^{k+1}, x^k-y^k \rangle \\
        &=
        \|x^k-y^k\|^2 + \|x^{k+1}-y^{k+1}\|^2 - 2 \langle \opT_{k+1} x^k - \opT_{k+1} y^k, x^k-y^k \rangle.
    \end{align*}
    From $(1-\lambda_{k+1})$-averagedness of $\opT_{k+1}$, \citep[Proposition~4.35(iv)]{bauschke2011convex} gives us
    \begin{align*}
        \|\opT_{k+1} x^k - \opT_{k+1} y^k\|^2 + (2\lambda_{k+1}-1) \|x^k-y^k\|^2 
        &\le
        2 \lambda_{k+1} \langle \opT_{k+1} x^k - \opT_{k+1} y^k, x^k-y^k \rangle.
    \end{align*}
    Then
    \begin{align*}
        &\lambda_{k+1} \|(x^k-x^{k+1}) - (y^k-y^{k+1})\|^2 \\
        &=
        \lambda_{k+1} \|x^k-y^k\|^2 
        + \lambda_{k+1}  \|x^{k+1}-y^{k+1}\|^2 
        - 2 \lambda_{k+1} \langle \opT_{k+1} x^k - \opT_{k+1} y^k,\, x^k-y^k \rangle \\
        &\le
        \lambda_{k+1} \|x^k-y^k\|^2 
        + \lambda_{k+1}  \|x^{k+1}-y^{k+1}\|^2 
        - \left\{
            (2\lambda_{k+1}-1) \|x^k-y^k\|^2 + \|x^{k+1}-y^{k+1}\|^2
        \right\} \\
        &=
        (1-\lambda_{k+1}) \left\{ \|x^k-y^k\|^2 - \|x^{k+1}-y^{k+1}\|^2 \right\}.
    \end{align*}
    As
    \begin{align*}
        x^k-x^{k+1}
        &=
        x^k - \{ (1 - \lambda_{k+1}) \opT x^k + \lambda_{k+1} x^k \}
        =
        (1 - \lambda_{k+1}) (x^k-\opT x^k),
    \end{align*}
    combining this fact with above inequality and dividing by $1-\lambda_{k+1}>0$.
    \begin{align*}
        \|x^k-y^k\|^2 - \|x^{k+1}-y^{k+1}\|^2
        &\ge
        \lambda_{k+1} (1-\lambda_{k+1}) \|(x^k-\opT x^k) - (y^k-\opT y^k)\|^2.
    \end{align*}
\end{proof}

We now prove the second main result of this section.
\begin{proof}[Proof of \cref{thm:kmfprconv}]
    Given $\varepsilon>0$, there exists $x_\varepsilon\in\cH$ such that
    \[
        \|x_\varepsilon-\opT x_\varepsilon-v\| \le \tilde{\varepsilon} = \min \left\{ \frac{\varepsilon^2}{2\|v\|+1}, 1, \varepsilon \right\}
    \]
    by \cref{lem:kmmonotonestar}.
    Let $\{x_\varepsilon^k\}_{k\in\mathbb{N}}$ be a sequence of iterates generated by \eqref{eq:kmiteration} starting from $x_\varepsilon$.
    With $y^0=x_\varepsilon$, summing up the inequality in \cref{lem:kmfejer} and removing the telescoping terms, we get
    \[
        \|x^0-x_\varepsilon\|^2 - \|x^{k+1}-x_\varepsilon^{k+1}\|^2
        \ge
        \sum_{i=0}^{k} \lambda_{i+1} (1-\lambda_{i+1}) \|(x^{i}-\opT x^{i}) - (x_\varepsilon^{i}-\opT x_\varepsilon^{i})\|^2
    \]
    for any $k\in\mathbb{N}$.
    Therefore,
    \begin{align*}
        &\frac{1}{\sum_{i=0}^{k} \lambda_{i+1} (1-\lambda_{i+1})} \|x^0 - x_\varepsilon\|^2 \\
        &\ge
        \sum_{i=0}^{k} \left( \frac{\lambda_{i+1} (1-\lambda_{i+1})}{\sum_{i=0}^{k} \lambda_{i+1} (1-\lambda_{i+1})} \right) \left\| (x^{i} - \opT x^{i}) - (x_\varepsilon^{i} - \opT x_\varepsilon^{i}) \right\|^2 \\
        &=
        \left\{ \sum_{i=0}^{k} \left( \frac{\lambda_{i+1}(1-\lambda_{i+1})}{\sum_{i=0}^{k} \lambda_{i+1} (1-\lambda_{i+1})} \right) \right\}
        \left\{ \sum_{i=0}^{k} \left( \frac{\lambda_{i+1} (1-\lambda_{i+1})}{\sum_{i=0}^{k} \lambda_{i+1} (1-\lambda_{i+1})} \right) \left\| (x^{i} - \opT x^{i}) - (x_\varepsilon^{i} - \opT x_\varepsilon^{i}) \right\|^2 \right\} \\
        &\ge
        \left\{ \sum_{i=0}^{k} \left( \frac{\lambda_{i+1}(1-\lambda_{i+1})}{\sum_{i=0}^{k} \lambda_{i+1} (1-\lambda_{i+1})} \right) \left\| (x^{i} - \opT x^{i}) - (x_\varepsilon^{i} - \opT x_\varepsilon^{i}) \right\| \right\}^2
        \tag{Cauchy-Schwarz}
    \end{align*}
    or equivalently,
    \[
        \sum_{i=0}^{k}
        \left( \frac{\lambda_{i+1} (1-\lambda_{i+1})}{\sum_{i=0}^{k} \lambda_{i+1} (1-\lambda_{i+1})} \right)
        \left\| (x^{i} - \opT x^{i}) - (x_\varepsilon^{i} - \opT x_\varepsilon^{i}) \right\|
        \le
        \frac{1}{\sqrt{ \sum_{i=0}^{k}
        \lambda_{i+1} (1-\lambda_{i+1}) }}
        \|x^0 - x_\varepsilon\|.
    \]
    Note that for any $x\in\cH$,
    \[
        \|x-\opT x-v\|^2 = \|x-\opT x\|^2 - 2\underbrace{\langle x-\opT x, v \rangle}_{\ge \|v\|^2} + \|v\|^2
        \le \|x-\opT x\|^2 - \|v\|^2.
    \]
    For any $i\in\mathbb{N}$,
    \begin{align*}
        \|x_\varepsilon^{i} - \opT x_\varepsilon^{i} - v\|^2
        &\le
        (\|x_\varepsilon^{i}-\opT x_\varepsilon^{i}\| - \|v\|) (\|x_\varepsilon^{i}-\opT x_\varepsilon^{i}\| + \|v\|) \\
        &\le
        (\|x_\varepsilon-\opT x_\varepsilon\| - \|v\|) (\|x_\varepsilon-\opT x_\varepsilon\| + \|v\|) \\
        &\le
        \tilde{\varepsilon} (2\|v\| + \tilde{\varepsilon}) \le \varepsilon^2,
    \end{align*}
    so
    \begin{align*}
        \left\| (x^{i} - \opT x^{i}) - (x_\varepsilon^{i} - \opT x_\varepsilon^{i}) \right\|
        &\le
        \left\| x^{i} - \opT x^{i} - v \right\| - \left\| x_\varepsilon^{i} - \opT x_\varepsilon^{i} - v \right\| \\
        &\le
        \left\| x^{i} - \opT x^{i} - v \right\| - \varepsilon.
    \end{align*}
    Therefore, we get
    \[
        \sum_{i=0}^{k}
        \left( \frac{\lambda_{i+1} (1-\lambda_{i+1})}{\sum_{i=0}^{k} \lambda_{i+1} (1-\lambda_{i+1})} \right)
        \left\| x^{i} - \opT x^{i} - v \right\|
        \le
        \frac{1}{\sqrt{ \sum_{i=0}^{k}
        \lambda_{i+1} (1-\lambda_{i+1}) }}
        \|x^0 - x_\varepsilon\| + \varepsilon.
    \]
    Also, note that for any $i$ such that $0 \le i \le k-1$,
    \[
        \left\| x^{i} - \opT x^{i} - v \right\|
        \ge
        \left\| x^{i} - \opT x^{i} \right\| - \|v\| 
        \ge
        \left\| x^k - \opT x^k \right\| - \|v\|
    \]
    where the first inequality comes from triangular inequality, and the last inequality comes from \cref{lem:kmmonotone}.
    Hence we get
    \[
        \|x^k - \opT x^k\| - \|v\|
        \le 
        \frac{1}{\sqrt{\sum_{i=0}^{k} \lambda_{i+1} (1-\lambda_{i+1})}} \|x^0 - x_\varepsilon\| + \varepsilon.
    \]

    If $v\in\cR(\opI-\opT)$, there exists $x_\star\in\cH$ such that $v=x_\star-\opT x_\star$.
    The proof above applies well with $\varepsilon=0$ and $x_\varepsilon = x_\star$, so we are done.
\end{proof}

According to \cref{thm:kmfprconv}, the fixed-point residual of \eqref{eq:kmiteration} converges to $v$ if $\sum_{i=1}^\infty \lambda_{i} (1-\lambda_{i}) = \infty$.
\begin{corollary}
    Let $\{x^k\}_{k\in\mathbb{N}}$ be the iterates of \eqref{eq:kmiteration} starting from $x^0\in\cH$.
    If $\sum_{i=1}^\infty \lambda_i (1-\lambda_i) = \infty$, then
    \[
        \lim_{k\to\infty} \left( x^k - \opT x^k \right) = v.
    \]
\end{corollary}

\begin{proof}
    Given $\sum_{i=1}^\infty \lambda_i(1-\lambda_i) = \infty$,
    \[
        0 \le \limsup_{k\to\infty} \|x^k-\opT x^k\| - \|v\| \le \varepsilon.
    \]
    Since above inequality holds for any choice of $\varepsilon>0$, $\lim_{k\to\infty} \|x^k-\opT x^k\| = \|v\|$.
    Since
    \[
        \|x^k-\opT x^k-v\|^2 \le \|x^k-\opT x^k\|^2 - \|v\|^2,
    \]
    taking limit on both sides, we get
    \[
        0 \le \limsup_{k\to\infty} \|x^k-\opT x^k-v\|^2 \le \lim_{k\to\infty} \|x^k-\opT x^k\|^2 - \|v\|^2 = \|v\|^2 - \|v\|^2 = 0.
    \]
\end{proof}

\subsection{Omitted proofs of \cref{subsec:halpern}}
\label{sec:proofhalpern}

Following lemmas will be used in the proof of \cref{thm:halaverconv} and \cref{thm:halfprconv}.

We first prove \cref{lem:haltheta}.
\begin{proof}[Proof of \cref{lem:haltheta}]
    If $k=0$, then
    \[
        \theta_1 = (1-\lambda_1) = (1-\lambda_1)(1+\underbrace{\theta_0}_{=0}).
    \]
    Suppose $k\ge1$.
    \begin{align*}
        \theta_{k+1}
        &=
        \sum_{n=1}^{k+1} (1-\lambda_{k+1})(1-\lambda_{k})\cdots(1-\lambda_{k-n+2}) \\
        &=
        (1-\lambda_{k+1})
        + (1-\lambda_{k+1}) \sum_{n=2}^{k+1} (1-\lambda_{k}) \cdots (1-\lambda_{k-n+2} \\
        &=
        (1-\lambda_{k+1}) 
        + (1-\lambda_{k+1}) \sum_{n=1}^k (1-\lambda_k) \cdots (1-\lambda_{k-n+1} \\
        &=
        (1-\lambda_{k+1}) (1+\theta_k).
    \end{align*}

    Suppose $\lambda_k\equiv 0$. Then
    \[
        \theta_k = 1 + \theta_{k-1} = 2 + \theta_{k-2} = \cdots = k + \theta_0 = k.
    \]

    If $\lambda_k = \frac{1}{k+1}$ for all $k\in\mathbb{N}$, then from $\theta_0 = 0$, suppose $\theta_{k-1} = \frac{k-1}{2}$.
    Then as
    \[
        \theta_k = \left(1 - \frac{1}{k+1} \right) (1+\theta_{k-1}) = \frac{k}{k+1} \frac{k+1}{2} = \frac{k}{2},
    \]
    the induction holds.
\end{proof}

\begin{remark}
    Let $\{x^k\}_{k\in\mathbb{N}}$ be the iterates of \eqref{eq:halpern} starting from $x^0\in\cH$.
    Then the $k$-th iterate $x^k$ of \eqref{eq:halpern} can be expressed as
    \[
        x^k - x^0 = - \sum_{i=0}^{k-1} \left\{ (1-\lambda_{k}) \cdots (1-\lambda_{i+1}) \right\} (x^i-\opT x^i).
    \]
    If $\lambda_k=\frac{1}{k+1}$ for $k\in\mathbb{N}$, the $k$-th iterate $x^k$ of \eqref{eq:halpern} can be expressed as
    \[
        x^k - x^0
        =
        - \sum_{i=0}^{k-1} \frac{i+1}{k+1} (x^i-\opT x^i)
    \]
    The sequence $\{\theta_k\}$ refers to the sum of all linear coefficients to $\{x^i-\opT x^i\}_{i=0,1,\ldots,k-1}$ used in the $x^k$-update of \eqref{eq:halpern}.
\end{remark}

Following lemma refers to the property that two independent iterates $\{x^k\}$ and $\{y^k\}$ generated by \eqref{eq:halpern} cannot be further than the distance between initial points $\|x^0-y^0\|$.
\begin{lemma}
    \label{lem:halmono}
    Let $\{x^k\}_{k\in\mathbb{N}}$ and $\{y^k\}_{k\in\mathbb{N}}$ be a sequence of iterates generated by \eqref{eq:halpern} starting from $x^0\in\cH$ and $y^0\in\cH$, respectively.
    Then
    \[
        \|x^k - y^k\| \le \|x^0 - y^0\|, \quad k=0,1,\dots.
    \]
\end{lemma}
\begin{proof}
    We prove by induction on $k$.
    If $k=0$, the claim automatically holds.
    Suppose $k\ge1$ and $\|x^{k-1}-y^{k-1}\| \le \|x^0-y^0\|$.
    Then
    \begin{align*}
        \|x^{k} - y^{k}\|
        &\le
        \left( 1 - \lambda_{k} \right) \|\opT x^{k-1} - \opT y^{k-1}\| + \lambda_{k} \|x^0 - y^0\| \\
        &\le
        \left( 1 - \lambda_{k} \right) \|x^{k-1} - y^{k-1}\| + \lambda_{k} \|x^0-y^0\| \\
        &\le
        \left( 1 - \lambda_{k} \right) \|x^0 - y^0\| + \lambda_{k} \|x^0-y^0\|
        =
        \|x^0 - y^0\|.
    \end{align*}
\end{proof}

\begin{lemma}
    \label{lem:halaverupperbound}
    If $\{x^k\}_{k\in\mathbb{N}}$ be a sequence of iterates generated by \eqref{eq:halpern} starting from $x^0\in\cH$, then
    \[
        \frac{\|x^k - x^0\|}{\theta_k} \le \|x^0-\opT x^0\|, \quad k=1,2,\dots.
    \]
\end{lemma}
\begin{proof}
    We prove by induction on $k$.
    \begin{itemize}
        \item[(i)] $k=1$.
        First of all,
        \[
            x^1 - x^0 = - (1-\lambda_1) (x^0-\opT x^0)
        \]
        so from $\theta_1 = 1 - \lambda_1$,
        \[
            \frac{\|x^1-x^0\|}{\theta_1} = \|x^0-\opT x^0\|.
        \]

        \item[(ii)] $k\ge2$.
        Suppose that the claim holds true for all $n$ such that $n<k$.
        \begin{align*}
            x^k-x^0 
            &=
            \left( 1 - \lambda_{k} \right) (\opT x^{k-1}-x^0) \\
            &=
            \left( 1 - \lambda_{k} \right) (\opT x^{k-1}-\opT x^0) + \left( 1 - \lambda_{k} \right) (\opT x^0-x^0) \\
            \|x^k-x^0\|
            &\le
            \left( 1 - \lambda_{k} \right) \|\opT x^{k-1}-\opT x^0\| + \left( 1 - \lambda_{k} \right) \|x^0-\opT x^0\| \\
            &\le
            \left( 1 - \lambda_{k} \right) \|x^{k-1}-x^0\| + \left( 1 - \lambda_{k} \right) \|x^0-\opT x^0\|
        \end{align*}
        Therefore,
        \begin{align*}
            \frac{\|x^k-x^0\|}{\theta_k}
            &\le
            \frac{1 - \lambda_{k}}{\theta_k} \|x^{k-1}-x^0\|
            + \frac{1 - \lambda_{k}}{\theta_k} \|x^0-\opT x^0\| \\
            &=
            \frac{\theta_{k-1}}{1+\theta_{k-1}} \frac{\|x^{k-1}-x^0\|}{\theta_{k-1}}
            + \frac{1}{1+\theta_{k-1}} \|x^0-\opT x^0\|
            \tag{$\because$ \cref{lem:haltheta}} \\
            &\le
            \frac{\theta_{k-1}}{1+\theta_{k-1}} \|x^0-\opT x^0\|
            + \frac{1}{1+\theta_{k-1}} \|x^0-\opT x^0\| \\
            &=
            \|x^0-\opT x^0\|.
        \end{align*}
    \end{itemize}
\end{proof}

Following lemma identifies the proper averaging of $x^k$ that resides in the closure of the range of $\opI-\opT$, which becomes the candidate for the sequence $\{v^k\}_{k\in\mathbb{N}}$ converging to $v$.
\begin{lemma} \label{lem:halaverspan}
     If $\{x^k\}_{k\in\mathbb{N}}$ be a sequence of iterates generated by \eqref{eq:halpern} starting from $x^0\in\cH$, then
    \[
        - \frac{x^k-x^0}{\theta_k} \in \overline{\cR(\opI-\opT)}
    \]
    for $k=1,2,\dots$.
\end{lemma}
\begin{proof}
    We prove by induction on $k$, using the convexity of $\overline{\cR(\opI-\opT)}$.
    \begin{itemize}
        \item[($i$)] $k=1$.
        \[
            - \frac{x^1-x^0}{\theta_1} 
            = x^0-\opT x^0
            \in \overline{\mathcal{R}(I-T)}.
        \]
        
        \item[($ii$)] $k\ge2$.
        Suppose that
        \[
            - \frac{x^{k-1}-x^0}{\theta_{k-1}} \in \overline{\cR(\opI-\opT)}.
        \]
        As
        \begin{align*}
            - \frac{x^k-x^0}{\theta_k}
            &=
            - \frac{(1-\lambda_k) \theta_{k-1}}{\theta_k} \frac{x^{k-1}-x^0}{\theta_{k-1}}
            + \frac{1 - \lambda_{k}}{\theta_k} (x^{k-1}-\opT x^{k-1}) \\
            &=
            \frac{\theta_{k-1}}{1+\theta_{k-1}} \left( - \frac{x^{k-1}-x^0}{\theta_{k-1}} \right)
            + \frac{1}{1+\theta_{k-1}} (x^{k-1}-\opT x^{k-1}),
        \end{align*}
        $-\frac{x^k-x^0}{\theta_k}$ is a convex combination of vectors in a convex set $\overline{\cR(\opI-\opT)}$, so it is also an element of $\overline{\cR(\opI-\opT)}$.
    \end{itemize}
\end{proof}

We now prove \cref{thm:halaverconv}.
\begin{proof}[Proof of \cref{thm:halaverconv}]
    From $v \in \overline{\cR(\opI-\opT)}$, we may choose a point $x_\varepsilon$ in $\cH$ such that
    \[
        \|x_\varepsilon - \opT x_\varepsilon\|^2 - \|v\|^2
        \le
        {\varepsilon^2}.
    \]
    Let $k\ge1$.
    From \cref{lem:halmono},
    \begin{align*}
        \left\| \frac{x^k-x^0}{\theta_k} - \frac{x^k_\varepsilon - x_\varepsilon}{\theta_k} \right\|
        &\le
        \left\| \frac{x^k - x^k_\varepsilon}{\theta_k} \right\|
        + \left\| \frac{x^0 - x_\varepsilon}{\theta_k} \right\| \\
        &\le
        \frac{2}{\theta_k} \|x^0-x_\varepsilon\|.
    \end{align*}
    % and 
    % \[
    %     \left\| \frac{x_\varepsilon^k - x_\varepsilon}{\theta_k} \right\|
    %     \le
    %     \| x_\varepsilon - \opT x_\varepsilon \|
    %     \le
    %     \| x_\varepsilon - \opT x_\varepsilon - v \| + \| v \|
    %     \le
    %     \varepsilon + \| v \|,
    % \]
    % we have
    % \begin{align*}
    %     \left\| \frac{x^k-x^0}{\theta_k} \right\| - \|v\|
    %     &\le
    %     \left\| \frac{x^k-x^0}{\theta_k} - \frac{x^k_\varepsilon-x_\varepsilon}{\theta_k} \right\|
    %     + \left\| \frac{x^k_\varepsilon-x_\varepsilon}{\theta_k} \right\| - \|v\| \\
    %     &\le
    %     \frac{2}{\theta_k} \|x^0-x_\varepsilon\| + \|x_\varepsilon-\opT x_\varepsilon\| - \|v\| \\
    %     &\le
    %     \frac{2}{\theta_k} \|x^0-x_\varepsilon\|
    %     + \|x_\varepsilon-\opT x_\varepsilon - v\| \\
    %     &\le
    %     \frac{2}{\theta_k} \|x^0-x_\varepsilon\|
    %     + \varepsilon.
    % \end{align*}
    Note that
    \begin{align*}
        \left\| \frac{x_\varepsilon^k - x_\varepsilon}{\theta_k} + v \right\|^2
        &\le
        \left\| \frac{x_\varepsilon^k - x_\varepsilon}{\theta_k} \right\|^2 - \|v\|^2 \\
        &\le
        \|x_\varepsilon - \opT x_\varepsilon\|^2 - \|v\|^2
        \le
        \varepsilon^2,
        \tag{$\because$ \cref{lem:halaverupperbound}}
    \end{align*}
    and from this we have
    \begin{align*}
        \left\| \frac{x^k-x^0}{\theta_k} + v \right\|
        &\le
        \left\| \frac{x^k-x^0}{\theta_k} - \frac{x^k_\varepsilon-x_\varepsilon}{\theta_k} \right\|
        + \left\| \frac{x^k_\varepsilon-x_\varepsilon}{\theta_k} + v \right\| \\
        &\le
        \frac{2}{\theta_k} \|x^0-x_\varepsilon\| + \varepsilon.
    \end{align*}
    This result holds for any $k\ge1$.

    If $v\in\cR(\opI-\opT)$, there exists $x_\star\in\cH$ such that $v=x_\star-\opT x_\star$.
    The proof above applies well with $\varepsilon=0$ and $x_\varepsilon = x_\star$, so we are done.
\end{proof}

According to \cref{thm:halaverconv}, the normalized iterate of \eqref{eq:halpern} converges to $-v$ when $\lim_{k\to\infty} \theta_k = \infty$.
\begin{corollary}
    Let $\{x^k\}_{k\in\mathbb{N}}$ be the iterates of \eqref{eq:halpern} starting form $x^0\in\cH$.
    If $\lim_{k\to\infty} \theta_k = \infty$,
    then
    \[
        \lim_{k\to\infty} \frac{x^k-x^0}{\theta_k} = -v.
    \]
\end{corollary}

\begin{proof}
    Further assume that $\lim_{k\to\infty} \theta_k = \infty$.
    Using triangle inequality,
    \[
        \left\| \frac{x^k-x^0}{\theta_k} \right\| - \|v\|
        \le
        \left\| \frac{x^k-x^0}{\theta_k} + v \right\|
        \le
        \frac{2}{\theta_k} \|x^0-x_\varepsilon\| + \varepsilon.
    \]
    From the fact that $- \frac{x^k-x^0}{\theta_k} \in \overline{\cR(\opI-\opT)}$ by \cref{lem:halaverspan} and the fact that $v$ is the minimum norm element in $\overline{\cR(\opI-\opT)}$,
    \[
        \left\| \frac{x^k-x^0}{\theta_k} \right\| \ge \|v\|.
    \]
    Then
    \begin{align*}
        \|v\|
        \le \liminf_{k\to\infty} \left\| \frac{x^k-x^0}{\theta_k} \right\|
        \le \limsup_{k\to\infty} \left\| \frac{x^k-x^0}{\theta_k} \right\|
        \le \|v\| + \varepsilon
    \end{align*}
    holds for any possible choice of $\varepsilon>0$, so
    \[
        \lim_{k\to\infty} \left\| \frac{x^k-x^0}{\theta_k} \right\| = \|v\|.
    \]
    We may conclude that, by the uniqueness of $v$ as a minimum norm element in $\overline{\cR(\opI-\opT)}$,
    \[
        \lim_{k\to\infty} \frac{x^k-x^0}{\theta_k} = - v.
    \]
\end{proof}

As in \cref{subsec:halpern}, we have a simpler condition for $\{\lambda_k\}_{k\in\mathbb{N}}$ to ensure the convergence of normalized iterate of Halpern iteration to $-v$.
\begin{lemma} \label{lem:halpernlambda}
    If
    \[
        \lim_{k\to\infty} \lambda_k = 0,
    \]
    then
    \[
        \lim_{k\to\infty} \theta_k = \infty.
    \]
\end{lemma}

\begin{proof}
    \[
        \lim_{k\to\infty} \frac{\theta_{k+1}}{1+\theta_k}
        = \lim_{k\to\infty} (1-\lambda_{k+1})
        = 1,
    \]
    and $\lambda_{k+1} \in [0,1]$, so for any $0<\varepsilon<1$, there exists $N_\varepsilon\in\mathbb{N}$ such that
    \[
        \frac{\theta_{k+1}}{1+\theta_k} \ge 1-\varepsilon, \quad \forall\, k\ge N_\varepsilon.
    \]
    Then
    \begin{align*}
        \theta_{k+N_\varepsilon}
        &\ge (1-\varepsilon) \theta_{k+N_\varepsilon-1} + (1-\varepsilon) \\
        &\ge (1-\varepsilon)^k \theta_{N_\varepsilon} + (1-\varepsilon) + \cdots + (1-\varepsilon)^k \\
        &= (1-\varepsilon)^k \theta_{N_\varepsilon} + \left(\frac{1}{\varepsilon} - 1 \right) \left\{ 1 - (1-\varepsilon)^k \right\}.
    \end{align*}
    As $k\to\infty$,
    \[
        \liminf_{k\to\infty} \theta_k \ge \frac{1}{\varepsilon} - 1
    \]
    holds for all $\varepsilon\in(0,1)$.
    As $\varepsilon\to 0$, $\liminf_{k\to\infty} \theta_k = \infty$, so we are done.
\end{proof}

% Halpern iteration for fixed-point residual
In order to prove \cref{thm:halfprconv}, we use the following fact to construct Lyapunov function.
\begin{lemma} \label{lem:Uk-mono}
    If $\{x^k\}_{k\in\mathbb{N}}$ is a sequence of iterates generated by \eqref{eq:halpern} starting from $x^0\in\cH$ with $\lambda_k = \frac{1}{k+1}$, then
    \begin{align*}
        &\|\opT x^k - \opT x^{k+1}\|^2 \le \|x^k-x^{k+1}\|^2 \\
        &~\Leftrightarrow~
        (k+2)\left\{
            (k+1)\|x^{k+1}-\opT x^{k+1}\|^2 + 2\langle x^{k+1}-\opT x^{k+1},\, x^{k+1}-x^0 \rangle
        \right\} \\
        &\qquad
        \le
        (k+1)\left\{
            k\|x^k-\opT x^k\|^2 + 2\langle x^k-\opT x^k,\, x^k-x^0 \rangle
        \right\}
    \end{align*}
\end{lemma}

\begin{proof}
    From
    \[
        x^{k+1} = \frac{k+1}{k+2}\opT x^k + \frac{1}{k+2} x^0, \quad k=0,1,\ldots,
    \]
    we have
    \begin{align*}
        &\|x^k-x^{k+1}\|^2 - \|\opT x^k-\opT x^{k+1}\|^2 \\
        &=
        \|(x^k - \opT x^k) - (x^{k+1}-\opT x^k)\|^2 - \|(x^{k+1}-\opT x^{k+1}) - (x^{k+1}-\opT x^k)\|^2 \\
        &=
        \|x^k-\opT x^k\|^2 - \|x^{k+1}-\opT x^{k+1}\|^2 - 2 \langle x^k-\opT x^k,\, x^{k+1}-\opT x^k \rangle + 2\langle x^{k+1}-\opT x^{k+1},\, x^{k+1}-\opT x^k \rangle \\
        &=
        \|x^k-\opT x^k\|^2 - \|x^{k+1}-\opT x^{k+1}\|^2
        - 2 \left\langle x^k-\opT x^k,\, \frac{1}{k+2}(x^k-\opT x^k) - \frac{1}{k+2}(x^k-x^0) \right\rangle \\
        &\quad
        + 2 \left\langle x^{k+1}-\opT x^{k+1},\, (x^{k+1}-x^0) - \frac{k+2}{k+1} (x^{k+1}-x^0) \right\rangle \\
        &=
        \frac{1}{k+2} \left\{ k \|x^k-\opT x^k\|^2 + 2 \langle x^k-\opT x^k,\, x^k-x^0 \rangle \right\} \\
        &\quad
        - \frac{1}{k+1} \left\{ (k+1) \|x^{k+1}-\opT x^{k+1}\|^2 + 2 \langle x^{k+1}-\opT x^{k+1},\, x^{k+1}-x^0 \rangle \right\}
    \end{align*}
    Equivalence follows immediately.
\end{proof}

We use the Lyapunov function $V^k$ for $k=1,2,\ldots$ of the following form.
\begin{align*}
    V^k
    &=
    (k+1)\left\{ k\|x^k-\opT x^k\|^2 + 2\langle x^k-\opT x^k, x^k-x^0 \rangle \right\}
     - \left( \sum_{n=1}^{k} \frac{1}{n}\right) \|x^0-x_\varepsilon\|^2 \\
    &\quad
    + k(k+1) \left\langle - \frac{2}{k}(x^k-x^0) - \left(x_\varepsilon-\opT x_\varepsilon\right), x_\varepsilon-\opT x_\varepsilon \right\rangle
    + \frac{2(k+1)}{k} \left\|x^k-x_\varepsilon + \frac{k}{2} \left(x_\varepsilon-\opT x_\varepsilon\right) \right\|^2
    \tag{Lyapunov function}
    \label{eq:lyapunov}
\end{align*}
$x_\varepsilon\in\cH$ is chosen to be the point which makes $x_\varepsilon-\opT x_\varepsilon$ very close to $v$.
In particular, if $v\in\cR(\opI-\opT)$, choose $x_\varepsilon$ such that $v = x_\varepsilon-\opT x_\varepsilon$.

Now, we show the monotonicity of $\{V^k\}_k$ in $k$.
\begin{lemma} \label{lem:lyap-mono}
    Let $\{x^k\}_{k\in\mathbb{N}}$ be a sequence of iterates generated by \eqref{eq:halpern} starting from $x^0\in\cH$ with $\lambda_k = \frac{1}{k+1}$, and define $\{V^k\}_{k\in\mathbb{N}\cup\{0\}}$ as \eqref{eq:lyapunov}.
    For any $k\in\mathbb{N}$,
    \[
        V^k \ge V^{k+1}.
    \]
\end{lemma}
\begin{proof}
    From \cref{lem:Uk-mono},
    \begin{align*}
        &V^k - V^{k+1} \\
        &\ge
        \frac{1}{k+1} \|x^0-x_\varepsilon\|^2
        + k(k+1) \left\langle -\frac{2}{k}(x^k-x^0) - (x_\varepsilon-\opT x_\varepsilon),\, x_\varepsilon-\opT x_\varepsilon \right\rangle \\
        &\quad
        - (k+1)(k+2) \left\langle - \frac{2}{k+1}(x^{k+1}-x^0) - (x_\varepsilon-\opT x_\varepsilon),\, x_\varepsilon-\opT x_\varepsilon \right\rangle \\
        &\quad
        + \frac{2(k+1)}{k} \left\| x^k-x_\varepsilon + \frac{k}{2} (x_\varepsilon-\opT x_\varepsilon) \right\|^2
        - \frac{2(k+2)}{k+1} \left\| x^{k+1}-x_\varepsilon + \frac{k+1}{2} (x_\varepsilon-\opT x_\varepsilon) \right\|^2 \\
        &=
        \frac{1}{k+1} \|x^0-x_\varepsilon\|^2
        + \left\{ -k(k+1) + (k+1)(k+2) + \frac{k(k+1)}{2} - \frac{(k+1)(k+2)}{2} \right\} \|x_\varepsilon-\opT x_\varepsilon\|^2 \\
        &\quad
        + \left\langle
            x_\varepsilon-\opT x_\varepsilon,\, -2(k+1)(x^k-x^0) + 2(k+2)(x^{k+1}-x^0) + 2(k+1) (x^k-x_\varepsilon) - 2(k+2) (x^{k+1}-x_\varepsilon)
        \right\rangle \\
        &\quad
        + \frac{2(k+1)}{k}\|x^k-x_\varepsilon\|^2 - \frac{2(k+2)}{(k+1)}\|x^{k+1}-x_\varepsilon\|^2 \\
        &=
        \frac{1}{k+1} \|x^0-x_\varepsilon\|^2
        + (k+1)\|x_\varepsilon-\opT x_\varepsilon\|^2
        - 2\left\langle x_\varepsilon-\opT x_\varepsilon,\, x^0-x_\varepsilon \right\rangle \\
        &\quad
        + \frac{2(k+1)}{k}\|x^k-x_\varepsilon\|^2 - \frac{2(k+2)}{(k+1)}\|x^{k+1}-x_\varepsilon\|^2.
    \end{align*}
    Using
    \[
        \|x^k-x_\varepsilon\|^2 \ge \|\opT x^k - \opT x_\varepsilon\|^2,
    \]
    we get
    \begin{align*}
        &V^k - V^{k+1} \\
        &\ge
        \frac{1}{k+1} \|x^0-x_\varepsilon\|^2
        + (k+1)\|x_\varepsilon-\opT x_\varepsilon\|^2
        - 2\left\langle x_\varepsilon-\opT x_\varepsilon,\, x^0-x_\varepsilon \right\rangle \\
        &\quad
        + \frac{2(k+1)}{k}\|\underbrace{\opT x^k - \opT x_\varepsilon}_{(\opT x^k-x_\varepsilon) + (x_\varepsilon-\opT x_\varepsilon)}\|^2
        - \frac{2(k+2)}{k+1}\|\underbrace{x^{k+1}-x_\varepsilon}_{=\frac{k+1}{k+2}(\opT x^k-x_\varepsilon) + \frac{1}{k+2}(x^0-x_\varepsilon)}\|^2 \\
        &=
        \frac{k}{(k+1)(k+2)}\|x^0-x_\varepsilon\|^2 + \frac{(k+1)(k+2)}{k}\|x_\varepsilon-\opT x_\varepsilon\|^2 - 2\langle x_\varepsilon-\opT x_\varepsilon, x^0-x_\varepsilon\rangle \\
        &\quad
        + \frac{4(k+1)}{k(k+2)}\|\opT x^k-x_\varepsilon\|^2 + \frac{4(k+1)}{k} \langle x_\varepsilon-\opT x_\varepsilon, \opT x^k-x_\varepsilon \rangle - \frac{4}{k+2} \langle \opT x^k-x_\varepsilon, x^0-x_\varepsilon \rangle \\
        &=
        \frac{1}{k(k+1)(k+2)}
        \left\|
            2(k+1)(\opT x^k-x_\varepsilon) - k(x^0-x_\varepsilon) + (k+1)(k+2)(x_\varepsilon-\opT x_\varepsilon)
        \right\|^2 \\
        &\ge 0.
    \end{align*}
\end{proof}

\begin{lemma} \label{lem:lyap-last}
    Let $\{x^k\}_{k\in\mathbb{N}}$ be a sequence of iterates generated by \eqref{eq:halpern} starting from $x^0\in\cH$ with $\lambda_k = \frac{1}{k+1}$ and $\{V^k\}_{k\in\mathbb{N}\cup\{0\}}$ be defined as \eqref{eq:lyapunov}.
    For $k\ge 1$,
    \begin{align*}
        V^{k}
        &\ge
        (k+1)^2 \|(x^k-\opT x^k) - (x_\varepsilon-\opT x_\varepsilon)\|^2
        + 2k(k+1) \langle (x^k-\opT x^k) - (x_\varepsilon-\opT x_\varepsilon),\, x_\varepsilon-\opT x_\varepsilon \rangle \\
        &\quad
        - 2(k+1) \langle (x^k-\opT x^k) - (x_\varepsilon-\opT x_\varepsilon),\, x^0-x_\varepsilon \rangle
        - \left( \sum_{n=1}^{k} \frac{1}{n}\right) \|x^0-x_\varepsilon\|^2.
    \end{align*}
\end{lemma}
\begin{proof}
    \begin{align*}
        V^k
        &=
        (k+1)\left\{ k\|x^k-\opT x^k\|^2 + 2\langle x^k-\opT x^k, x^k-x^0 \rangle \right\}
         - \left( \sum_{n=1}^{k} \frac{1}{n}\right) \|x^0-x_\varepsilon\|^2 \\
        &\quad
        + k(k+1) \left\langle - \frac{2}{k}(x^k-x^0) - \left(x_\varepsilon-\opT x_\varepsilon\right), x_\varepsilon-\opT x_\varepsilon \right\rangle
        + \frac{2(k+1)}{k} \left\|x^k-x_\varepsilon + \frac{k}{2} \left(x_\varepsilon-\opT x_\varepsilon\right) \right\|^2 \\
        &\ge
        (k+1)\left\{ k\|x^k-\opT x^k\|^2 + 2\langle x^k-\opT x^k, x^k-x^0 \rangle \right\}
        - \left( \sum_{n=1}^{k} \frac{1}{n}\right) \|x^0-x_\varepsilon\|^2 \\
        &\quad
        + k(k+1) \left\langle - \frac{2}{k}(x^k-x^0) - \left(x_\varepsilon-\opT x_\varepsilon\right), x_\varepsilon-\opT x_\varepsilon \right\rangle \\
        &=
        k(k+1) \left( \|x^k-\opT x^k\|^2 - \|x_\varepsilon-\opT x_\varepsilon\|^2 \right)
        + 2(k+1) \langle (x^k-\opT x^k) - (x_\varepsilon-\opT x_\varepsilon),\, x^k-x^0 \rangle \\
        &\quad
        - \left( \sum_{n=1}^{k} \frac{1}{n}\right) \|x^0-x_\varepsilon\|^2 \\
        &=
        k(k+1) \|(x^k-\opT x^k) - (x_\varepsilon-\opT x_\varepsilon)\|^2
        + 2k(k+1) \langle (x^k-\opT x^k) - (x_\varepsilon-\opT x_\varepsilon),\, x_\varepsilon-\opT x_\varepsilon \rangle \\
        &\quad
        + 2(k+1) \langle (x^k-\opT x^k) - (x_\varepsilon-\opT x_\varepsilon),\, x^k-x_\varepsilon \rangle
        - 2(k+1) \langle (x^k-\opT x^k) - (x_\varepsilon-\opT x_\varepsilon),\, x^0-x_\varepsilon \rangle \\
        &\quad
        - \left( \sum_{n=1}^{k} \frac{1}{n}\right) \|x^0-x_\varepsilon\|^2.
    \end{align*}
    $\opT$ is nonexpansive, from
    \[
        \|\opT x^k - \opT x_\varepsilon\|^2 \le \|x^k - x_\varepsilon\|^2,
    \]
    we get
    \begin{align*}
        \|x^k - x_\varepsilon\|^2 - \|\opT x^k - \opT x_\varepsilon\|^2
        &=
        \left\langle
            (x^k - x_\varepsilon) - (\opT x^k - \opT x_\varepsilon),\,
            (x^k - x_\varepsilon) + (\opT x^k - \opT x_\varepsilon)
        \right\rangle \\
        &=
        \left\langle
            (x^k - \opT x^k) - (x_\varepsilon - \opT x_\varepsilon),\,
            2(x^k - x_\varepsilon) - \{ (x^k - \opT x^k) - (x_\varepsilon - \opT x_\varepsilon) \}
        \right\rangle \\
        &=
        2 \left\langle (x^k - \opT x^k) - (x_\varepsilon - \opT x_\varepsilon),\, x^k - x_\varepsilon \right\rangle
        - \left\| (x^k - \opT x^k) - (x_\varepsilon - \opT x_\varepsilon) \right\|^2 \\
        &\ge 0
    \end{align*}
    so
    \[
        \langle (x^k-\opT x^k) - (x_\varepsilon-\opT x_\varepsilon),\, x^k-x_\varepsilon \rangle
        \ge \frac{1}{2} \| (x^k-\opT x^k) - (x_\varepsilon-\opT x_\varepsilon) \|^2.
    \]
    From this, we get
    \begin{align*}
        V^k
        &\ge
        (k+1)^2 \|(x^k-\opT x^k) - (x_\varepsilon-\opT x_\varepsilon)\|^2
        + 2k(k+1) \langle (x^k-\opT x^k) - (x_\varepsilon-\opT x_\varepsilon),\, x_\varepsilon-\opT x_\varepsilon \rangle \\
        &\quad
        - 2(k+1) \langle (x^k-\opT x^k) - (x_\varepsilon-\opT x_\varepsilon),\, x^0-x_\varepsilon \rangle
        - \left( \sum_{n=1}^{k} \frac{1}{n}\right) \|x^0-x_\varepsilon\|^2.
    \end{align*}
\end{proof}

\begin{lemma} \label{lem:lyap-first}
    Let $\{x^k\}_{k\in\mathbb{N}}$ be a sequence of iterates generated by \eqref{eq:halpern} starting from $x^0\in\cH$ with $\lambda_k = \frac{1}{k+1}$ and $\{V^k\}_{k\in\mathbb{N}\cup\{0\}}$ be defined as \eqref{eq:lyapunov}.
    Then
    \[
        V^1 \le 3\|x^0-x_\varepsilon\|^2.
    \]
\end{lemma}
\begin{proof}
    \begin{align*}
        &V^1 \\
        &=
        2\left\{
            \|x^1-\opT x^1\|^2 + 2\langle x^1-\opT x^1,\, x^1-x^0 \rangle
        \right\}
        - \|x^0-x_\varepsilon\|^2 \\
        &\quad
        + 2 \left\langle -2(x^1-x^0) - (x_\varepsilon-\opT x_\varepsilon),\, x_\varepsilon-\opT x_\varepsilon \right\rangle
        + 4 \left\| x^1-x_\varepsilon + \frac{1}{2}(x_\varepsilon-\opT x_\varepsilon) \right\|^2 \\
        &\le
        0 - 2 \left\langle 2(x^1-x^0) + (x_\varepsilon-\opT x_\varepsilon),\, x_\varepsilon-\opT x_\varepsilon \right\rangle
        + 4 \left\| (x^1-x^0) + (x^0-x_\varepsilon) + \frac{1}{2}(x_\varepsilon-\opT x_\varepsilon) \right\|^2 - \|x^0-x_\varepsilon\|^2 \\
        &=
        \left\| \{2(x^1-x^0) + (x_\varepsilon-\opT x_\varepsilon)\} + 2(x^0-x_\varepsilon) \right\|^2 
        - 2 \langle 2(x^1-x^0) + (x_\varepsilon-\opT x_\varepsilon),\, x_\varepsilon-\opT x_\varepsilon \rangle - \|x^0-x_\varepsilon\|^2 \\
        &=
        \|- \{ (x^0-\opT x^0) - (x_\varepsilon-\opT x_\varepsilon) \} + 2(x^0-x_\varepsilon) \|^2
        + 2 \langle (x^0-\opT x^0) - (x_\varepsilon-\opT x_\varepsilon),\, x_\varepsilon-\opT x_\varepsilon \rangle
        - \|x^0-x_\varepsilon\|^2 \\
        &=
        \|(x^0-\opT x^0) - (x_\varepsilon-\opT x_\varepsilon)\|^2
        - 4 \langle (x^0-\opT x^0) - (x_\varepsilon-\opT x_\varepsilon),\, x^0-x_\varepsilon \rangle \\
        &\quad
        + 3\|x^0-x_\varepsilon\|^2
        + 2 \langle (x^0-\opT x^0) - (x_\varepsilon-\opT x_\varepsilon),\, x_\varepsilon-\opT x_\varepsilon \rangle \\
        &\le
        2 \langle (x^0-\opT x^0) - (x_\varepsilon-\opT x_\varepsilon),\, x_\varepsilon-\opT x_\varepsilon \rangle
        + 3\|x^0-x_\varepsilon\|^2
        - \|(x^0-\opT x^0) - (x_\varepsilon-\opT x_\varepsilon)\|^2 \\
        &=
        3 \|x^0-x_\varepsilon\|^2 - \|x^0-\opT x^0\|^2 - 3\|x_\varepsilon-\opT x_\varepsilon\|^2 \\
        &\le
        3 \|x^0-x_\varepsilon\|^2.
    \end{align*}    
    First inequality comes from \cref{lem:Uk-mono} with $k=0$.
\end{proof}

\begin{theorem} \label{thm:lyap-connect}
    Let $\{x^k\}_{k\in\mathbb{N}}$ be a sequence of iterates generated by \eqref{eq:halpern} starting from $x^0\in\cH$ with $\lambda_k = \frac{1}{k+1}$ and $\{V^k\}_{k\in\mathbb{N}\cup\{0\}}$ be defined as \eqref{eq:lyapunov}.
    For any $k\ge1$,
    \begin{align*}
        &(k+1)^2 \|(x^k-\opT x^k) - (x_\varepsilon-\opT x_\varepsilon)\|^2
        + 2k(k+1) \langle (x^k-\opT x^k) - (x_\varepsilon-\opT x_\varepsilon),\, x_\varepsilon-\opT x_\varepsilon \rangle \\
        &
        - 2(k+1) \langle (x^k-\opT x^k) - (x_\varepsilon-\opT x_\varepsilon),\, x^0-x_\varepsilon \rangle
        - \left( \sum_{n=1}^{k} \frac{1}{n}\right) \|x^0-x_\varepsilon\|^2 \\
        &\le
        3 \|x^0-x_\varepsilon\|^2.
    \end{align*}
\end{theorem}
\begin{proof}
    Direct application of Lemma \ref{lem:lyap-mono}, Lemma \ref{lem:lyap-last} and Lemma \ref{lem:lyap-first}.
\end{proof}

\begin{lemma} \label{lem:fpr-mono}
    Let $\{x^k\}_{k\in\mathbb{N}}$ be a sequence of iterates generated by \eqref{eq:halpern} starting from $x^0\in\cH$ with $\lambda_k = \frac{1}{k+1}$.
    For any $k\in\mathbb{N}$,
    \[
        \|x^k-\opT x^k\| \le \|x^0-\opT x^0\|.
    \]
\end{lemma}
\begin{proof}
    We use \cref{lem:Uk-mono}, the definition of $x^{k+1}$-update and and that $\theta_k = \frac{k}{2}$, which is from \cref{lem:haltheta}.
    Dividing by $\frac{(k+1)(k+2)}{2}$, we have
    \begin{align*}
        &\frac{2k}{k+2} \|x^k-\opT x^k\|^2 + 4 \langle x^k - \opT x^k,\, x^k-x^0 \rangle \\
        &\ge
        2 \|x^{k+1}-\opT x^{k+1}\|^2 + \frac{4}{k+1} \langle x^{k+1} - \opT x^{k+1},\, x^{k+1} - x^0 \rangle \\
        &=
        \|x^{k+1}-\opT x^{k+1}\|^2 + \left\| (x^{k+1}-\opT x^{k+1}) + \frac{x^{k+1}-x^0}{\theta_{k+1}} \right\|^2 - \left\| \frac{x^{k+1}-x^0}{\theta_{k+1}} \right\|^2.
    \end{align*}
    Since
    \begin{align*}
        \frac{x^{k+1}-x^0}{\theta_{k+1}}
        &=
        \frac{k}{k+2} \left( \frac{x^k-x^0}{\theta_k} \right) - \frac{2}{k+2} (x^k-\opT x^k),
    \end{align*}
    we have
    \begin{align*}
        &\|x^{k+1}-\opT x^{k+1}\|^2
        + \left\| (x^{k+1}-\opT x^{k+1}) + \frac{x^{k+1}-x^0}{\theta_{k+1}} \right\|^2 \\
        &\le
        \frac{2k}{k+2} \|x^k-\opT x^k\|^2 
        + 4 \langle x^k-\opT x^k,\, x^k-x^0 \rangle
        + \left\| \frac{k}{k+2} \left( \frac{x^k-x^0}{\theta_k} \right) - \frac{2}{k+2} (x^k-\opT x^k) \right\|^2 \\
        &=
        \|x^k-\opT x^k\|^2 + \left( \frac{k}{k+2} \right)^2 \left\| (x^k-\opT x^k) + \frac{x^k-x^0}{\theta_k} \right\|^2
    \end{align*}
    hold for all $k=0,1,\ldots$.
    Therefore, for any $k\in\mathbb{N}$, we get
    \begin{align*}
        \|x^0-\opT x^0\|^2
        &\ge
        \|x^1-\opT x^1\|^2
        + \left\| (x^1-\opT x^1) + \frac{x^1-x^0}{\theta_1} \right\|^2 \\
        &\ge
        \|x^1-\opT x^1\|^2 
        + \left( \frac{1}{1+2} \right)^2 \left\| (x^1-\opT x^1) + \frac{x^1-x^0}{\theta_1} \right\|^2 \\
        &\ge
        \|x^2-\opT x^2\|^2
        + \left\| (x^2-\opT x^2) + \frac{x^2-x^0}{\theta_2} \right\|^2 \\
        &\ge
        \cdots \\
        &\ge
        \|x^k-\opT x^k\|^2 + \left\| (x^k-\opT x^k) + \frac{x^k-x^0}{\theta_k} \right\|^2 \\
        &\ge
        \|x^k-\opT x^k\|^2.
    \end{align*}
\end{proof}

Now we find some relation between $x_\varepsilon-\opT x_\varepsilon$ and $v$.

\begin{lemma} \label{lem:vconversion}
    Let $\{x^k\}_{k\in\mathbb{N}}$ be a sequence of iterates generated by \eqref{eq:halpern} starting from $x^0\in\cH$ with $\lambda_k = \frac{1}{k+1}$.
    For any $\varepsilon>0$, there exists $x_\varepsilon\in\dom\opT$ such that
    \[
        \|x_\varepsilon-\opT x_\varepsilon-v\|\le\varepsilon,
    \]
    and from this,
    \[
        \|x^k-\opT x^k-(x_\varepsilon-\opT x_\varepsilon)\|^2 
        \ge
        \|x^k-\opT x^k-v\|^2 - 2\|x^0-\opT x^0\| \varepsilon
    \]
    and
    \[
        \langle x^k-\opT x^k - (x_\varepsilon-\opT x_\varepsilon), x_\varepsilon - \opT x_\varepsilon \rangle
        \ge
        \langle x^k-\opT x^k-v, v \rangle 
        - \left\{
            \|x^0-\opT x^0\| + 2\|v\| + \varepsilon
        \right\} \varepsilon.
    \]
    Furthermore, if $v\in\cR(\opI-\opT)$, then there exists $x_\star\in\dom\opT$ such that $x_\star-\opT x_\star = v$.
\end{lemma}
\begin{proof}
    \begin{align*}
        &\|x^k-\opT x^k-(x_\varepsilon-\opT x_\varepsilon)\|^2 - \|x^k-\opT x^k-v\|^2 \\
        &=
        - 2 \langle x^k-\opT x^k, x_\varepsilon-\opT x_\varepsilon-v \rangle
        + \underbrace{\|x_\varepsilon-\opT x_\varepsilon\|^2 - \|v\|^2}_{\ge0} \\
        &\ge
        - 2\|x^k-\opT x^k\|\|x_\varepsilon-\opT x_\varepsilon-v\| \\
        &\ge
        - 2\|x^0-\opT x^0\|\varepsilon
    \end{align*}
    where the last inequality comes from \cref{lem:fpr-mono}.
    Also,
    \begin{align*}
        &\langle x^k-\opT x^k-(x_\varepsilon-\opT x_\varepsilon), x_\varepsilon-\opT x_\varepsilon \rangle \\
        &=
        \langle (x^k-\opT x^k-v) - (x_\varepsilon-\opT x_\varepsilon-v), (x_\varepsilon-\opT x_\varepsilon-v) + v \rangle \\
        &=
        \langle x^k-\opT x^k-v, v \rangle + \langle x^k-\opT x^k-v, x_\varepsilon-\opT x_\varepsilon-v \rangle - \|x_\varepsilon-\opT x_\varepsilon-v\|^2 \\
        &\quad
        - \langle x_\varepsilon-\opT x_\varepsilon-v, v\rangle \\
        &\ge
        \langle x^k-\opT x^k-v, v \rangle - \|x^k-\opT x^k-v\|\|x_\varepsilon-\opT x_\varepsilon-v\| \\
        &\quad
        -\|x_\varepsilon-\opT x_\varepsilon-v\|^2 - \|x_\varepsilon-\opT x_\varepsilon-v\|\|v\|.
    \end{align*}
    Then
    \begin{align*}
        &\langle x^k-\opT x^k-(x_\varepsilon-\opT x_\varepsilon), x_\varepsilon-\opT x_\varepsilon \rangle
        - \langle x^k-\opT x^k-v, v \rangle \\
        &\ge
        - \|x_\varepsilon-\opT x_\varepsilon-v\| \left\{
            \|x^k-\opT x^k-v\| + \|x_\varepsilon-\opT x_\varepsilon-v\| + \|v\|
        \right\} \\
        &\ge
        - \varepsilon \left\{
            (\|x^k-\opT x^k\| + \|v\|) + \varepsilon + \|v\|
        \right\} \\
        &=
        - \varepsilon \left\{
            \|x^0-\opT x^0\| + 2\|v\| + \varepsilon
        \right\}
    \end{align*}
    where the last inequality comes from \cref{lem:fpr-mono}.
\end{proof}

We now prove the convergence rate result of \eqref{eq:halpern} with $\lambda_k=\frac{1}{k+1}$.

\begin{theorem}
    \label{thm:halfprconvalpha}
    Let $\{x^k\}_{k\in\mathbb{N}}$ be a sequence of iterates generated by \eqref{eq:halpern} starting from $x^0\in\cH$ with $\lambda_k = \frac{1}{k+1}$.
    For any $\varepsilon>0$ and $0<\alpha<1$, there exists $x_\varepsilon\in\dom\opT$ such that
    \[
        \|x^{k}-\opT x^{k}-v\|^2
        \le
        \frac{1}{(1-\alpha)(k+1)^2} \left( \sum_{n=1}^{k} \frac{1}{n} + 3 + \frac{1}{\alpha} \right)\|x^0-x_\varepsilon\|^2
        + \varepsilon.
    \]
    If we further assume that $v\in\cR(\opI-\opT)$, there exists $x_\star\in\cH$ such that $v = x_\star-\opT x_\star$ and
    \[
        \|x^{k}-\opT x^{k}-v\|^2
        \le
        \frac{1}{(1-\alpha)(k+1)^2} \left( \sum_{n=1}^{k} \frac{1}{n} + 3 + \frac{1}{\alpha} \right)\|x^0-x_\star\|^2.
    \]
\end{theorem}
\begin{proof}
    For $\varepsilon>0$ and $0<\alpha<1$, consider $x_\varepsilon\in\dom\opT$ such that
    \[
        \|x_\varepsilon-\opT x_\varepsilon-v\|
        \le
        \tilde{\varepsilon}
    \]
    where
    \[
        \tilde{\varepsilon} = \min\left\{ \left( 2\|x^0-\opT x^0\| + \frac{2}{1-\alpha} (\|x^0-\opT x^0\|+2\|v\|+1) \right)^{-1}{\varepsilon}, 1, \varepsilon \right\}.
    \]
    
    According to Theorem \ref{thm:lyap-connect},
    \begin{align*}
        &(k+1)^2 \|(x^k-\opT x^k) - (x_\varepsilon-\opT x_\varepsilon)\|^2
        + 2k(k+1) \langle (x^k-\opT x^k) - (x_\varepsilon-\opT x_\varepsilon),\, x_\varepsilon-\opT x_\varepsilon \rangle \\
        &
        - 2(k+1) \langle (x^k-\opT x^k) - (x_\varepsilon-\opT x_\varepsilon),\, x^0-x_\varepsilon \rangle
        - \left( \sum_{n=1}^{k} \frac{1}{n}\right) \|x^0-x_\varepsilon\|^2 \\
        &\le
        3 \|x^0-x_\varepsilon\|^2.
    \end{align*}
    For any $\alpha\in(0,1)$, 
    \begin{align*}
        3 \|x^0-x_\varepsilon\|^2
        &\ge
        (1-\alpha)(k+1)^2 \|(x^k-\opT x^k) - (x_\varepsilon-\opT x_\varepsilon)\|^2
        + 2k(k+1) \langle (x^k-\opT x^k) - (x_\varepsilon-\opT x_\varepsilon),\, x_\varepsilon-\opT x_\varepsilon \rangle \\
        &\quad
        + \alpha (k+1)^2 \|(x^k-\opT x^k) - (x_\varepsilon-\opT x_\varepsilon)\|^2 + 2(k+1) \langle (x^k-\opT x^k)-(x_\varepsilon-\opT x_\varepsilon),\, x^0-x_\varepsilon \rangle \\
        &\quad
        + \frac{1}{\alpha} \|x^0-x_\varepsilon\|^2
        - \left( \frac{1}{\alpha} + \sum_{n=1}^k \frac{1}{n} \right) \|x^0-x_\varepsilon\|^2 \\
        &=
        (1-\alpha)(k+1)^2 \|(x^k-\opT x^k) - (x_\varepsilon-\opT x_\varepsilon)\|^2
        + 2k(k+1) \langle (x^k-\opT x^k) - (x_\varepsilon-\opT x_\varepsilon),\, x_\varepsilon-\opT x_\varepsilon \rangle \\
        &\quad
        + \frac{1}{\alpha} \left\| \alpha(k+1) \{ (x^k-\opT x^k) - (x_\varepsilon-\opT x_\varepsilon) \} + (x^0-x_\varepsilon) \right\|^2
        - \left( \frac{1}{\alpha} + \sum_{n=1}^k \frac{1}{n} \right) \|x^0-x_\varepsilon\|^2 \\
        &\ge
        (1-\alpha)(k+1)^2 \|(x^k-\opT x^k) - (x_\varepsilon-\opT x_\varepsilon)\|^2
        - \left( \frac{1}{\alpha} + \sum_{n=1}^k \frac{1}{n} \right) \|x^0-x_\varepsilon\|^2 \\
        &\quad
        + 2k(k+1) \langle (x^k-\opT x^k) - (x_\varepsilon-\opT x_\varepsilon),\, x_\varepsilon-\opT x_\varepsilon \rangle.
    \end{align*}
    Rearranging the terms, we get
    \begin{align*}
        &\frac{1}{(1-\alpha)(k+1)^2} \left( 3 + \frac{1}{\alpha} + \sum_{n=1}^k \frac{1}{n} \right) \|x^0-x_\varepsilon\|^2 \\
        &\ge
        \|(x^k-\opT x^k) - (x_\varepsilon-\opT x_\varepsilon)\|^2
        + \frac{2}{1-\alpha}\frac{k}{k+1} \langle (x^k-\opT x^k) - (x_\varepsilon-\opT x_\varepsilon),\, x_\varepsilon-\opT x_\varepsilon \rangle \\
        &\ge
        \|x^k-\opT x^k-v\|^2 - 2\|x^0-\opT x^0\|\tilde{\varepsilon}
        + \frac{2}{1-\alpha}\frac{k}{k+1} \langle x^k-\opT x^k-v,\, v \rangle \\
        &\quad
        - \frac{2}{1-\alpha}\frac{k}{k+1} \left\{ \|x^0-\opT x^0\| + 2\|v\| + \tilde{\varepsilon} \right\} \tilde{\varepsilon} \\
        &\ge
        \|x^k-\opT x^k-v\|^2
        - \frac{2}{1-\alpha} \left\{
            (2-\alpha) \|x^0-\opT x^0\| + 2\|v\| + \tilde{\varepsilon}
        \right\} \tilde{\varepsilon} \\
        &\ge
        \|x^k-\opT x^k-v\|^2 - \varepsilon.
    \end{align*}
    The second inequality comes from \cref{lem:vconversion}, the third inequality comes from \cref{lem:normconvequiv}, and the last inequality comes from the definition of $\tilde{\varepsilon}>0$.
\end{proof}

\begin{proof}[Proof of \cref{thm:halfprconv}]
    With $\lambda_k = \frac{1}{k+1}$, we have $\theta_k = \frac{k}{2}$ from \cref{lem:haltheta}.
    Using \cref{lem:Uk-mono}, we get
    \begin{align*}
        &(k+2) \left\{ (k+1)\|x^{k+1} - \opT x^{k+1}\|^2 + 2 \langle x^{k+1}-\opT x^{k+1}, x^{k+1}-x^0 \rangle \right\} \\
        &\le
        (k+1) \left\{ k \|x^k-\opT x^k\|^2 + 2 \langle x^k-\opT x^k, x^k-x^0 \rangle \right\}
    \end{align*}
    for all $k=0,1,\dots$.
    Therefore, for any $k\in\mathbb{N}$,
    \[
        k \|x^k-\opT x^k\|^2 + 2 \langle x^k-\opT x^k, x^k-x^0 \rangle \le 0.
    \]
    Using the Cauchy-Schwarz inequality, we get
    \[
        \|x^k-\opT x^k\| \le \left\| \frac{2}{k} (x^k-x^0) \right\| = \left\| \frac{x^k-x^0}{\theta_k} \right\|
    \]
    for any $k\in\mathbb{N}$.
    Therefore, for any $\varepsilon>0$ and $x_\varepsilon\in\cH$ such that $\|x_\varepsilon - \opT x_\varepsilon\|^2 - \|v\|^2 \le \varepsilon^2$, we have
    \[
        \|x^k-\opT x^k\| - \|v\|
        \le
        \left\| \frac{2}{k} (x^k-x^0) \right\| - \|v\|
        \le
        \left\| \frac{2}{k} (x^k-x^0) + v \right\|
        \le
        \frac{4}{k} \|x^0-x_\varepsilon\| + \varepsilon
    \]
    for any $k\in\mathbb{N}$, where the second from last inequality comes from \cref{thm:halaverconv}.

    From \cref{thm:halfprconvalpha}, given an arbitrary $\varepsilon>0$ and $x_\varepsilon$ such that
    \[
        \|x_\varepsilon-\opT x_\varepsilon-v\|
        \le
        \tilde{\varepsilon}
    \]
    where
    \[
        \tilde{\varepsilon} = \min\left\{ \left( 2\|x^0-\opT x^0\| + \frac{2}{1-\alpha} (\|x^0-\opT x^0\|+2\|v\|+1) \right)^{-1}{\varepsilon}, 1, \varepsilon \right\}
        = \cO(\varepsilon),
    \]
    we get
    \[
        \|x^k-\opT x^k-v\|^2
        \le
        \frac{1}{(1-\alpha)(k+1)^2} \left( 3 + \frac{1}{\alpha} + \sum_{n=1}^k \frac{1}{n} \right) \|x^0-x_\varepsilon\|^2 + \varepsilon
    \]
    for any $0 < \alpha < 1$.
    Now we find a minimizer $\alpha^\star$ of 
    \[
        \frac{1}{1-\alpha} \left( c_k + \frac{1}{\alpha} \right)
    \]
    where
    \[
        c_k = 3 + \sum_{n=1}^k \frac{1}{n}
    \]
    is a positive constant.
    \begin{align*}
        \frac{1}{1-\alpha} \left(c_k + \frac{1}{\alpha} \right)
        &=
        \frac{c_k}{1-\alpha} + \frac{1}{\alpha(1-\alpha)} \\
        &=
        \frac{c_k+1}{1-\alpha} + \frac{1}{\alpha} \\
        &=
        \left( \frac{c_k+1}{1-\alpha} + \frac{1}{\alpha} \right) \left( (1-\alpha) + \alpha \right) \\
        &\ge
        (\sqrt{c_k+1} + 1)^2
        \tag{$\because$ Cauchy-Schwarz.}
    \end{align*}
    and the equality holds if and only if
    \[
        \frac{c_k+1}{(1-\alpha)^2} = \frac{1}{\alpha^2}
        \quad\Leftrightarrow\quad
        \alpha = \frac{1}{\sqrt{c_k+1}+1} = \frac{\sqrt{c_k+1}-1}{c_k}.
    \]
    With such $\alpha$, we get
    \[
        \frac{1}{1-\alpha} \left(c_k+\frac{1}{\alpha}\right)
        = (\sqrt{c_k+1}+1)^2.
    \]
    Therefore,
    \[
        \|x^k-\opT x^k-v\|^2 \le \left( \frac{\sqrt{\sum_{n=1}^k \frac{1}{n} + 4} + 1}{k+1} \right)^2 \|x^0-x_\varepsilon\|^2 + \varepsilon.
    \]
    
    If $v=x_\star-\opT x_\star$, we follow the same steps and get
    \[
        \|x^k-\opT x^k-v\|^2 \le \left( \frac{\sqrt{\sum_{n=1}^k \frac{1}{n} + 4} + 1}{k+1} \right)^2 \|x^0-x_\star\|^2.
    \]

\end{proof}

We now prove the equivalence of the normalized iterate $-\frac{x^{k+1}-x^0}{k+1}$ of Picard iteration and the fixed-point residual $x^k-\opT x^k$ of \eqref{eq:halpern} with $\lambda_k = \frac{1}{k+1}$ for affine $\opT$, which was discussed in the last part of \cref{subsec:halpern}.
Let $\opT\colon\cH\to\cH$ be an affine operator, i.e., $\opT x = A x + b$ where $A\colon\cH\to\cH$ is a linear operator and $b\in\cH$.

\begin{lemma} \label{lem:picardequiv}
    Suppose $\opT\colon\cH\to\cH$ is an affine operator.
    Let $\{x^k\}_{k\in\mathbb{N}}$ and $\{y^k\}_{k\in\mathbb{N}}$ be the sequences of iterates generated by \eqref{eq:halpern} with $\lambda_k = \frac{1}{k+1}$ and Picard iteration with $\opT$, respectively, starting from the same initial point $x^0 = y^0$.
    Then for any $k\in\mathbb{N}\cup\{0\}$,
    \[
        x^k - \opT x^k = - \frac{y^{k+1} - y^0}{k+1}.
    \]
\end{lemma}

\begin{proof}
    First, note that when $\opT$ is an affine operator, i.e., $\opT x = Ax + b$ for any $x\in\cH$,
    \[
        \opT \left( \sum_{i=1}^k \nu_i \opT x_i \right)
        = 
        \sum_{i=1}^k \nu_i \opT x_i
    \]
    for any $x_i \in \cH$ and $\nu_i \in [0,1]$ such that $\sum_{i=1}^k \nu_i = 1$.
    
    We see that for Picard iteration,
    \[
        - \frac{y^{k+1} - y^0}{k+1} = - \frac{\opT^{k+1} y^0 - y^0}{k+1} = - \frac{\opT^{k+1} x^0 - x^0}{k+1}.
    \]
    Considering the \eqref{eq:halpern} iterates $\{x^k\}_{k\in\mathbb{N}}$ starting from $x^0$, we claim by induction on $k$ that
    \[
        x^k = \frac{1}{k+1} \sum_{i=0}^k \opT^i x^0.
    \]
    When $k=1$,
    \[
        x^1 = \frac{1}{2} \opT x^0 + \frac{1}{2} x^0.
    \]
    Now, suppose the claim holds for $k=n$.
    \begin{align*}
        x^{n+1}
        &=
        \frac{n+1}{n+2} \opT x^n + \frac{1}{n+2} x^0 \\
        &=
        \frac{n+1}{n+2} \opT \left( \frac{1}{n+1} \sum_{i=0}^n \opT^i x^0 \right) + \frac{1}{n+2} x^0 \\
        &=
        \frac{n+1}{n+2} \left( \frac{1}{n+1} \sum_{i=0}^n \opT^{i+1} x^0 \right) + \frac{1}{n+2} x^0 \\
        &=
        \frac{1}{n+2} \sum_{i=0}^{n+1} \opT^i x^0.
    \end{align*}
    Therefore,
    \begin{align*}
        x^k - \opT x^k
        &=
        \frac{1}{k+1} \sum_{i=0}^k \opT^i x^0 - \opT \left( \frac{1}{k+1} \sum_{i=0}^k \opT^i x^0 \right) \\
        &=
        \frac{1}{k+1} \sum_{i=0}^k \opT^i x^0 - \frac{1}{k+1} \sum_{i=0}^k \opT^{i+1} x^0 \\
        &=
        \frac{1}{k+1} x^0 - \frac{1}{k+1} \opT^{k+1} x^0 \\
        &=
        - \frac{\opT^{k+1} x^0 - x^0}{k+1}.
    \end{align*}
\end{proof}

Due to \cref{lem:picardequiv}, when $\opT$ is an affine nonexpansive operator, \eqref{eq:halpern} with $\lambda_k = \frac{1}{k+1}$ is exactly optimal with matching lower bound, for fixed-point residual.

\subsection{Omitted proofs of \cref{subsec:mann}}
\label{subsec:proofmann}
Consider a \emph{Mann iteration}
\begin{equation*}
    x^k = \sum_{i=0}^{k} \nu^k_i \opT x^{i-1}
    \tag{Mann} \label{eq:manniterationappendix}
\end{equation*}
where $\nu_i \ge 0$, $\sum_{i=0}^k \nu^k_i = 1$ and $\opT x^{-1} := x^0$.

\begin{lemma} \label{lem:mannnormalizingfactor}
    Let $\alpha_0 = 0$ and $\{x^k\}_{k\in\mathbb{N}}$ be a sequence of iterates generated by \eqref{eq:manniterationappendix} starting from $x^0\in\cH$.
    If the sequence of real numbers $\{\alpha_k\}_{k\in\mathbb{N}\cup\{0\}}$ is defined recursively from the equation
    \[
        \alpha_{k} = (1-\nu^k_0) + \sum_{i=1}^{k} \nu^k_i \alpha_{i-1}, \quad k=1,\ldots,
    \]
    and $\alpha_k > 0$ for all $k\in\mathbb{N}$, then
    \[
        - \frac{x^k-x^0}{\alpha_k} \in \overline{\cR(\opI-\opT)}.
    \]
\end{lemma}

\begin{proof}
    Note that $\alpha_k$ for $k\ge1$ can also be written as
    \[
        \alpha_k
        =
        \sum_{i=1}^{k} \nu_i^k + \sum_{i=1}^k \nu^k_i \alpha_{i-1} 
    \]
    since $\sum_{i=0}^k \nu_i^k = 1$.

    Let $k=0$.
    Then from the definition of \eqref{eq:manniterationappendix}, $\nu_0^0=1$.
    If $k=1$, $\alpha_1 = \nu^1_1$ and
    \[
        x^1 = \nu^1_0 \opT x^{-1} + \nu^1_1 \opT x^0 = x^0 - \nu^1_1 (x^0-\opT x^0).
    \]
    so
    \[
        - \frac{x^1-x^0}{\alpha_1} = - \frac{x^1-x^0}{\nu^1_1} = x^0-\opT x^0 \in \overline{\cR(\opI-\opT)}.
    \]
    Now, fix $k>1$ and suppose that
    \[
        - \frac{x^i-x^0}{\alpha_i} \in \overline{\cR(\opI-\opT)}, \quad \forall \, i < k.
    \]
    Then from
    \begin{align*}
        x^{k} - x^0
        &=
        \sum_{i=0}^{k} \nu_i^{k} (\opT x^{i-1} - x^0) \\
        &=
        \sum_{i=1}^k \nu_i^k (\opT x^{i-1} - x^0) \\
        &=
        - \sum_{i=1}^{k} \nu_i^{k} (x^{i} - \opT x^{i})
        + \sum_{i=1}^{k} \nu_i^{k} (x^{i-1} - x^0),
    \end{align*}
    we get
    \begin{align*}
        - \frac{x^{k}-x^0}{\alpha_{k}}
        &=
        \frac
        {\sum_{i=1}^{k} \nu_i^{k} (x^{i-1}-\opT x^{i-1}) + \sum_{i=1}^{k} \nu_i^{k} \alpha_{i-1} \left( - \frac{x^{i-1}-x^0}{\alpha_{i-1}} \right)}
        {\sum_{i=1}^{k} \nu_i^{k} + \sum_{i=1}^{k} \nu_i^{k} \alpha_{i-1}}.
    \end{align*}
    Since $\overline{\cR(\opI-\opT)}$ is a closed convex set, it is closed under convex combination.
    Therefore, $-\frac{x^k-x^0}{\alpha_k} \in \overline{\cR(\opI-\opT)}$.
\end{proof}

\begin{remark}
    Note that $\{\alpha_k\}_{k\in\mathbb{N}\cup\{0\}}$ of \eqref{lem:mannnormalizingfactor} recovers $\sum_{i=1}^{k} (1-\lambda_i)$ of \eqref{eq:kmiteration} and $\theta_k$ of \eqref{eq:halpern}.

    First of all, \eqref{eq:kmiteration} is defined as
    \[
        x^{k} = (1-\lambda_{k}) \opT x^{k-1} + \lambda_{k} x^{k-1},
    \]
    so $\nu_k^k = 1-\lambda_k$.
    From recursively applying the same identity as above, we get 
    \[
        \nu_i^k = \begin{cases}
            1-\lambda_k & \text{if } i=k \\
            \lambda_k \cdots \lambda_{i+1} (1-\lambda_{i}) & \text{if } 1 \le i < k \\
            \lambda_{k+1} \lambda_k \cdots \lambda_1 & \text{if } i=0
        \end{cases}
    \]
    From \cref{lem:mannnormalizingfactor}, as
    \[
        \alpha_k = \sum_{i=1}^{k} \nu_i^k + \sum_{i=1}^k \nu_i^k \alpha_{i-1}
    \]
    with $\alpha_0 = 0$, we get
    \[
        \alpha_k = \sum_{i=1}^k (1-\lambda_i)
    \]
    from plugging $\nu_i^k$ above.
    
    Next, \eqref{eq:halpern} is defined as
    \[
        x^k = (1-\lambda_k) \opT x^k + \lambda_k x^0
    \]
    so 
    \[
        \nu_i^k = \begin{cases}
            1 - \lambda_k & \text{if } i=k \\
            0 & \text{if } 1\le i < k \\
            \lambda_k & \text{if } i=0
        \end{cases}
    \]
    Then
    \begin{align*}
        \alpha_k 
        &=
        \sum_{i=1}^{k} \nu_i^k + \sum_{i=1}^k \nu_i^k \alpha_{i-1}
        =
        (1-\lambda_k) + (1-\lambda_k) \alpha_{k-1}
        = 
        (1-\lambda_k) (1+\alpha_{k-1}).
    \end{align*}
    This recursive formula is exactly the same as the recursive formula in \cref{lem:haltheta}, so $\alpha_k = \theta_k$.
\end{remark}

We elaborate on some properties of \eqref{eq:manniterationappendix} which will be used in our main result, \cref{thm:mannaverconv}.
\begin{lemma} \label{lem:mannupperbound}
    Let $\{x^k\}_{k\in\mathbb{N}}$ and $\{y^k\}_{k\in\mathbb{N}}$ be the sequences of iterates generated by \eqref{eq:manniterationappendix} starting from $x^0\in\cH$ and $y^0\in\cH$, respectively.
    Define $\{\alpha_k\}_{k\in\mathbb{N}\cup\{0\}}$ as in \cref{lem:mannnormalizingfactor}.
    Then
    \[
        \left\| \frac{x^k-x^0}{\alpha_k} \right\| \le \|x^0-\opT x^0\|, \quad k=1,2,\ldots
    \]
    and
    \[
        \|x^k-y^k\| \le \|x^0-y^0\|, \quad k=1,2,\ldots.
    \]
\end{lemma}

\begin{proof}
    For $k=1$, $ - \frac{x^1-x^0}{\alpha_1} = x^0-\opT x^0$ so the claim is trivial.
    Now, let 
    \[
        \left\| \frac{x^i-x^0}{\alpha_i} \right\| \le \|x^0-\opT x^0\|
    \]
    for all $i < k$. Then
    \begin{align*}
        x^{k} - x^0
        &=
        \sum_{i=1}^{k} \nu_i^{k} (\opT x^{i-1} - x^0) \\
        &=
        - \sum_{i=1}^{k} \nu_i^{k} (x^0 - \opT x^0) + \sum_{i=1}^k \nu_i^k (\opT x^{i-1} - \opT x^0) \\
        &=
        - (1-\nu_0^k) (x^0 - \opT x^0) + \sum_{i=1}^k \nu_i^k (\opT x^{i-1} - \opT x^0),
    \end{align*}
    so
    \begin{align*}
        \left\| \frac{x^k-x^0}{\alpha_k} \right\|
        &=
        \frac{1}{\alpha_k}
        \left\|
            (1 - \nu_0^k) (x^0-\opT x^0) 
            - \sum_{i=2}^k \nu_i^k \alpha_{i-1} \left( \frac{\opT x^{i-1} - \opT x^0}{\alpha_{i-1}} \right)
        \right\| \\
        &\le
        \frac{1}{\alpha_k} \left\{
            (1-\nu_0^k) \|x^0-\opT x^0\|
            + \sum_{i=2}^k \nu_i^k \alpha_{i-1} \left\| \frac{\opT x^{i-1} - \opT x^0}{\alpha_{i-1}} \right\|
        \right\} \\
        &\le
        \frac{1}{\alpha_k} \left\{
            (1-\nu_0^k) \|x^0-\opT x^0\|
            + \sum_{i=2}^k \nu_i^k \alpha_{i-1} \left\| \frac{x^{i-1} - x^0}{\alpha_{i-1}} \right\|
        \right\} \\
        &\le
        \frac{1}{\alpha_k} \left\{
            (1-\nu_0^k) \|x^0-\opT x^0\|
            + \sum_{i=2}^k \nu_i^k \alpha_{i-1} \| x^0 - \opT x^0 \|
        \right\} \\
        &=
        \frac{1}{\alpha_k} \left\{ (1-\nu_k^k) + \sum_{i=2}^k \nu_i^k \alpha_{i-1} \right\} \|x^0-\opT x^0\| = \|x^0-\opT x^0\|.
    \end{align*}

    Now we prove the second claim.
    First of all,
    \begin{align*}
        \|x^1-y^1\|
        &= \|\nu_0^1 (x^0-y^0) + \nu_1^1 (\opT x^0 - \opT y^0)\| \\
        &\le \nu_0^1 \| x^0 - y^0 \| + \nu_1^1 \| \opT x^0 - \opT y^0 \| \\
        &\le
        \nu_0^1 \|x^0-y^0\| + \nu_1^1 \|x^0-y^0\|
        = \|x^0-y^0\|.
    \end{align*}
    Suppose $\|x^i-y^i\| \le \|x^0-y^0\|$ for all $i < k$.
    Then
    \begin{align*}
        \|x^k-y^k\|
        &=
        \left\| \sum_{i=0}^k \nu_i^k (\opT x^{i-1} - \opT y^{i-1}) \right\| \\
        &\le
        \sum_{i=0}^k \nu_i^k \|\opT x^{i-1} - \opT y^{i-1}\| \\
        &\le
        \nu_0^0 \|x^0-y^0\| + \sum_{i=1}^k \nu_i^k \|x^{i-1}-y^{i-1}\| \\
        &\le
        \nu_0^0 \|x^0-y_0\| + \sum_{i=1}^k \nu_i^k \|x^0-y^0\|
        = \|x^0-y^0\|.
    \end{align*}
\end{proof}

We can extend \cref{thm:kmaverconv} and \cref{thm:halaverconv} to cover the case of general Mann iteration.
\begin{theorem} \label{thm:mannaverconv}
    Let $\{x^k\}_{k\in\mathbb{N}}$ be a sequence in $\cH$ generated by \eqref{eq:manniterationappendix} starting from $x^0\in\cH$ and $\{\alpha_k\}_{k\in\mathbb{N}\cup\{0\}}$ be a sequence of positive numbers defined in \cref{lem:mannnormalizingfactor}.
    Then for any $\varepsilon>0$, there exists $x_\varepsilon\in\cH$ such that $\|x_\varepsilon - \opT x_\varepsilon - v\| < \varepsilon$ and
    \[
        \left\| \frac{x^k-x^0}{\alpha_k} + v \right\|
        \le
        \frac{2}{\alpha_k} \|x^0-x_\varepsilon\| + \varepsilon.
    \]
    If we further assume that $v \in \cR(\opI-\opT)$, then there exists $x_\star\in\cH$ such that
    \[
        \left\| \frac{x^k-x^0}{\alpha_k} + v \right\|
        \le
        \frac{2}{\alpha_k} \|x^0-x_\star\|.
    \]
    Therefore, if $\lim_{k\to\infty} \alpha_k = \infty$, then
    \[
        \lim_{k\to\infty} \frac{x^k-x^0}{\alpha_k} = -v.
    \]
\end{theorem}

\begin{proof}
    Fix $\varepsilon>0$.
    Let $x_\varepsilon\in\cH$ be a point such that $\|x_\varepsilon - \opT x_\varepsilon\|^2  - \|v\|^2 < \varepsilon^2$.
    Then
    \begin{align*}
        \|x_\varepsilon - \opT x_\varepsilon - v\|^2 
        &= \|x_\varepsilon - \opT x_\varepsilon\|^2 - \|v\|^2 - 2 \underbrace{ \langle x_\varepsilon - \opT x_\varepsilon - v,\, v \rangle }_{\ge 0} \\
        &\le \|x_\varepsilon - \opT x_\varepsilon\|^2  - \|v\|^2 \\
        &< \varepsilon^2.
    \end{align*}
    Suppose $\{x_\varepsilon^k\}_{k\in\mathbb{N}}$ be a sequence in $\cH$ generated by \eqref{eq:manniteration} starting from $x_\varepsilon$.
    Since $- \frac{x^k_\varepsilon-x_\varepsilon}{\alpha_k} \in \overline{\cR(\opI-\opT)}$ by \cref{lem:mannnormalizingfactor}, 
    \[
        \left\langle - \frac{x^k_\varepsilon-x_\varepsilon}{\alpha_k},\, v \right\rangle \ge \|v\|^2
    \]
    for any $k\in\mathbb{N}$, so we get
    \begin{align*}
        \left\| \frac{x^k_\varepsilon - x_\varepsilon}{\alpha_k} + v \right\|^2
        &=
        \left\| \frac{x^k_\varepsilon - x_\varepsilon}{\alpha_k} \right\|^2
        + \|v\|^2
        - 2 \left\langle - \frac{x^k_\varepsilon-x_\varepsilon}{\alpha_k},\, v \right\rangle \\
        &\le
        \left\| \frac{x^k_\varepsilon - x_\varepsilon}{\alpha_k} \right\|^2
        - \|v\|^2 \\
        &\le
        \| x_\varepsilon - \opT x_\varepsilon \|^2 - \|v\|^2 \le \varepsilon^2
    \end{align*}
    or $\left\| \frac{x^k_\varepsilon - x_\varepsilon}{\alpha_k} + v \right\| \le \varepsilon$.
    \begin{align*}
        \left\| \frac{x^k - x^0}{\alpha_k} + v \right\|
        &\le
        \left\| \frac{x^k - x^0}{\alpha_k} - \frac{x_\varepsilon^k - x_\varepsilon}{\alpha_k} \right\|
        + \left\| \frac{x_\varepsilon^k - x_\varepsilon}{\alpha_k} + v \right\| \\
        &\le
        \frac{\|x^k-x^k_\varepsilon\|}{\alpha_k}
        + \frac{\|x^0-x_\varepsilon\|}{\alpha_k}
        + \left\| \frac{x_\varepsilon^k - x_\varepsilon}{\alpha_k} + v \right\| \\
        &\le
        \frac{2}{\alpha_k} \|x^0-x_\varepsilon\| 
        + \left\| \frac{x_\varepsilon^k - x_\varepsilon}{\alpha_k} + v \right\| \\
        &\le
        \frac{2}{\alpha_k} \|x^0-x_\varepsilon\| + \varepsilon
    \end{align*}
    holds, where the third inequality comes from \cref{lem:mannupperbound}.
    If $v\in \cR(\opI-\opT)$, there exists $x_\star\in\cH$ such that $v=x_\star-\opT x_\star$, and the above proof stil holds with $\varepsilon=0$ and $x_\varepsilon = x_\star$.
    Therefore,
    \[
        \left\| \frac{x^k - x^0}{\alpha_k} + v \right\|
        \le \frac{2}{\alpha_k} \|x^0-x_\star\|.
    \]
    If $\lim_{k\to\infty}\alpha_k = \infty$, for any $\varepsilon>0$,
    \[
        \limsup_{k\to\infty} \left\| \frac{x^k-x^0}{\alpha_k} + v\right\| \le \varepsilon,
    \]
    so we get $\lim_{k\to\infty} \frac{x^k-x^0}{\alpha_k} = - v$.
\end{proof}

\begin{remark}
    By obtaining the upper bound to $\alpha_k$, we may optimize the upper bound of \cref{thm:mannaverconv}.
    From the definition of $\{\alpha_k\}_{k\in\mathbb{N}\cup\{0\}}$,
    \begin{align*}
        \alpha_{k} 
        &= 
        (1-\nu^k_0) + \sum_{i=2}^{k} \nu^k_i \alpha_{i-1}
        =
        \sum_{i=1}^k \nu_i^k (1 + \alpha_{i-1})
    \end{align*}
    with $\alpha_0 = 0$.

    Consider an extreme case of \eqref{eq:fpi}, which corresponds to the choice of $\{\nu_i^k\}_{i=1,\ldots,k}$ for $k\in\mathbb{N}\cup\{0\}$ as
    \[
        \nu_i^k 
        =
        \begin{cases}
            0 & (0 \le i \le k-1) \\
            1 & (i = k)
        \end{cases}
    \]
    In this case, $\alpha_k = k$.
    We claim that this is the biggest possible value for $\alpha_k$ for any $k\in\mathbb{N}$, using induction.
    First, $\alpha_0 = 0$.
    Suppose $\alpha_i \le i$ for all $i$ such that $0 \le i < k$.
    Then
    \begin{align*}
        \alpha_k
        &=
        \sum_{i=1}^k \nu_i^k (1+\alpha_{i-1}) \\
        &\le
        \sum_{i=1}^k \nu_i^k \{ 1 + (i-1) \} \\
        &\le
        \sum_{i=1}^k \nu_i^k \{ 1 + (k-1) \} \\
        &=
        k \sum_{i=1}^k \nu_i^k \le k.
    \end{align*}
    Therefore, $\alpha_k \le k$ for all $k\in\mathbb{N}$.
    Hence \eqref{eq:fpi} yields optimal upper bound, which is the same optimal upper bound as in \cref{thm:mannaverconv}.
\end{remark}

\section{Omitted proofs of \cref{sec:pep}}
\label{sec:proofpep}

\subsection{Omitted proofs of \cref{subsec:interpolation}}
\label{sec:proofinterpolation}

Below results will be used to prove \cref{thm:interpolate}.

\begin{theorem}[Projection theorem, Theorem~3.16, \citet{bauschke2011convex}] \label{thm:projection}
    Let $C$ be a nonempty closed convex subset of $\cH$.
    Then for every $x$ and $p$ in $\cH$,
    \[
        p = \Proj_C x
        ~\Leftrightarrow~
        \Big[
            \, \langle y-p,\, x-p \rangle \le 0, \quad \forall\, y\in C \,
        \Big]
    \]
\end{theorem}

\begin{theorem}[Corollary~5, \citet{bauschke2007fenchel}] \label{thm:interpolatefn}
    Let $D$ be a nonempty subset of $\cH$ and let $\opT\colon D\to\cH$ be firmly-nonexpansive operator.
    Then there exists a firmly-nonexpansive operator $\widetilde{\opT}\colon\cH\to\cH$ such that $\widetilde{\opT}\mid_D = \opT$ and $\cR(\widetilde{\opT}) \subset \clconv\cR(\opT)$.
\end{theorem}

\begin{lemma} \label{lem:projconv}
    Let $R$ be a nonempty set in $\cH$.
    Suppose that $v\in R$ is a vector such that
    \[
        \langle x - v,\, v \rangle \ge 0, \quad \forall\, x \in R.
    \]
    Then
    \[
        \langle x - v,\, v \rangle \ge 0, \quad \forall\, x \in \clconv R.
    \]
    % or in other words, $v = P_{\clconv R}(0)$, a minimum norm element in $\clconv R$.
\end{lemma}
\begin{proof}
    Let $x\in\clconv R$.
    Then there exists $\{x_k\}_{k\in\mathbb{N}}$ such that $x_k\in \conv R$ for all $k$ and $\lim_{k\to\infty} x_k = x$.
    
    Since $x_k\in\conv R$, for each $k$, there exist $n_k\in\mathbb{N}$, $\alpha_i^k\in(0,1]$ and $x_i^k\in R$ for $i=1,\ldots,n_k$ such that
    \[
        x_k = \sum_{i=1}^{n_k} \alpha_i^k x_i^k
    \]
    and $\sum_{i=1}^{n_k} \alpha_i^k = 1$.
    \begin{align*}
        \langle x_k,\, v \rangle
        &=
        \left\langle \sum_{i=1}^{n_k} \alpha_i^k x_i^k,\, v \right\rangle
        =
        \sum_{i=1}^{n_k} \alpha_i^k \langle x_i^k,\, v \rangle
        \ge
        \sum_{i=1}^{n_k} \alpha_i^k \|v\|^2
        =
        \|v\|^2.
    \end{align*}
    So $\langle x_k-v,\, v \rangle \ge 0$ for all $x_k \in\conv R$.
    Then
    \[
        \langle x,\, v \rangle
        = \left\langle \lim_{k\to\infty} x_k,\, v \right\rangle
        = \lim_{k\to\infty} \langle x_k,\, v \rangle
        \ge \|v\|^2.
    \]
\end{proof}

\begin{lemma} \label{lem:interpolatewithrange}
    Let $\{(x_i,\, y_i)\}_{i\in I} \subset \cH \times \cH$ be a set of vectors with index set $I$ such that
    \[
        \|y_i - y_j\| \le \|x_i - x_j\|, \quad \forall\, i,j \in I
    \]
    and define $D = \{x_i\}_{i\in I} \subset \cH$.
    Then there exists a nonexpansive operator $\widetilde{\opT}\colon\cH\to\cH$ such that $\widetilde{\opT} \mid_D = \opT$ and
    \[
        \overline{\cR(\opI-\widetilde{\opT})} = \clconv \cR(\opI-\opT).
    \]
\end{lemma}

\begin{proof}
    Define an operator $\opT\colon D\to\cH$ as
    \[
        \opT x_i = y_i, \quad i\in I
    \]
    where $D = \{x_i\}_{i\in I} \subset \cH$.
    Then $\opS\colon D\to\cH$ defined as $\opS = \frac{\opI-\opT}{2}$ is a firmly-nonexpansive operator.
    According to \cref{thm:interpolatefn}, there exists a firmly-nonexpansive extension $\widetilde{\opS}\colon\cH\to\cH$ of $\opS$ such that $\widetilde{\opS} \mid_D = \opS$ and $\cR(\widetilde{\opS}) \subset \clconv \cR(\opS)$.
    If $\widetilde{\opT} = \opI - 2\widetilde{\opS}$, then $\widetilde{\opT}\colon\cH\to\cH$ becomes a nonexpansive extension of $\opT$ such that $\widetilde{\opT} \mid_D = \opT$ and $\cR(\opI - \widetilde{\opT}) = 2 \cR(\widetilde{\opS}) \subset 2 \clconv \cR(\opS) = \clconv \cR(2 \opS) = \clconv \cR(\opI-\opT)$.
    Obviously,
    \[
        \cR(\opI-\opT) \subseteq \cR(\opI-\widetilde{\opT}) \subseteq \clconv\cR(\opI-\opT).
    \]
    Since $\overline{\cR(\opI-\widetilde{\opT})}$ is a convex set,
    \[
        \conv \cR(\opI-\opT) \subseteq \overline{\cR(\opI-\widetilde{\opT})} \subseteq \clconv \cR(\opI-\opT),
    \]
    and as it is also a closed set,
    \[
        \overline{\cR(\opI-\widetilde{\opT})} = \clconv \cR(\opI-\opT).
    \]
\end{proof}

We now prove \cref{thm:interpolate}.

\begin{proof}[Proof of \cref{thm:interpolate}]
    Let $C = \clconv \{x_i - y_i\}_{i\in I}$.
    \begin{itemize}
    \item[($i$)]
        From \cref{lem:interpolatewithrange}, there exists a nonexpansive operator $\widetilde{\opT}\colon\cH\to\cH$ such that $y_i = \widetilde{\opT} x_i,\, \forall\, i\in I$ and $\overline{\cR(\opI-\widetilde{\opT})} = C$.
        Then $v = \Proj_C(0) = \Proj_{\overline{\cR(\opI-\widetilde{\opT})}}(0)$ is an infimal displacement vector of $\widetilde{\opT}$.
        
    \item[($ii$)]
        Further assume that $v = x_\star - y_\star$ with $\star\in I$ and
        \[
            \langle x_i - y_i,\, v \rangle \ge \|v\|^2, \quad \forall\, i \in I.
        \]
        Then \cref{lem:interpolatewithrange} asserts that there exists a nonexpansive operator $\widetilde{\opT}\colon\cH\to\cH$ such that $y_i = \widetilde{\opT} x_i,\, \forall\, i\in I$ and $\overline{\cR(\opI-\widetilde{\opT})} = C$.
        According to \cref{lem:projconv}, $\langle z,\, v \rangle \ge \|v\|^2$ for all $z \in C$.
        Then from \cref{thm:projection}, $v = \Proj_C(0) = \Proj_{\overline{\cR(\opI-\widetilde{\opT})}}(0)$ so it is an infimal displacement vector of $\widetilde{\opT}$.
    \end{itemize}
\end{proof}

\subsection{Omitted proofs of \cref{subsec:pepformulation}}
\label{subsec:proofpepformulation}

Problem we want to solve is a maximization problem in following form, given $k\in\mathbb{N}$ and an index set $I=\{0,1,\ldots,k,\star\}$.
As pointed out, we restrict the choice of nonexpansive operator $\opT$ to be the ones where $v$ actually lies in the range of $\opI-\opT$.
\begin{align*}
    \underset{\opT}{\mbox{maximize}}
    &\quad
    \|x^k - \opT x^k - v\|^2 \\
    \mbox{subject to} 
    &\quad
    \text{$\opT\colon\cH\to\cH$ is nonexpansive} \\
    &\quad
    v = \Pi_{\overline{\cR(\opI-\opT)}}(0) = x_\star - \opT x_\star \\
    &\quad
    x^{n+1} = \frac{n+1}{n+2} \opT x^n + \frac{1}{n+2} x^0, \quad n=0,1,\ldots,k-1 \\
    &\quad
    \|x^0-x_\star\|^2 \le R^2
\end{align*}
Without loss of generality, we may only consider the case of $R=1$, which can be rescaled by $R$ to obtain original problem.
\begin{align*}
    \underset{\opT}{\mbox{maximize}}
    &\quad
    \|x^k - \opT x^k - v\|^2 \\
    \mbox{subject to} 
    &\quad
    \text{$\opT\colon\cH\to\cH$ is nonexpansive} \\
    &\quad
    v = \Pi_{\overline{\cR(\opI-\opT)}}(0) = x_\star - \opT x_\star \\
    &\quad
    x^{n+1} = \frac{n+1}{n+2} \opT x^n + \frac{1}{n+2} x^0, \quad n=0,1,\ldots,k-1 \\
    &\quad
    \|x^0-x_\star\|^2 \le 1
\end{align*}

Above problem is an infinite-dimensional problem, which is possibly an intractable problem.
Such dimensionality stems from the variable of this problem, $\opT$, lying in a function space which cannot be finite-dimensional.

We use \cref{thm:interpolate} reduce the problem dimension by not considering the whole function space of nonexpansive operators any more.
According to \cref{thm:interpolate}, the existence of nonexpansive $\opT$ with $v = x_\star - \opT x_\star$ is equivalent to the existence of iterates $\{(x^i,\, y^i)\}_{i \in I}$ satisfying the following inequalities.
\begin{align*}
    \|y^i - y^j\|^2 \le \|x^i - x^j\|^2,
    &\quad \forall\, i,j\in I,\, i\neq j \\
    \langle x^i - y^i,\, v \rangle \ge \|v\|^2,
    &\quad \forall\, i \in I
\end{align*}

Therefore, the problem can be reformulated as
\begin{align*}
    \underset{\{(x^i,\, y^i)\}_{i\in I}}{\mbox{maximize}}
    &\quad
    \|x^k - y^k - v\|^2 \\
    \mbox{subject to} 
    &\quad
    \left( \exists \cH \right)\,
    \left( x^i,\, y^i \in \cH, \quad \forall\, i \in I \right) \\
    &\quad
    \|y^i - y^j\|^2 \le \|x^i - x^j\|^2, \quad \forall\, i,j \in I,\, i \neq j \\
    &\quad
    v = x_\star - y_\star \\
    &\quad
    \langle x^i - y^i,\, v \rangle \ge \|v\|^2, \quad \forall\, i \in I \\
    &\quad
    x^{n+1} = \frac{n+1}{n+2} y^n + \frac{1}{n+2} x^0, \quad n=0,1,\ldots,k-1 \\
    &\quad
    \|x^0-x_\star\|^2 \le 1
\end{align*}
However, this problem is still intractable in a sense that the iterates $\{(x^i,\, y^i)\}_{i\in I}$ needs to be searched within any choice of real Hilbert space $\cH$.
We remove such dependency using the semidefinite formulation of PEP.

Consider a gram matrix $Z \in \SS^{k+3}$ defined as
\begin{align*}
    Z 
    &=
    \begin{bmatrix}
        \|v^0\|^2 & \langle v^0,\, v^1 \rangle & \cdots & \langle v^0,\, v^k \rangle & \langle v^0,\, v \rangle & \langle v^0,\, x^0 - x_\star \rangle \\
        \langle v^1,\, v^0 \rangle & \|v^1\|^2 & \cdots & \langle v^1,\, v^k \rangle & \langle v^1,\, v \rangle & \langle v^1,\, x^0-x_\star \rangle \\
        \vdots & \vdots & \ddots & \vdots & \vdots & \vdots \\
        \langle v^k,\, v^0 \rangle & \langle v^k,\,v^1 \rangle & \cdots & \|v^k\|^2 & \langle v^k,\, v \rangle & \langle v^k,\, x^0-x_\star \rangle \\
        \langle v,\, v^0 \rangle & \langle v,\, v^1 \rangle & \cdots & \langle v,\, v^k \rangle & \|v\|^2 & \langle v,\, x^0-x_\star \rangle \\
        \langle x^0-x_\star, v^0 \rangle & \langle x^0-x_\star, v^1 \rangle & \cdots & \langle x^0-x_\star,\, v^k \rangle & \langle x^0-x_\star,\, v \rangle & \| x^0-x_\star \|^2 \\
    \end{bmatrix} \\
    &=
    \begin{bmatrix}
        v^0 & v^1 & \cdots & v^k & v & x^0 - x_\star
    \end{bmatrix}^\intercal
    \begin{bmatrix}
        v^0 & v^1 & \cdots & v^k & v & x^0 - x_\star
    \end{bmatrix}
\end{align*}
where $v^i = x^i - y^i$ for $i\in I\setminus\{\star\}$ and $v = x_\star - y_\star$.
Let $G$ denote the horizontal stack of vectors
\[
    G = \begin{bmatrix}
        v^0 & v^1 & \cdots & v^k & v & x_\star & x^0 - x_\star
    \end{bmatrix},
\]
then $Z = G^\intercal G$.
From
\[
    x^{n+1} = \frac{n+1}{n+2} y^n + \frac{1}{n+2} x^0, \quad n=0,1,\ldots,k-1 
\]
being equivalent to
\begin{equation}
    x^{n+1} = x^0 - \sum_{i=0}^{n} \frac{i+1}{n+2}\, v^i , \quad n=0,1,\ldots,k-1 ,
    \label{eq:halpernspan}
\end{equation}
and we use this fact for our semidefinite PEP formulation.

For notational simplicity, let $\ve_i \in \RR^{k+3}$ denote the $i$-th canonical basis vector, i.e., only the $i$-th entry of $(k+3)$-dimensional real vector is $1$ and all the other entries are $0$, and let $a \odot b = \frac{1}{2} (a b^\intercal + b a^\intercal)$.
\begin{itemize}
\item[(i)] Objective function.

\begin{align*}
    \| x^k - y^k - v \|^2
    &=
    \|v^k - v\|^2 \\
    &=
    \left( G (\ve_{k+1} - \ve_{k+2}) \right)^\intercal \left( G (\ve_{k+1} - \ve_{k+2}) \right) \\
    &=
    \tr \left( (\ve_{k+1} - \ve_{k+2})(\ve_{k+1} - \ve_{k+2})^\intercal Z \right)
\end{align*}
Letting $C_k = (\ve_{k+1} - \ve_{k+2})(\ve_{k+1} - \ve_{k+2})^\intercal$, $\|x^k-y^k-v\|^2 = \tr(C_k Z)$.

\item[(ii)] Interpolation condition on nonexpansiveness.

Using \eqref{eq:halpernspan},
\begin{align*}
    &\| x^i - x^j \|^2 - \| y^i - y^j \|^2 \\
    &=
    2 \langle (x^i - y^i) - (x^j - y^j),\, x^i - x^j \rangle - \| (x^i - y^i) - (x^j - y^j) \|^2 \\
    &=
    2 \langle v^i - v^j,\, x^i - x^j \rangle - \| v^i - v^j \|^2 \\
    &=
    2 \left\langle G (\ve_{i+1} - \ve_{j+1}),\, G \left( - \sum_{l=0}^{i-1} \frac{l+1}{i+1} \ve_{l+1} + \sum_{m=0}^{j-1} \frac{m+1}{j+1} \ve_{j+1} \right) \right\rangle
    % \\ &\quad
    - \left(G (\ve_{i+1} - \ve_{j+1}) \right)^\intercal \left(G (\ve_{i+1} - \ve_{j+1}) \right) \\
    &=
    \tr \left[ \left\{
        - 2 (\ve_{i+1} - \ve_{j+1}) \odot \left( \sum_{l=0}^{i-1} \frac{l+1}{i+1} \ve_{l+1} - \sum_{m=0}^{j-1} \frac{m+1}{j+1} \ve_{j+1} \right)
        + (\ve_{i+1} - \ve_{j+1}) (\ve_{i+1} - \ve_{j+1})^\intercal
        \right\}
    Z \right]
\end{align*}
Letting
\[
    A_{i,j}
    = 
    - 2 (\ve_{i+1} - \ve_{j+1}) \odot \left( \sum_{l=0}^{i-1} \frac{l+1}{i+1} \ve_{l+1} - \sum_{m=0}^{j-1} \frac{m+1}{j+1} \ve_{j+1} \right)
    + (\ve_{i+1} - \ve_{j+1}) (\ve_{i+1} - \ve_{j+1})^\intercal,
\]
the inequality condition $\|y^i-y^j\|^2 \le \|x^i-x^j\|^2$ for $i,j\in I\setminus\{\star\}$ is equivalent to $\tr(A_{i,j} Z) \ge 0$.

Also,
\begin{align*}
    &\|x^i-x_\star\|^2 - \|y^i-y_\star\|^2 \\
    &=
    2 \langle (x^i - y^i) - (x_\star - y_\star),\, x^i - y^i \rangle - \| (x^i - y^i) - (x_\star - y_\star) \|^2 \\
    &=
    2 \langle v^i - v,\, x^i - x_\star \rangle - \|v^i - v\|^2 \\
    &=
    2 \left\langle
        G(\ve_{i+1} - \ve_{k+2}),\,
        G \left( \ve_{k+3} - \sum_{l=0}^{i-1} \frac{l+1}{i+1} \ve_{l+1} \right)
    \right\rangle
    - \left( G(\ve_{i+1} - \ve_{k+2}) \right)^\intercal \left( G(\ve_{i+1} - \ve_{k+2}) \right) \\
    &=
    \tr \left[ \left\{
        - 2 (\ve_{i+1} - \ve_{k+2}) \odot \left( \sum_{l=0}^{i-1} \frac{l+1}{i+1} \ve_{l+1} - \ve_{k+3} \right)
        - (\ve_{i+1} - \ve_{k+2}) (\ve_{i+1} - \ve_{k+2})^\intercal
        \right\} Z
    \right].
\end{align*}
Letting 
\[
    A_{i,\star} = 
    - 2 (\ve_{i+1} - \ve_{k+2}) \odot \left( \sum_{l=0}^{i-1} \frac{l+1}{i+1} \ve_{l+1} - \ve_{k+3} \right)
    - (\ve_{i+1} - \ve_{k+2}) (\ve_{i+1} - \ve_{k+2})^\intercal,
\]
the inequality condition $\|y^i - y_\star\|^2 \le \|x^i - x_\star\|^2$ is equivalent to $\tr(A_{i,\star} Z) \ge 0$.

\item[(iii)] Interpolation condition on infimal displacement vector.

\begin{align*}
    \langle v^i,\, v \rangle - \|v\|^2 
    &=
    \langle v^i - v,\, v \rangle \\
    &=
    \langle G(\ve_{i+1} - \ve_{k+2}),\, G \ve_{k+2} \rangle \\
    &=
    \tr \left[ \left( (\ve_{i+1} - \ve_{k+2} ) \odot \ve_{k+2} \right) Z \right]
\end{align*}
Therefore, letting $B_i = (\ve_{i+1} - \ve_{k+2})\odot\ve_{k+2}$, the inequality condition $\langle v^i,\, v \rangle \ge \|v\|^2$ is equivalent to $\tr(B_i Z) \ge 0$.

\item[(iv)] Initial point condition.

\begin{align*}
    \|x^0-x_\star\|^2
    &=
    \left( G \ve_{k+3} \right)^\intercal \left( G \ve_{k+3} \right) \\
    &=
    \tr \left( \ve_{k+3}{\ve_{k+3}}^\intercal Z \right)
\end{align*}
so if $D_0 = \ve_{k+3}{\ve_{k+3}}^\intercal$, the inequality condition $\|x^0-x_\star\|^2 \le 1$ is equivalent to $\tr(D_0 Z) \le 1$.
\end{itemize}

Gathering all these facts, the problem at hand can be reformulated into the semidefinite program
\[
    \begin{array}{ll}
        \underset{Z\in\SS^{k+3}_+}{\mbox{maximize}} & \tr(C_k Z) \\
        \mbox{subject to} & \tr(A_{i,j} Z) \ge 0, \quad \forall\, i,j \in I\setminus\{\star\},\, i\neq j \\
        & \tr(A_{i,\star} Z) \ge 0, \quad \forall\, i\in I \setminus\{\star\} \\
        & \tr(B_{i} Z) \le 0, \quad \forall\, i\in I\setminus\{\star\} \\
        & \tr(D_0 Z) \le 1
    \end{array}
\]
Here, the condition on which real Hilbert space $\cH$ and that the iterates $x^i$'s and $y^i$'s must be defined can be ignored, and this problem indeed can be solved with numerical solvers.
The equivalence of the last reformulation comes from \cref{lem:interpolatedim}.

\begin{lemma} \label{lem:interpolatedim}
    If $\dim\cH \ge k+3$, $Z$ is a positive-semidefinite $(k+3)\times(k+3)$ matrix if and only if
    there exist $x^0-x_\star$, $v$, and $v^i = x^i-y^i$ for $i=0,1,\ldots,k$ in $\cH$ such that $G$ is defined as in \eqref{eq:g-def} and $Z=G^\intercal G$.
\end{lemma}

\subsection{Numerical result of PEP}

We numerically solved the SDP formulated in \cref{subsec:pepformulation} to obtain a numerical guarantee on the rate of convergence to $\|x^k - \opT x^k - v\|^2$ for \eqref{eq:halpern} with $\lambda_k = \frac{1}{k+1}$.
We used MOSEK \cite{aps2019mosek} with $k = 1,2,\dots,100$.
We observe that the numerical solution of PEP to indicate an optimal rate of $\tilde{\mathcal{O}}(1/k^2)$ but not $\mathcal{O}(1/k^2)$.

\begin{figure*}[ht]
    \centering
    \includegraphics[width=.48\textwidth]{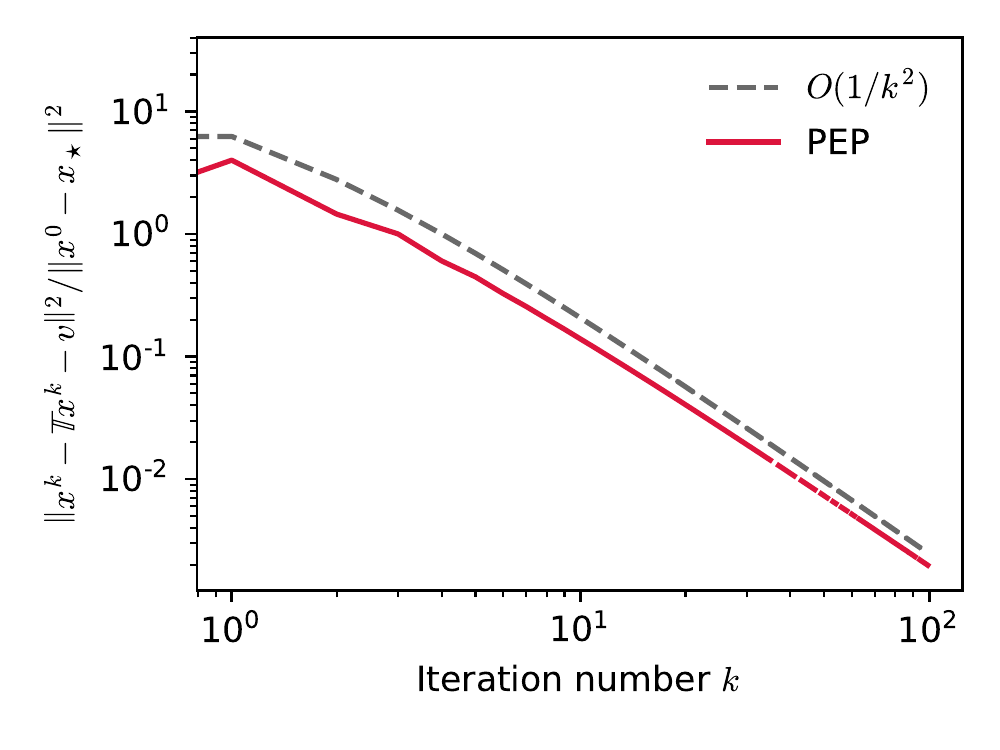}
    \hfill
    \includegraphics[width=.48\textwidth]{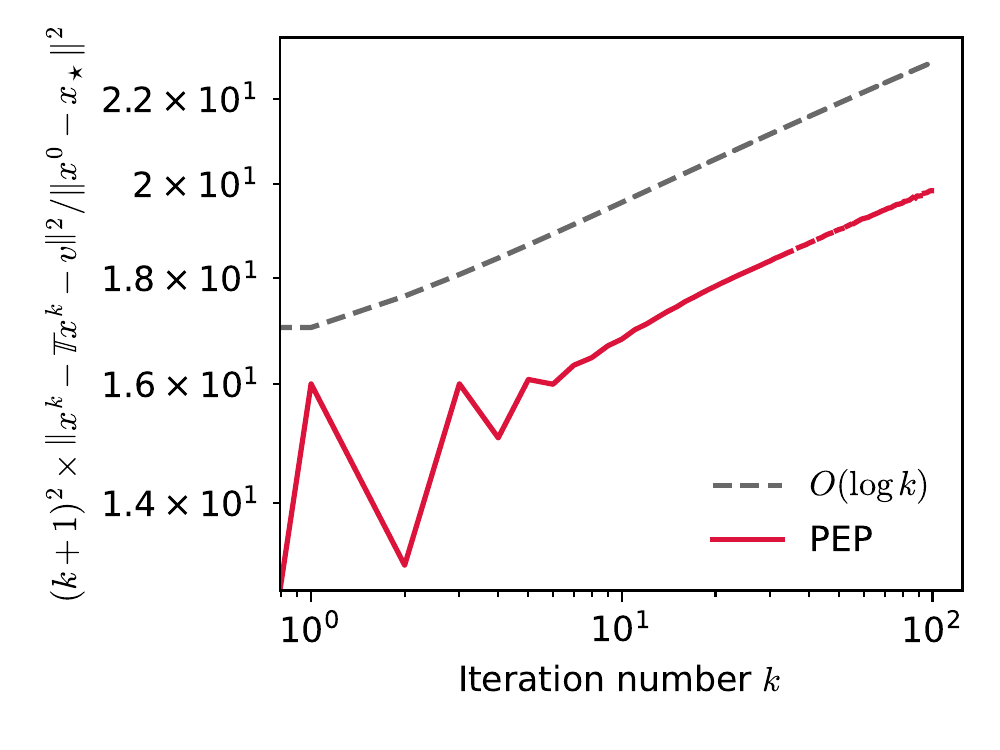}
    \caption{
        We solved the problem for iteration numbers $k=1$ through $k=100$.
        (Left) Plot of $\|x^k - \opT x^k - v\|^2 / \|x^0 - x_\star\|^2$ for $k=1,\dots,100$.
        (Right) Plot of $(k+1)^2 \cdot \|x^k - \opT x^k - v\|^2 / \|x^0 - x_\star\|^2$ for $k=1,\dots,100$.
        }
    \label{fig:pep_experiment}
\end{figure*}

\section{Omitted proofs of \cref{sec:lowerbound}}
\label{sec:prooflowerbound}

We use following lemma in the proof of \cref{thm:lowerboundspan} and \cref{thm:lowerboundgeneral}.
\begin{lemma} \label{lem:zeroresp}
    Consider any orthogonal matrix $U\colon\RR^m\to\RR^n$ such that $m\le n$ and $U^\intercal U = I_m$.
    For any nonexpansive operator $\opT\colon\RR^m\to\RR^m$ and any $x_0\in\RR^n$, define $\opT_U\colon\RR^n\to\RR^n$ as $\opT_U(\cdot) = U\opT U^\intercal (\cdot-x_0) + x_0$.
    Then,
    \begin{itemize}
        \item[$(i)$] $\|Ux\| = \|x\|$ for any $x\in\RR^m$ and $\|U^\intercal x\| \le \|x\|$ for any $x\in\RR^n$.
        \item[$(ii)$] $\opT_U\colon\RR^n\to\RR^n$ is a nonexpansive operator.
        \item[$(iii)$] $U^\intercal \cR(\opI-\opT_U) = U^\intercal U\cR(\opI-\opT) = \cR(\opI-\opT)$
        \item[$(iv)$] $\tilde{v} = \Proj_{\overline{\cR(\opI-\opT)}}(0)$ if and only if $v = U \tilde{v} = \Proj_{\overline{\cR(\opI-\opT_U)}}(0)$.
        If there exists $x_\star\in\cH_1$ such that $\tilde{v} = x_\star-\opT x_\star$, then $y_\star = x_0 + Ux_\star$ implies $v = y_\star - \opT_U y_\star$.
        If there exists $y_\star\in\cH_2$ such that $v = y_\star - \opT_U y_\star$, then $x_\star = U^\intercal (y_\star-x_0)$ implies $x_\star-\opT x_\star = \tilde{v}$.
        This implies $\fix\opT=\emptyset$ if and only if $\fix\opT_U=\emptyset$ as well.
    \end{itemize}
\end{lemma}

\begin{proof}
    From orthogonality of $U$, $U^\intercal U = I_m$ and $UU^\intercal$ is an orthogonal projection onto the range of $U$.
    \[
        \|U x\|^2
        = \langle Ux, Ux \rangle
        = \langle x, U^\intercal U x \rangle
        = \langle x, x \rangle
        = \|x\|^2,
        \qquad \forall\, x \in \RR^m
    \]
    and
    \[
        \|U^\intercal x\|^2
        = \|UU^\intercal x\|^2
        \le \|x\|^2,
        \qquad \forall\, x \in \RR^n.
    \]
    Also,
    \begin{align*}
        \|\opT_U y - \opT_U z\|
        &=
        \|U (\opT U^\intercal (y-x_0) - \opT U^\intercal (z-x_0))\| \\
        &=
        \|\opT U^\intercal (y-x_0) - \opT U^\intercal (z-x_0)\| \\
        &\le
        \| U^\intercal (y-x_0) - U^\intercal (z-x_0) \| \\
        &\le
        \| y-z \|, \quad \forall\, y,z \in \RR^n
    \end{align*}
    so $\opT_U$ is a nonexpansive operator.

    Finally,
    \[
        U(x-\opT x) = Ux - U\opT x = Ux - U\opT U^\intercal Ux = (\opI-\opT_U)(Ux+x_0),
        \quad \forall\, x\in\RR^m,
    \]
    so $U \cR(\opI-\opT) \subseteq \cR(\opI-\opT_U)$, and $U^\intercal U \cR(\opI-\opT) = \cR(\opI-\opT) \subseteq U^\intercal \cR(\opI-\opT_U)$.
    \[
        U^\intercal (y-\opT_Uy)
        =
        U^\intercal y - U^\intercal U \opT U^\intercal y
        =
        (\opI-\opT) (U^\intercal y+x_0),
        \quad \forall\, y\in\RR^n,
    \]
    so $U^\intercal \cR(\opI-\opT_U) \subseteq \cR(\opI-\opT)$ and $(iii)$ holds true.

    In order to prove $(iv)$, note that from $(i)$, $\|\tilde{v}\|^2 = \|U\tilde{v}\|^2$.
    Suppose $\tilde{v} = \Proj_{\overline{\cR(\opI-\opT_U)}}(0)$.
    Then from $(iii)$, $U^\intercal v \in \overline{\cR(\opI-\opT)}$.
    $\tilde{v}$ is the minimum norm element of $\overline{\cR(\opI-\opT)}$, so
    \[
        \|U\tilde{v}\|
        = \|\tilde{v}\|
        \le \|U^\intercal v\|
        \le \|v\|.
    \]
    $U \tilde{v} \in \overline{\cR(\opI-\opT_U)}$ and its norm is smaller than or equal to that of ${v}$. Therefore, $U \tilde{v}$ must also be the minimum norm element of $\overline{\cR(\opI-\opT_U)}$, which leads to $U \tilde{v} = v$ by the uniqueness of such element.
    
    If $x_\star\in\RR^m$ is a point where $\tilde{v} = x_\star-\opT x_\star$, then setting $y_\star = x_0 + Ux_\star$ leads to
    \[
        v
        = U \tilde{v} 
        = Ux_\star-U\opT x_\star 
        = (Ux_\star + x_0) - \opT_U (Ux_\star + x_0) 
        = y_\star - \opT_U y_\star.
    \]

    Now, if $y_\star\in\RR^n$ is a point where $v = y_\star-\opT_U y_\star$, then setting $x_\star = U^\intercal (y_\star-x_0)$ leads to
    \begin{align*}
        x_\star - \opT x_\star
        &= 
        U^\intercal (y_\star - x_0) - U^\intercal U \opT U^\intercal (y_\star - x_0) \\
        &= U^\intercal (y_\star - \opT_U y_\star)
        = U^\intercal \tilde{v}
        = U^\intercal U \tilde{v}
        = \tilde{v}.
    \end{align*}    
\end{proof}

\subsection{Fixed-Point Iteration with Span Assumption}

\begin{proof}[Proof of \cref{thm:lowerboundspan}]
    Let $\ve_i\in\RR^{k+1}$ denote an $i$-th canonical basis vector whose $i$-th entry is $1$ and the other entries are all zero.
    As we discussed in the outline of proof, we only consider the case $\tilde{v}=(0,\dots,0,\|v\|) = \|v\| \ve_{k+1}$.
    Define $\opT\colon\RR^{k+1}\to\RR^{k+1}$ as
    \[
        x - \opT x = 
        \underbrace{\begin{bmatrix}
            1 & 0 & 0 & \hdots & 0 & 1 & 0 \\
            -1 & 1 & 0 & \hdots & 0 & 0 & 0 \\
            0 & -1 & 1 & \hdots & 0 & 0 & 0 \\
            \vdots & \vdots & \vdots & \ddots & \vdots & \vdots & \vdots \\
            0 & 0 & 0 & \hdots & 1 & 0 & 0 \\
            0 & 0 & 0 & \hdots & -1 & 1 & 0 \\
            0 & 0 & 0 & \hdots & 0 & 0 & 0
        \end{bmatrix}}_{= M \in \mathbb{R}^{(k+1)\times (k+1)}} x
        + \alpha \ve_1 - \|v\| \ve_{k+1}, \quad \forall x \in \cH_1
    \]
    for $\alpha\neq0$.
    Note that $M$ is an invertible matrix.
    If $v=0$, $\opT$ has fixed-points of the form
    \[
        - \frac{\alpha}{2} \sum_{i=1}^k \ve_i + (x_\star)_{k+1} \ve_{k+1}, \quad (x_\star)_{k+1} \in \RR.
    \]
    If $v\neq 0$, $\opT$ does not have a fixed point, and its infimal displacement vector $\tilde{v} = \|v\| \ve_{k+1} \neq 0$.

    Let the iterates $\{x^n\}_{n=0}^k$ satisfy the linear span assumption \eqref{def:span}.
    Since $x^0 = 0$, $x^0-\opT x^0\in\Span\{\ve_1, \tilde{v}\}$ and 
    \[
        x^1 \in x^0 + \Span\{x^0-\opT x^0\} \subseteq \Span\{\ve_1, \tilde{v} \}.
    \]
    Then $x^1-\opT x^1 \in \Span\{\ve_1, \ve_2, \tilde{v}\}$ and we also have
    \[
        x^2
        =
        x^0 + \Span\{x^0-\opT x^0, x^1-\opT x^1\}
        \in \Span\{\ve_1,\ve_2,\tilde{v} \}.
    \]
    
    From the observation above, we claim that
    \begin{align*}
        x^n &\in \Span\{ \ve_1,\ve_2,\ldots,\ve_n,\tilde{v}  \} \\
        x^n - \opT x^n &\in \alpha \ve_1 + \tilde{v}  + \Span\{ M\ve_1,M\ve_2,\ldots,M\ve_{n}\}
        \subseteq
        \Span\{ \ve_1,\ve_2,\ldots,\ve_{n+1}, \tilde{v} \}
    \end{align*}
    for $n=1,\dots,k-1$.
    
    We have already proven the case of $n=1$.
    Let $n < k-1$ and assume that the claim above hold for all $m$ such that $m \le n$.
    Then
    \begin{align*}
        x^{n+1}
        &\in
        x^0 + \Span\{x^0-\opT x^0, \ldots, x^n-\opT x^n\}
        \subseteq
        \Span\{\ve_1,\ldots,\ve_{n+1},\tilde{v} \}.
    \end{align*}
    Also,
    \begin{align*}
        x^{n+1} - \opT x^{n+1}
        &=
        M x^{n+1} + \alpha \ve_1 + \tilde{v} \\
        &\in
        \alpha \ve_1 + \tilde{v} + \Span\{M\ve_1,\ldots,M\ve_{n+1}\} \\
        &\subseteq
        \Span\{ \ve_1,\ve_2,\ldots,\ve_{n+2}, \tilde{v} \}.
    \end{align*}
    The claim holds for $n+1$ as well.

    From above claim and its proof, we have
    \[
        \sum_{i=0}^{k-1} \nu_i (x^i-\opT x^i) - \tilde{v} 
        \in
        \alpha \ve_1 + \underbrace{\Span\{ M\ve_1,\ldots,M\ve_{k-1}\}}_{=: R_{k-1}},
    \]
    so
    \begin{align*}
        \left\| \sum_{i=0}^{k-1} \nu_i (x^i-\opT x^i) - \tilde{v} \right\|^2
        &\ge
        \|\alpha \ve_1\|^2 - \|\Proj_{R_{k-1}} (\alpha \ve_1)\|^2 \\
        &=
        \alpha^2 \left\| \Proj_{R_{k-1}^\perp} (\ve_1) \right\|^2.
    \end{align*}
    Since $R_{k-1}^\perp = \Span \left\{ \sum_{i=1}^{k} \ve_i,\, \tilde{v} \right\}$,
    \begin{align*}
        \left\| \sum_{i=0}^{k-1} \nu_i (x^i-\opT x^i) - \tilde{v} \right\|^2
        &\ge
        \alpha^2 \left\| \Proj_{\Span\{ \sum_{i=1}^{k} \ve_i,\, \tilde{v} \}}(\ve_1) \right\|^2 \\
        &=
        \alpha^2 \left\| \Proj_{\Span\{ \sum_{i=1}^{k} \ve_i \}}(\ve_1) \right\|^2 \\
        &=
        \alpha^2 \left\|
            {\left\langle \ve_1, \frac{\sum_{i=1}^{k} \ve_i}{\|\sum_{i=1}^{k} \ve_i \|} \right\rangle} \frac{\sum_{i=1}^{k} \ve_i}{\|\sum_{i=1}^{k} \ve_i \|}
        \right\|^2
        =
        \frac{\alpha^2}{k}
    \end{align*}
    
    We know that the set of possible choice of $x_\star\in\cH$ is
    \[
        \left\{
            x_\star \in \RR^{k+1}
            \,\mid\,
            x_\star =
            - \frac{\alpha}{2}\sum_{i=1}^{k} \ve_i + (x_\star)_{k+1} \ve_{k+1},
            \, \forall (x_\star)_{k+1} \in \RR
        \right\}.
    \]
    As $x^0 = 0$, $\|x^0-x_\star\|^2 \ge \frac{k \alpha^2}{4}$ and equality holds when $(x_\star)_{k+1} = 0$.

    Gathering the facts above, we may conclude that there exists $\opT\colon\RR^{k+1}\to\RR^{k+1}$ with infimal displacement vector $\tilde{v}$ and corresponding $x_\star\in\RR^{k+1}$ such that $x_\star - \opT x_\star = \tilde{v}$ and
    \[
        \left\| \sum_{i=0}^{k-1} \nu_i (x^i-\opT x^i) - \tilde{v} \right\|^2
        \ge
        \frac{4}{k^2} \|x^0-x_\star\|^2
    \]
    for any iterates $\{x^n\}_{n=0}^{k-1}$ satisfying the linear span assumption, starting from $x^0=0$.

    Now, we may use the same operator $\opT$ to prove that
    \[
        \left( \left\| \sum_{i=0}^{k-1} \nu_i (x^i-\opT x^i) \right\| - \|\tilde{v}\| \right)^2
        \ge
        \frac{1}{2k^2} \|x^0-x_\star\|^2
    \]
    for any iterates $\{x^n\}_{n=0}^k$ satisfying the linear span assumption, starting from $x^0=0$.

    If $v = 0$, then $\tilde{v} = 0$ so
    \begin{align*}
        \left\| \sum_{i=0}^{k-1} \nu_i (x^i-\opT x^i) \right\| - \|\tilde{v}\|
        &= 
        \left\| \sum_{i=0}^{k-1} \nu_i (x^i-\opT x^i) - \tilde{v} \right\| \\
        &\ge
        \frac{2}{k} \|x^0-x_\star\|^2 \\
        &\ge
        \frac{1}{\sqrt{2}k} \|x^0-x_\star\|^2
    \end{align*}
    so the desired inequality
    \[
        \left( \left\| \sum_{i=0}^{k-1} \nu_i (x^i-\opT x^i) \right\| - \|\tilde{v}\| \right)^2
        \ge
        \frac{1}{2k^2} \|x^0-x_\star\|^2
    \]
    holds true.
    
    Now suppose $v\neq0$.
    Following the calculations above,
    \begin{align*}
        \left\| \sum_{i=0}^k \nu_i (x^i-\opT x^i) \right\| - \|\tilde{v}\|
        &=
        \sqrt{\|\tilde{v}\|^2 + \|\alpha\ve_1\|^2 - \|\Proj_{R_{k-1}}(\alpha\ve_1)\|^2} - \|\tilde{v}\| \\
        &=
        \sqrt{\|\tilde{v}\|^2 + \alpha^2 \left\| \Proj_{R_{k-1}^\perp}(\ve_1) \right\|^2} - \|\tilde{v}\|.
    \end{align*}
    With a real function $f(t) = \sqrt{t}$ defined on $[0,\infty)$, we may use the fact that $f$ is a concave function, therefore
    \[
        f(t+h) - f(t)
        \ge h f'(t+h)
        = \frac{h}{2\sqrt{t+h}},
        \quad \forall\, t,h>0.
    \]
    Substituting $t$ and $h$ by $\|v\|^2>0$ and $\alpha^2 \|\Proj_{R_{k-1}^\perp}(\ve_1)\|^2>0$,
    \begin{align*}
        \left\| \sum_{i=0}^{k-1} \nu_i (x^i-\opT x^i) \right\| - \|\tilde{v}\|
        &\ge
        \frac{\alpha^2 \|\Proj_{R_{k-1}^\perp}(\ve_1)\|^2}{2\sqrt{\|\tilde{v}\|^2 + \alpha^2 \|\Proj_{R_{k-1}^\perp}(\ve_1)\|^2}} \\
        &=
        \frac{\alpha\|\Proj_{R_{k-1}^\perp}(\ve_1)\|}{2\sqrt{1 + \frac{\|\tilde{v}\|^2}{\alpha^2 \|\Proj_{R_{k-1}^\perp}(\ve_1)\|^2}}}.
    \end{align*}
    $\alpha>0$ is a positive real number which was left unspecified, so we may just assign any positive value for the calculation.
    Let $\alpha = \frac{\|\tilde{v}\|}{\|\Proj_{R_{k-1}^\perp}(\ve_1)\|}=\sqrt{k}\|\tilde{v}\|>0$, then as $\|x^0-x_\star\| = \frac{\alpha}{2} \sqrt{k} = \frac{k\|\tilde{v}\|}{2}$, we get
    \begin{align*}
        \left(
            \left\| \sum_{i=0}^{k-1} \nu_i (x^i-\opT x^i) \right\| - \|\tilde{v}\|
        \right) / \|x^0-x_\star\|
        &\ge
        \frac{\|\tilde{v}\|}{2 \sqrt{2}} \frac{2}{k \|\tilde{v}\|}
        =
        \frac{1}{\sqrt{2}{k}}.
    \end{align*}
    so
    \[
        \left( \left\| \sum_{i=0}^{k-1} \nu_i (x^i-\opT x^i) \right\| - \|\tilde{v}\| \right)^2
        \ge
        \frac{1}{2k^2} \|x^0-x_\star\|^2.
    \]

    We now extend the result to the arbitrarily given $v$, not $\tilde{v}$.
    The case $v=0$ is trivial, so suppose $v\neq0$.
    Let $U\in\RR^{(k+1)\times(k+1)}$ be an orthogonal matrix such that $U^\intercal U = I_{k+1}$ and $U \tilde{v} = v$.
    Since $U$ is a square matrix, $U^\intercal U = U U^\intercal = I_{k+1}$ so $\|U^\intercal x\| = \|x\|$ for all $x\in\RR^{k+1}$.
    According to \cref{lem:zeroresp}, $\opT_U$ defined as $\opT_U(\cdot) = U\opT U^\intercal (\cdot - y^0) + y^0$ is a nonexpansive operator with $v = U \tilde{v} = \Pi_{\overline{\cR(\opI-\opT_U)}}(0)$.
    Let $\{y^n\}_{n=0}^k$ is a sequence of iterates satisfying the linear span assumption \eqref{def:span} with $\opT_U$.
    \begin{align*}
        U^\intercal (y^{n} - \opT_U y^n)
        &=
        U^\intercal (y^n - y^0) - U^\intercal U \opT U^\intercal (y^n - y^0) \\
        &=
        U^\intercal (y^n - y^0) - \opT U^\intercal (y^n - y^0)
    \end{align*}
    and as
    \[
        y^{n+1} \in y^0 + \Span\{ y^0 - \opT_U y^0,\, \dots,\, y^n - \opT_U y^n\}
    \]
    implies
    \begin{align*}
        U^\intercal (y^{n+1} - y^0)
        &\in 
        \Span\{ U^\intercal y^0 - U^\intercal \opT_U y^0,\, \dots,\, U^\intercal y^n - U^\intercal \opT_U y^n\} \\
        &=
        \Span\{ U^\intercal (y^0-y^0) - \opT U^\intercal (y^0-y^0),\, \dots, U^\intercal (y^n-y^0) - \opT U^\intercal (y^n-y^0) \},
    \end{align*}
    $\left\{ U^\intercal (y^n-y^0) \right\}_{n=0}^{k-1}$ satisfies linear span assumption \eqref{def:span} with $\opT$ where $\tilde{v} = U^\intercal {v}$ is an infimal displacement vector of $\opT$.
    Therefore,
    \begin{align*}
        \left\| \sum_{i=0}^{k-1} \nu_i (y^i-\opT_U y^i) - v \right\|
        &\ge
        \left\| U^\intercal \sum_{i=0}^{k-1} \nu_i (y^i - \opT_U y^i) - U^\intercal v \right\| \\
        &=
        \left\| \sum_{i=0}^{k-1} \nu_i \left\{ U^\intercal (y^n-y^0) - \opT U^\intercal (y^n-y^0) \right\} - \tilde{v} \right\| 
        \tag{$\because$ $U^\intercal v = U^\intercal U \tilde{v} = \tilde{v}$} \\
        &\ge
        \frac{4}{k^2} \|U^\intercal (y^0-y^0) - U^\intercal (y_\star - y^0) \|^2 \\
        &=
        \frac{4}{k^2} \|y^0 - y_\star\|^2.
    \end{align*}
    Similarly,
    \begin{align*}
        \left\| \sum_{i=0}^{k-1} \nu_i (y^i-\opT_U y^i) \right\| - \|v\|
        &\ge
        \left\| U^\intercal \sum_{i=0}^{k-1} \nu_i (y^i - \opT_U y^i) \right\| - \|v\| \\
        &=
        \left\| \sum_{i=0}^{k-1} \nu_i \left\{ U^\intercal (y^n-y^0) - \opT U^\intercal (y^n-y^0) \right\} \right\| - \|\tilde{v}\|
        \tag{$\because$ $v = U\tilde{v}$} \\
        &\ge
        \frac{1}{\sqrt{2}k} \|U^\intercal (y^0-y^0) - U^\intercal (y_\star - y^0) \| \\
        &=
        \frac{1}{\sqrt{2}k} \|y^0 - y_\star\|.
    \end{align*}

    Proof with real Hilbert space $\cH$ can be done in the same manner as the proof of $\cH=\RR^{k+1}$.
    Set $\{\ve_i\}_{i=1}^{k+1}$ as a set of orthonormal basis vectors of $\cH$ where $\ve_{k+1} = \frac{v}{\|v\|}$ if $v\neq0$ and arbitrary if $v=0$, and proceed with the same proof as above.
\end{proof}

\begin{remark}
    We calculate the complexity lower bound of various quantities including fixed-point residual $x^k-\opT x^k$ and normalized iterate $(x^k-x^0)/\alpha_k$ by appropriately choosing the convex coefficients $\nu_i$'s.
    In order to measure the lower bound of fixed-point residual $x^{k-1} - \opT x^{k-1}$ converging to $v$, choose $\nu_{k-1} = 1$ and choose all other $\nu_i$'s as $0$.
    For normalized iterate, we use the fact that KM and Halpern choose the iterate $x^{k}$ to be in form of $x^0 + \sum_{i=0}^{k-1} \lambda_i^k (x^i-\opT x^i)$.
    Therefore, if we choose $\nu_i = \frac{\lambda_i^k}{\sum_{i=0}^{k-1} \lambda_i^k}$,
    we obtain normalized iterates $\frac{x^{k}-x^0}{\sum_{i=1}^{k} (1-\lambda_i)}$ for KM and $\frac{x^{k}-x^0}{\theta_{k}}$ for Halpern.
    This scheme can be extended to calculate the lower bound of Mann iteration as well.
\end{remark}

\subsection{General Fixed-Point Iterations}

We follow the general complexity lower bound result of \citet{ParkRyu2022} for operators \emph{with} fixed points, and extend their result to the case where fixed point might not exist.

\begin{definition}[Section~D.2,\& D.4, \citet{ParkRyu2022}]
    Let $\cH$ be a real Hilbert space and $\opT\colon\cH\to\cH$ be a nonexpansive operator.
    Let $\{\ve_i\}_{i\in I}$ with index set $I$ denote a set of orthonormal basis of $\cH$.

    A \emph{deterministic fixed-point iteration} $\al$ is defined as a mapping of the point $x^0\in\cH$ and an operator $\opT$ to sequences of iterates $\{x^k\}_{k\in\mathbb{N}}$ and $\{\bar{x}^k\}_{k\in\mathbb{N}}$.
    Here, $x^k$ is the $k$-th \emph{query point} and $\bar{x}^k$ is the $k$-th \emph{approximation point}, and we consider the setup with $x^k = \bar{x}^k$ so we omit $\bar{x}^k$.
    Actually, $\al$ is defined as a sequence of mappings $\{\al_k\}_{k\in\mathbb{N}}$, where $x^k$ is an output of $\al_k$ given the point $x^0$ and the operator $\opT$ where
    \[
        x^k 
        =
        \al_k [x^0, \opT]
        =
        \al_k [x^0, \cO_\opT(x^0), \cO_\opT(x^1), \dots, \cO_\opT(x^{k-1})]
    \]
    for any $k\in\mathbb{N}$.
    Here, $\{x^k\}_{k\in\mathbb{N}}$ only depends on $x^0$ and $\opT$ via the \emph{fixed-point residual oracle} $\cO_\opT (x) = x - \opT x$.
    $\al$ is deterministic in a sense that, when provided with the same point $x^0$ and the same oracle queries $\cO_\opT (x^k) = x^k - \opT x^k$ for $k\in\mathbb{N}$, $\al$ will give the same sequence of iterates $\{x^k\}$ as an output of $\al$.

    For $z\in\cH$, denote by $\supp\{z\}$ the \emph{support} of $z$, i.e., 
    \[
        \supp\{z\} = \left\{
            i \in I \mid \langle z,\, \ve_i \rangle \neq 0
        \right\}.
    \]
    We say a sequence $\{z^t\}_{t\in\mathbb{N}\cup\{0\}}$ is \emph{zero-respecting with respect to $\opT$} if
    \[
        \supp \{z^t\} \subseteq \bigcup_{s < t} \, \supp \{ z^s - \opT z^s \}
    \]
\end{definition}

If $\{z^t\}_{t\in\mathbb{N}\cup\{0\}}$ is a zero-respecting sequence with respect to $\opT$, then by the definition, $\supp\{z^0\} \subseteq \emptyset$ so $z^0 = 0$.

As we did for fixed-point iterations with linear span assumption \eqref{def:span}, from now on, we prove the result for Euclidean spaces, since the proof naturally extends to Hilbert spaces $\cH$ with $\dim\cH\ge 2k-1$ and its set of orthonormal vectors $\{\ve_i\}_{i=1}^{2k-1}$.

\begin{lemma} \label{lem:zeroresprate}
    Given $k\in\mathbb{N}$ and $v\in\RR^{k+1}$, let $\opT\colon\RR^{k+1}\to\RR^{k+1}$ be the worst-case operator defined in the proof of \cref{thm:lowerboundspan}, along with its infimal displacement vector $v$ and $x_\star$.
    If the iterates $\{x^n\}_{n=0}^{k-1}$ are zero-respecting with respect to $\opT$,
    \[  
        \left\| \sum_{i=0}^{k-1} \nu_i (x^i-\opT x^i) - v \right\| \ge \frac{2}{k+1} \|x^0-x_\star\|
    \]
    and
    \[
        \left\| \sum_{i=0}^{k-1} \nu_i (x^i-\opT x^i) \right\| - \|v\| \ge \frac{1}{\sqrt{2}(k+1)} \|x^0-x_\star\|
    \]
    for any choice of real numbers $\{\nu_i\}_{i=0}^{k-1}$ such that $\sum_{i=0}^{k-1} \nu_i = 1$.
\end{lemma}

\begin{proof}
    We claim that any zero respecting sequence $\{x^n\}_{n=0}^{k-1}$ satisfies
    \begin{align*}
        x^n
        &\in \Span\{\ve_1,\dots,\ve_n\} \\
        x^n-\opT x^n
        &\in \Span\{\ve_1,\dots,\ve_{n+1}\}
    \end{align*}
    for $n=0,\dots,k$.
    If this holds, then the proof of \cref{thm:lowerboundspan} is still applicable to zero-respecting sequences, leading to the desired result.
    
    If $n=0$, $x^0 = 0$ and $x^0-\opT x^0\in\Span\{\ve_1\}$, so the case of $n=0$ holds true.
    Now, suppose that $n<k-1$ and the claim above holds for all $m$ such that $m\le n$.
    Since the iterates form a zero-respecting sequence with respect to $\opT$, $\supp \{x^{n+1}\} \subseteq \cup_{m\le n} \supp \{x^m-\opT x^m\}$ and therefore $x^{n+1}\in\Span\{\ve_1,\dots,\ve_{n+1}\}$.
    Using this fact, $x^{n+1}-\opT x^{n+1} \in \Span\{\ve_1,\dots,\ve_{n+2}\}$ easily follows.
\end{proof}

\begin{lemma} \label{lem:existU}
    Let $\al$ be a general deterministic fixed-point iteration and $\opT\colon\RR^{n+1}\to\RR^{n+1}$ be a nonexpansive operator defined as in the proof of \cref{thm:lowerboundspan}.
    For any arbitrary $x^0 \in \RR^d$ and $v \in \RR^d$ with $d\ge n + N$, there exists an orthogonal matrix $U\in\RR^{d\times (n+1)}$ and the iterates $\{x^k\}_{k=1}^{N} = \al[x^0, \opT_U]$ such that $x^{(k)} = U^\intercal (x^k - x^0)$, $\left\{ x^{(k)} \right\}_{k=0}^N$ is zero-respecting with respect to $\opT$, and $v$ becomes an infimal displacement vector of $\opT_U$.
\end{lemma}

\begin{proof}
    We prove the existence of an orthogonal matrix $U\in\reals^{d\times n}$ such that $\left\{ x^k \right\}^N_{k=1} = \al[x^0,\, \opT_U]$, $x^{(k)} = U^\intercal (x^k-x^0)$ and $v$ becomes an infimal displacement vector of $\opT_U$.
    Constructing such orthogonal matrix is equivalent to choosing appropriate set of orthonormal vectors $\{u_i\}_{i=1}^{n+1}$, whose $i$-th vector $u_i$ becomes an $i$-th column of matrix $U$, i.e.,
    \[
        U = 
        \begin{bmatrix}
            \vert & \dots & \vert \\
            u_1 & \dots & u_n \\
            \vert & \dots & \vert
        \end{bmatrix}.
    \]
    We modify the proof of \citep[Lemma~D.4(i)]{ParkRyu2022} to cover the case when $v\neq0$, as the original proof is restricted to the case where the fixed point exists, or in other words, the case where $v=0$.
    
    We provide the scheme which inductively finds the column $u_i$'s, given an arbitrary nonzero vector $v\in\reals^d$.
    Define the set of indices $S_t$ for $t\in\{ 1,\dots,N \}$ as
    \[
        S_t = \bigcup_{s<t} \supp \left\{ x^{(s)} - \opT x^{(s)} \right\} \subseteq \{1,2,\dots,n+1\}
    \]
    and note that $S_0 = \emptyset \subseteq S_1 \subseteq \dots \subseteq S_t$.
    As $v\neq 0$, $x^{(0)} - \opT x^{(0)} \neq 0$, and $S_1 \neq \emptyset$.
    Choose a set of vectors $\{u_i\}_{i\in S_1}$ to be any unit vectors which are orthonormal to each other.
    The precise choice of $\{u_i\}_{i\in S_1}$ will be later specified, and it will make $v$ to be an infimal displacement vector of $\opT_U$.
    
    Now, suppose that for $t \ge 2$, $\{u_i\}_{i\in S_{t-1}}$ is already chosen.
    Choose a set of unit vectors $\{u_i\}_{i \in S_t \setminus S_{t-1}}$ from the orthogonal complement of
    \[
        W_t := \Span \left( \{ x^1 - x^0,\, \cdots,\, x^{t-1} - x^0 \} \cup \{ u_i \}_{i\in S_{t-1}} \right)
        \subseteq \reals^d
    \]
    and let them be orthogonal to each other.
    When $t=N$ and $S_N \neq \{1,\dots,n+1\}$, choose any $\{u_i\}_{i \in \{1,\dots,n+1\}\setminus S_N}$ which makes $U$ orthogonal.
    Above scheme is well-defined if 
    \[
        \dim W_t^\perp \ge |S_t \setminus S_{t-1}|,
    \]
    and this is guaranteed from the fact that $d \ge n+N$ and $t \le N$ since
    \[
        \dim W_t^\perp = d - \dim W_t
        \ge d - \{(t-1) + |S_{t-1}|\}
        \ge
        (n+1) - |S_{t-1}|
        \ge
        |S_t \setminus S_{t-1}|.
    \]
    Since $\langle u_i,\, y_t - y_0 \rangle$ for $i \notin S_t$ for $t=1,\dots,N$,
    \[
        x^{(t)} = U^\intercal (x^t - x^0) \in \Span\{e_i\}_{i\in S_t}
    \]
    leads to $\supp \{x^{(t)}\} \subseteq S_t$.
    This proves that there exists an orthogonal matrix $U\in\reals^{d\times (n+1)}$ such that $\{x^k\}_{k=1}^N = \al[x^0, \opT_U]$ and $x^{(k)} = U^\intercal (x^k-x^0)$ for all $k=1,\dots,N$.

    Now, it remains to show that certain choice of $\{u_i\}_{i\in S_1}$ implies that $\opT_U$ has $v$ as its infimal displacement vector.
    First, observe that for any arbitrary choice of $\{u_i\}_{i\in S_1}$,
    \[
        S_1 = \supp \left\{ x^{(0)} - \opT x^{(0)} \right\}
    \]
    and
    \[
        x^{(0)} - \opT x^{(0)}
        =
        0 - \opT 0 
        =
        - \alpha \ve_1 + \|v\| \ve_n
    \]
    so $S_1 = \{1, n\}$.
    Note that the infimal displacement vector of $\opT$ is $\tilde{v} = \|v\| \ve_{n+1}$.
    From \cref{lem:zeroresp}, $U\tilde{v} = \|v\| \cdot U \ve_{n+1} = \|v\| u_n$ is an infimal displacement vector of $\opT_U$.
    As $n+1\in S_1$, we may choose $u_{n+1} = \frac{v}{\|v\|}$ so that $v$ is an infimal displacement vector of $\opT_U$.
\end{proof}

\begin{lemma}
    Consider the setup of \cref{lem:existU} with $U\in\RR^{m\times n}$ and the iterates $\{x^k\}_{k=1}^{N} = \al[x^0, \opT_U]$.
    Then
    \[
        \left\| \sum_{i=0}^k \nu_i (x^{(i)}-\opT x^{(i)}) - \tilde{v} \right\|
        \le
        \left\| \sum_{i=0}^k \nu_i (x^i-\opT_U x^i) - {v} \right\|
    \]
    and
    \[
        \left\| \sum_{i=0}^k \nu_i (x^{(i)}-\opT x^{(i)}) \right\| - \|\tilde{v}\|
        \le
        \left\| \sum_{i=0}^k \nu_i (x^i-\opT_U x^i) \right\| - \|{v}\|
    \]
    for any $\nu_i\in\reals$ such that $\sum_{i=0}^k \nu_i = 1$, where $v$ is an infimal displacement vector of $\opT_U$ and $\tilde{v}$ is an infimal displacement vector of $\opT$.
\end{lemma}

\begin{proof}
    According to \cref{lem:zeroresp}, $v = U\tilde{v}$.

    \begin{align*}
        \left\|
            \sum_{i=0}^k \nu_i (x^{(i)}-\opT x^{(i)}) - \tilde{v}
        \right\|
        &=
        \left\|
            \sum_{i=0}^k \nu_i \left( U^\intercal (x^i-x^0) - \opT U^\intercal (x^i-x^0) \right) - U^\intercal v
        \right\| \\
        &=
        \left\|
            U^\intercal \sum_{i=0}^k \nu_i \left( (x^i-x^0) - U \opT U^\intercal (x^i-x^0) \right) - U^\intercal v
        \right\| \\
        &=
        \left\|
            U^\intercal \sum_{i=0}^k \nu_i (x^i - \opT_U x^i) - U^\intercal v
        \right\| \\
        &\le
        \left\|
            \sum_{i=0}^k \nu_i (x^i-\opT_U x^i) - v
        \right\|.
    \end{align*}
    Note that $\|\tilde{v}\| = \|Uv\| = \|v\|$.
    Therefore,
    \begin{align*}
        \left\|
            \sum_{i=0}^k \nu_i (x^{(i)}-\opT x^{(i)})
        \right\|
        - \|\tilde{v}\|
        &=
        \left\|
            U^\intercal \sum_{i=0}^k \nu_i (x^i - \opT_U x^i)
        \right\|
        - \|{v}\| \\
        &\le
        \left\|
            \sum_{i=0}^k \nu_i (x^i-\opT_U x^i)
        \right\| - \|{v}\|.
    \end{align*}
\end{proof}

We now prove the main result.
\begin{theorem} \label{thm:lowerboundgeneralR}
    Let $d \ge 2k-1$ for $k\in\mathbb{N}$.
    For any deterministic fixed-point iteration $\al$, any initial point $x^0\in\RR^d$ and any vector ${v}\in\RR^d$, there exists a nonexpansive operator $\opT\colon\RR^d\to\RR^d$ such that
    \[
        \left\| \sum_{i=0}^{k-1} \nu_i (x^i-\opT x^i) - {v} \right\|^2
        \ge
        \frac{4}{k^2} \|x^0-x_\star\|^2
    \]
    and
    \[
        \left( \left\| \sum_{i=0}^{k-1} \nu_i (x^i-\opT x^i) \right\| - \|{v}\| \right)^2
        \ge
        \frac{1}{2k^2} \|x^0-x_\star\|^2
    \]
    where $v = x_\star-\opT x_\star = \Proj_{\overline{\cR(\opI-\opT)}}(0)$ and $\nu_i\in\reals$ with $\sum_{i=0}^{k-1} \nu_i = 1$.
\end{theorem}

\begin{proof}
    Let $\opS\colon\RR^{k+1}\to\RR^{k+1}$ be a worst-case operator in the proof of \cref{thm:lowerboundspan}.
    From \cref{lem:existU}, there exists an orthogonal matrix $U\in\RR^{d\times (k+1)}$ such that $d \ge k + (k-1) = 2k-1$,  $\left\{ x^{(t)} \right\}_{t=0}^{k-1}$ a sequence of iterates defined as $x^{(t)} = U^\intercal (x^t-x^0)$ is zero-respecting with respect to $\opS$,
    \[
        \left\| \sum_{i=0}^{k-1} \nu_i (x^{(i)}-\opS x^{(i)}) - \tilde{v} \right\|
        \le
        \left\| \sum_{i=0}^{k-1} \nu_i (x^i-\opS_U x^i) - {v} \right\|,
    \]
    \[
        \left\| \sum_{i=0}^{k-1} \nu_i (x^{(i)}-\opS x^{(i)}) \right\| - \|\tilde{v}\|
        \le
        \left\| \sum_{i=0}^{k-1} \nu_i (x^i-\opS_U x^i) \right\| - \|{v}\|
    \]
    for any $\nu_i\in\reals$ such that $\sum_{i=0}^{k-1} \nu_i = 1$, and $v = U\tilde{v}$ or $\tilde{v} = U^\intercal v$ by \cref{lem:zeroresp}.
    
    % Meanwhile, if you look into the proof of \cref{lem:existU}, for specific operator $\opT$ as in our assumption, we are able to choose $U$ such that $U \tilde{v} = v$, therefore the given $v$ becomes the infimal displacement vector of $\opT_U$.
    % This is possible by, following the representation in the proof of Lemma~D.4(i) of \citet{ParkRyu2022}, $S_1 = \{0,k+1\}$ and we may choose the $(k+1)$-th column $u_{k+1} = \frac{v}{\|v\|}$ in case of $v\neq0$.
    
    Since $\left\{ x^{(t)} \right\}_{t=0}^{k-1}$ is zero-respecting with respect to $\opT$, \cref{lem:zeroresprate} implies
    \[
        \left\| \sum_{i=0}^{k-1} \nu_i (x^{(i)}-\opS x^{(i)}) - \tilde{v} \right\|^2
        \ge
        \frac{4}{k^2} \|x^{(0)}-x^{(\star)}\|^2
    \]
    and
    \[
        \left( \left\| \sum_{i=0}^{k-1} \nu_i (x^{(i)}-\opS x^{(i)}) \right\| - \|\tilde{v}\| \right)^2
        \ge
        \frac{1}{2k^2} \|x^{(0)}-x^{(\star)}\|^2
    \]
    where $x^{(\star)}\in\cH_0$ is a point such that $x^{(\star)}-\opS x^{(\star)}=v$.
    If $x_\star = x^0 + Ux^{(\star)}$,
    \[
        \|x^{(0)}-x^{(\star)}\|^2
        =
        \| -x^{(\star)} \|^2
        =
        \| - Ux^{(\star)} \|^2
        =
        \|x^0-x_\star\|^2
    \]
    so we may conclude that
    \[
        \left\| \sum_{i=0}^{k-1} \nu_i (x^i-\opS_U x^i) - {v} \right\|^2
        \ge 
        \frac{4}{k^2} \|x^0-x_\star\|^2
    \]
    and
    \[
        \left( \left\| \sum_{i=0}^{k-1} \nu_i (x^i - \opS_U x^i) \right\| - \|v\| \right)^2 \ge \frac{1}{2k^2} \|x^0 - x_\star\|^2.
    \]
    Therefore, $\opT = \opS_U$ is our desired worst-case operator.
\end{proof}

\begin{proof}[Proof of \cref{thm:lowerboundgeneral}]
    Use the worst-case nonexpansive operator of \cref{thm:lowerboundgeneralR} and construct the nonexpansive operator with the orthonormal basis set $\{\ve_i\}_{i=1}^{2k-1}$ of $\cH$ with $\dim\cH = 2k-1$.
\end{proof}

\section{Details of experiment in \cref{sec:experiments}}
\label{sec:proofexperiments}

Consider a semidefinite problem (SDP)
\[
    \begin{array}{ll}
        \underset{x\in\RR^d}{\mbox{minimize}} &
        \sum_{i=1}^p c_i^\intercal x \\
        \mbox{subject to} &
        \mathcal{A}_i[x] = \sum_{j=1}^{d} A_i^j x_j \preceq B_i, 
        \quad 1\le i\le p,
    \end{array}
\]
and PG-EXTRA applied on this problem.
\begin{align*}
     \!\!U_i^{k+1} &\!= \Pi_{-\SS^n_+} \left( U_i^k + \beta(B_i - \mathcal{A}_i[x^k]) \right) \\
     \!\!\vw^{k+1} &\!= \vw^k + \frac{1}{2}(I-W)\vx^k
    \tag{PG-EXTRA}
    \\
     \!\!x_i^{k+1} &\!= x_i^k - \alpha\beta (2 w_i^{k+1} - w_i^k)+\alpha \left( 
    \mathcal{A}_i^*
    [ 2 U_i^{k+1} - U_i^k] - c_i \right)
\end{align*}

\subsection{Deriving PG-EXTRA for SDP}

Consider
\[
    \begin{array}{ll}
        \underset{x\in\reals^m}{\mbox{minimize}} & \sum_{i=1}^p \langle c_i, x \rangle_{\reals^m} \\
        \mbox{subject to} & \sum_{j=1}^m x_j A_i^j \preceq_{S^n_+} B_i, \quad i=1,\cdots,p
    \end{array}
\]
or in other words,
\[
    \begin{array}{ll}
        \underset{x_i\in\reals^n}{\mbox{minimize}} & \sum_{i=1}^p \langle c_i, x_i \rangle_{\reals^m} \\
        \mbox{subject to} & L_i(x_i) := B_i - \sum_{j=1}^m (x_i)_j A_i^j \succeq_{S^n_+} 0, \quad i=1,\cdots,p \\
        & (I-W)\vx = 0 ~\Leftrightarrow~ U\vx = 0
    \end{array}.
\]
Defining $L\,\vx = (L_1 (x_1) - B_1, \cdots, L_p (x_p) - B_p, U\vx)$, which is a linear map from $\reals^{p\times m}$ to $(S^n)^p \times \reals^{p\times m}$,
\[
    \begin{array}{ll}
        \underset{x_i\in\reals^n}{\mbox{minimize}} & \sum_{i=1}^p \langle c_i, x_i \rangle_{\reals^m} + \delta_{(S^n_+)^p\times\{0\}}( L\,\vx + \mathbf{B})
    \end{array}
\]
Corresponding Lagrangian is
\[
    \lagrange(\vx, \vu)
    =
    \underbrace{\sum_{i=1}^p \langle c_i,\, x_i \rangle_{\reals^m}}_{:= \langle \vc,\, \vx \rangle_{\reals^{p\times m}}}
    + \langle \vu, L\,\vx + \mathbf{B} \rangle_{(S^n)^p \times \reals^{p\times m}}
    - \left(\delta_{(S^n_+)^p\times\{0\}}\right)^* (\vu)
\]
where $\vu = (u_1,\cdots,u_p,\vy) \in (S^n)^p \times \reals^{p\times m}$
and its saddle subdifferential is
\[
    \partial\lagrange(\vx,\vu)
    =
    \begin{bmatrix}
        \vc + L^*\vu \\
        - L\,\vx - \mathbf{B} + \partial \left(\delta_{(S^n_+)^p\times\{0\}}\right)^* (\vu)
    \end{bmatrix}
\]
where
\[
    \vc =
    % \begin{bmatrix}
    %     \hline & c_1^\intercal & \hline \\
    %     & \hdots & \\
    %     \hline & c_p^\intercal & \hline 
    % \end{bmatrix}
    \begin{bmatrix}
        c_1^\intercal \\
        \vdots \\
        c_p^\intercal
    \end{bmatrix}
    \in \reals^{p\times m}.
\]
Note that from
\[
    L\,\vx = \left( - \sum_{j=1}^m (x_1)_j A_1^j,\, \cdots,\, - \sum_{j=1}^m (x_p)_j A_p^j,\, U\vx \right),
\]
$L^*\colon (S^n)^p\times \reals^{p\times m} \to \reals^{p\times m}$ is
\[
    L^*(y_1,\cdots,y_p,\,\vz)
    = U^\intercal \vz - \begin{bmatrix} (A_1^* y_1)^\intercal \\ \vdots \\ (A_p^* y_p)^\intercal \end{bmatrix}.
\]
Then
\begin{align*}
    \partial\lagrange(\vx,\vu)
    &=
    \begin{bmatrix}
        \vc - U\vy + \begin{bmatrix} (A_1^* u_1)^\intercal \\ \vdots \\ (A_p^* u_p)^\intercal \end{bmatrix} \\
        - L\,\vx - \mathbf{B} + \partial \left(\delta_{(S^n_+)^p\times\{0\}}\right)^* (\vu)
    \end{bmatrix} \\
    &=
    \underbrace{\begin{bmatrix} \vc \\ - \mathbf{B} \end{bmatrix}}_{=:\opH(\vx,\vu)}
    + \underbrace{\begin{bmatrix} 0 & L^* \\ -L & \partial \left(\delta_{(S^n_+)^p\times\{0\}}\right)^* \end{bmatrix} \begin{bmatrix} \vx \\ \vu \end{bmatrix}}_{=:\opF(\vx,\vu)}.
\end{align*}
Note that
\[
    \delta_{\{0\}}^*(x)
    = \sup_{y} \left( \langle x, y \rangle - \delta_{\{0\}}(y) \right)
    = 0.
\]
and
\[
    \delta_{S^n_+}^*(X)
    = \sup_{Y} \left( \langle X, Y \rangle - \delta_{S^n_+}(Y) \right)
    = \sup_{Y\in S^n_+} \langle X, Y \rangle
    =
    \begin{cases}
        0 & -X \in S^n_{+} \\
        \infty & \mathrm{o.w.}
    \end{cases}
    = \delta_{-S^n_+}(X).
\]

Let
\[
    M = 
    \begin{bmatrix}
        (1/\alpha)\opI & L^* \\
        L & (1/\beta)\opI
    \end{bmatrix},
\]
then the FPI of forward-backward splitting $(\vx^{k+1},\, \vu^{k+1}) = (M + \opF)^{-1} (M - \opH)(\vx^k,\,\vu^k)$ is
\begin{align*}
    &\begin{bmatrix}
        (1/\alpha)\opI & 2L^* \\
        0 & (1/\beta)\opI + \partial \left(\delta_{(S^n_+)^p\times\{0\}}\right)^*
    \end{bmatrix} \begin{bmatrix} \vx^{k+1} \\ \vu^{k+1} \end{bmatrix}
    \ni
    \begin{bmatrix}
        (1/\alpha)\vx^k + L^*\vu^k - \vc \\
        L\vx^k + (1/\beta)\vu^k + \mathbf{B}
    \end{bmatrix} \\
    &\Leftrightarrow
    \vx^{k+1} + 2\alpha L^* \vu^{k+1} = \vx^k + \alpha (L^* \vu^k - \vc) \\
    &\qquad\quad
    \vu^{k+1} = \prox_{\beta(\delta_{(S^n_+)^p\times\{0\}})^*} \left( \vu^k + \beta(L\vx^k + \mathbf{B}) \right) \\
    &\Leftrightarrow
    \vu^{k+1} = \prox_{\beta(\delta_{(S^n_+)^p\times\{0\}})^*} \left( \vu^k + \beta(L\vx^k + \mathbf{B}) \right) \\
    &\qquad\quad
    \vx^{k+1} = \vx^k - \alpha (2L^*\vu^{k+1} - L^*\vu^k + \vc) \\
    &\Leftrightarrow
    u_i^{k+1} = \Pi_{-S^n_+} \left( u_i^k + \beta(B_i - {\scriptstyle\sum_{j=1}^n} (x_i^k)_j A_i^j ) \right) \\
    &\qquad\quad
    \vy^{k+1} = \vy^k + \beta U\vx^k \\
    &\qquad\quad
    \vx^{k+1} = \vx^k - \alpha (2L^*\vu^{k+1} - L^*\vu^k + \vc)
\end{align*}
Note that $y_i$ and $\vy$ are not related to each other.
If we let $\vw^0 = 0$, $\vx^0$ arbitrary, and
\[
    \vw^k = \frac{1}{\beta} U\vy^k = \frac{1}{2} (I-W)\sum_{j=0}^{k-1} \vx^j,
\]
then
\begin{align*}
    & u_i^{k+1} = \Pi_{-S^n_+} \left( u_i^k + \beta(B_i - {\scriptstyle\sum_{j=1}^n} (x_i^k)_j A_i^j ) \right) \\
    & \vx^{k+1} = \vx^k - \alpha \beta (2 \vw^{k+1} - \vw^k)
    + \alpha \left(
        2\begin{bmatrix} (A_1^* u_1^{k+1})^\intercal \\ \vdots \\ (A_p^* u_p^{k+1})^\intercal \end{bmatrix}
        - \begin{bmatrix} (A_1^* u_1^{k})^\intercal \\ \vdots \\ (A_p^* u_p^{k})^\intercal \end{bmatrix}
        - \begin{bmatrix} c_1^\intercal \\ \vdots \\ c_p^\intercal \end{bmatrix}
    \right) \\
    & \vw^{k+1} = \vw^k + \frac{1}{2}(I-W)\vx^k
\end{align*}
or
\begin{align*}
    & u_i^{k+1} = \Pi_{-S^n_+} \left( u_i^k + \beta(B_i - {\scriptstyle\sum_{j=1}^n} (x_i^k)_j A_i^j ) \right) \\
    & \vw^{k+1} = \vw^k + \frac{1}{2}(I-W)\vx^k \\
    & x_i^{k+1} = x_i^k - \alpha\beta (2 w_i^{k+1} - w_i^k) 
    + \alpha \left( A_i^* ( 2 u_i^{k+1} - u_i^k) - c_i \right)
\end{align*}
where $U$ and $u_i$ are irrelevant and
\[
    A^* u =
    \begin{bmatrix}
        \tr(A_1 u_1) \\ 
        \vdots \\
        \tr(A_n u_n)
    \end{bmatrix}.
\]

Above solves decentralized semidefinite problem, when $\alpha,\beta>0$ are chosen to define a metric on $\reals^{n\times p} \times (S^n)^p \times \reals^{n\times p}$.
\[
    M =
    \begin{bmatrix}
        (1/\alpha)\opI & L^* \\
        L & (1/\beta)\opI
    \end{bmatrix}
\]

\subsection{Measuring fixed-point residual in $M$-norm}
Although the algorithm itself does not keep the iterate $\vy^k$ such that
\[
    \vw^k = \frac{1}{\beta} U\vy^k,
\]
we need $\vy^k$-iterates in order to calculate the fixed-point residual $\|(\vx^k,\vu^k)\|^2_{M}$ where $M\colon\reals^{p\times m} \times ((S^n)^p \times \reals^{p\times m}) \to \reals^{p\times m} \times ((S^n)^p \times \reals^{p\times m})$ is a linear map defined as
\begin{align*}
    M = \begin{bmatrix}
        (1/\alpha)\opI & L^* \\
        L & (1/\beta)\opI
    \end{bmatrix}.
\end{align*}
Then for any $\vx\in\reals^{p\times m}$ and $\vu =(u_1,\cdots,u_p,\vy) \in(S^n)^p\times\reals^{p\times m}$,
\begin{align*}
    \| (\vx, \vu) \|_M^2
    &=
    \frac{1}{\alpha} \|\vx\|^2_{\reals^{p\times m}}
    + \frac{1}{\beta} \|\vu\|^2_{(S^n)^p \times {\reals^{p\times m}}}
    + \langle \vx,\, L^* \vu \rangle_{\reals^{p\times m}}
    + \langle L \vx,\, \vu \rangle_{(S^n)^p \times \reals^{p\times m}} \\
    &=
    \frac{1}{\alpha} \|\vx\|^2_{\reals^{p\times m}}
    + \frac{1}{\beta} \sum_{i=1}^p \|u_i\|^2_{S^n}
    + \frac{1}{\beta} \|\vy\|^2_{\reals^{p\times m}}
    + 2 \underbrace{\langle L \vx,\, \vu \rangle_{(S^n)^p \times \reals^{p\times m}}}_{(\star)}.
\end{align*}
Then
\begin{align*}
    (\star)
    &=
    \left\langle
        \left( -\sum_{j=1}^m (x_1)_j A_1^j,\, -\sum_{j=1}^m (x_2)_j A_2^j, \cdots, - \sum_{j=1}^m (x_p)^j A_p^j,\, U\vx \right),
        (u_1,u_2,\cdots,u_p,\vy)
    \right\rangle \\
    &=
    \langle \vx,\, \beta\vw \rangle
    - \sum_{k=1}^p \sum_{j=1}^m (x_k)_j \tr (A_k^j u_k),
\end{align*}
so
\begin{align*}
    \| (\vx, \vu) \|_M^2
    &=
    \frac{1}{\alpha} \|\vx\|^2_{\reals^{p\times m}}
    + \frac{1}{\beta} \sum_{i=1}^p \|u_i\|^2_{S^n}
    + \frac{1}{\beta} \|\vy\|^2_{\reals^{p\times m}}
    + 2 \beta \langle \vx,\, \vw \rangle
    - 2 \sum_{k=1}^p \sum_{j=1}^m (x_k)_j \tr (A_k^j u_k).
\end{align*}

Now, $\vy$ can be calculated as follows.
Consider a eigenvalue decomposition $(I-W) = V\Sigma V^\intercal$.
Let $v_i$ be the $i$-th column of $V$, $\sigma_i$ be the $i$-th eigenvalue corresponding to $v_i$.
Suppose $\sigma_p = 0$ with $v_p=\ones$.
As $\vy \perp \ones$, $\vy = \sum_{i=1}^{p-1} y_i v_i$.
Then
\[
    \beta \vw = U \vy = U \sum_{i=1}^{p-1} y_i v_i.
\]
As $U = V\Sigma^{1/2}V^\intercal$, 
\[
    \beta\vw = V\Sigma^{1/2}V^\intercal \sum_{i=1}^{p=1} y_i v_i
    = V \Sigma^{1/2} \begin{bmatrix}
        y_1 \\ y_2 \\ \vdots \\ y_{p-1} \\ y_p = 0
    \end{bmatrix}
    = V \begin{bmatrix}
        \sqrt{\sigma_1}y_1 \\ \sqrt{\sigma_2}y_2 \\ \vdots \\ \sqrt{\sigma_{p-1}}y_{p-1} \\ 0
    \end{bmatrix}
    = \sum_{i=1}^{p-1} \sqrt{\sigma_i} y_i v_i
\]
Calculate $y_i$ from taking inner product of $\beta\vw$ and $v_i$.

\subsection{Experiment settings and additional plots}

In this experiment, we use the parameters $\alpha=\beta=0.01$ with $n=10$, $m=11$, $p=10$, and $\varepsilon=0.5$.
These numbers come from the infeasible linear matrix inequality (LMI) designed for this experiment, which we state below.
\begin{align*}
    \begin{bmatrix}
        x_1 & x_2 \\
        x_2 & \varepsilon
    \end{bmatrix}
    &\succeq 0 \\
    \begin{bmatrix}
        x_2 & x_3 \\
        x_3 & \varepsilon
    \end{bmatrix}
    &\succeq 0 \\
    &\vdots \\
    \begin{bmatrix}
        x_k & x_{k+1} \\
        x_{k+1} & \varepsilon
    \end{bmatrix}
    &\succeq 0,
\end{align*}
with $\varepsilon>0$.
Then the set of inequalities above is a subset of
\[
    \left\{
        (x_1,\,x_2,\,\dots,\,x_k,\,x_{k+1}) \in \RR^{k+1}
        \,\mid\,
        \frac{x_1}{\varepsilon} \ge \left( \frac{x_{k+1}}{\varepsilon} \right)^{2^k}
    \right\}.
\]
If we add another LMI
\[
    \begin{bmatrix}
        - x_1 & x_2 \\
        x_2 & \varepsilon
    \end{bmatrix}
    \succeq 0,
\]
The feasible region is also a subset of 
\[
    \left\{
        (x_1,\,x_2,\,\dots,\,x_k,\,x_{k+1}) \in \RR^{k+1}
        \,\mid\,
        \frac{x_1}{\varepsilon} \le - \left( \frac{x_{k+1}}{\varepsilon} \right)^{2^k}
    \right\}.
\]
Reversing the sign of the first $(1,1)$-entry of each LMI results in the only feasible region $\{(0,\,\dots,\,0)\}$.
Then, if we additionally impose an LMI such as
\[
    \begin{bmatrix}
        x_1 & 0 \\
        0 & x_{k+1}
    \end{bmatrix}
    \succeq
    \begin{bmatrix}
        1 & 0 \\
        0 & 1
    \end{bmatrix},
\]
Then the origin $\{(0,\,\dots,\,0)\}$ is never in a feasible region of the set of all LMIs, so the SDP becomes infeasible.
The value of $\|v\|^2$ has been numerically calculated using the normalized iterate of Picard iteration after $200,000$ iterations.

Additionally, we draw plots of the difference of fixed-point residual or normalized iterate between $v$ and $-v$, respectively.
\begin{figure}[h]
    \centering
    \includegraphics[width=.45\textwidth]{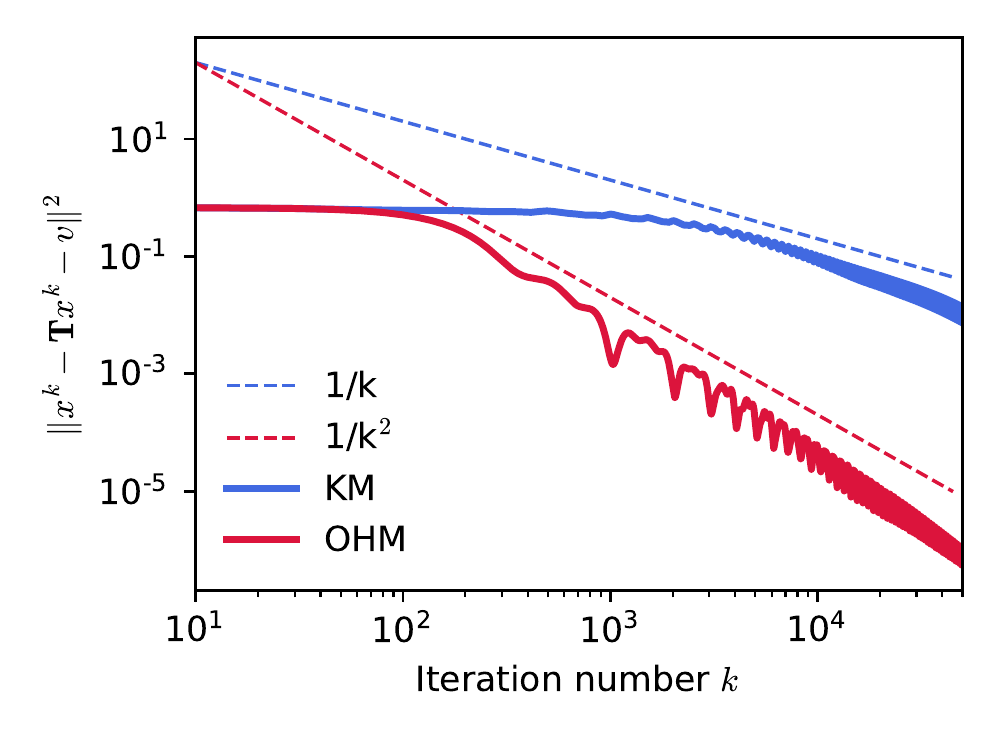}
    \hfill
    \includegraphics[width=.45\textwidth]{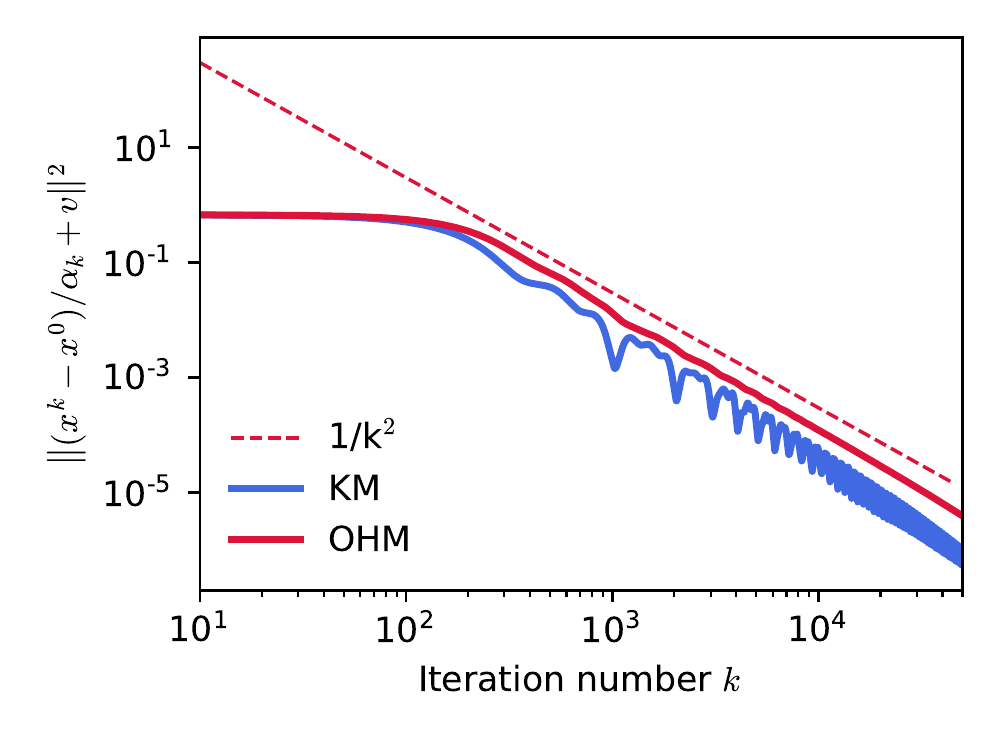}
    \caption{
        (Left) Squared difference between fixed-point residual and $v$, $\|x^k-\opT x^k-v\|^2$.
        (Right) Squared difference between normalized iterate and $-v$, $\|(x^k-x^0)/\alpha_k + v\|^2$.
    }
    \label{fig:my_label}
\end{figure}

%%%%%%%%%%%%%%%%%%%%%%%%%%%%%%%%%%%%%%%%%%%%%%%%%%%%%%%%%%%%%%%%%%%%%%%%%%%%%%%
%%%%%%%%%%%%%%%%%%%%%%%%%%%%%%%%%%%%%%%%%%%%%%%%%%%%%%%%%%%%%%%%%%%%%%%%%%%%%%%

\end{document}